\documentclass[12pt]{amsart}
\usepackage{amsfonts,amsthm,amssymb}

\advance\textheight\topskip
\usepackage[OT2,T1]{fontenc}
\newcommand{\mifody}{%
  \renewcommand\rmdefault{wncyr}%
  \renewcommand\sfdefault{wncyss}%
  \renewcommand\encodingdefault{OT2}%
  \normalfont
  \selectfont}

\usepackage[mathscr]{eucal}
\usepackage{mathabx}
\usepackage{times}
\usepackage{mathptmx}
\newcommand{\myboldsymbol}{\pmb}
\allowdisplaybreaks
\usepackage{t-angles}

\usepackage[matrix,arrow,curve,dvips]{xy} 
\UsePSspecials{dvips}
\usepackage[dvips, bookmarks=false, hyperfigures=false,
setpagesize=false]{hyperref}

\DeclareMathAlphabet{\mathpzc}{OT1}{pzc}{m}{it}

\allowdisplaybreaks

\numberwithin{equation}{section}
\numberwithin{figure}{section}
\makeatletter
\@addtoreset{equation}{section}
\@addtoreset{figure}{section}
\@addtoreset{subsubsection}{section}

\def\@secnumfont{\bfseries}
\def\subsubsection{\@startsection{subsubsection}{3}%
  \z@{.5\linespacing\@plus.7\linespacing}{-.5em}%
  {\normalfont\bfseries}}
\def\paragraph{\@startsection{paragraph}{4}%
  \z@\z@{-\fontdimen2\font}%
  \normalfont\bfseries}
\def\subparagraph{\@startsection{subparagraph}{5}%
  \z@\z@{-\fontdimen2\font}%
  \normalfont\bfseries}
\makeatother

\newcommand{\maps}[1]{\rho_{\symm}[#1]}
\newcommand{\mapa}[1]{\rho_{\asymm}[#1]}
\newcommand{\mapsinv}[1]{\rho^{-1}_{\symm}[#1]}

\newcommand{\ssigma}{\widehat{\sigma}}

\newcommand{\bbigcap}{\mbox{\small${}\bigcap{}$}}

\newcommand{\pone}{(p,1)}
\newcommand{\Wpone}{W_{p,1}}

\newcommand{\W}{\mathscr{W}}
\newcommand{\jW}{\mathsf{W}}
\newcommand{\WW}{\textup{\textsf{\textbf{W}}}}
\newcommand{\WWtritri}{\WW(2,(2p)^{\times 3\times 3})}
\newcommand{\bWWtritri}{\WW\myboldsymbol{(2,(2p)^{\times 3\times 3})}}
\newcommand{\WWtrimult}{\WW(2,(4 j p^2 + 2 p)^{\times 3\times(4 j p + 3)})}

\newcommand{\Deltamatter}{\Delta^{\text{m}}}

\newcommand{\bK}{\mathsf{K}}
\newcommand{\bGamma}{\overline{\Gamma}}

\newcommand{\rightder}{\overleftarrow{\mathscr{D}}}

\newcommand{\KK}{\mathcal{K}\kern-5.7pt\raisebox{-3.9pt}{\footnotesize\textit{2}}\,}

\newcommand{\ccirc}{\mathbin{\raisebox{1pt}{\,$\scriptscriptstyle\circ$\,}}}

\newcommand{\WWtridwa}{\WW(2,(2p)^{\times 3\times 2})}
\newcommand{\WWW}{\boldsymbol{\mathscr{W}}}
\newcommand{\bWWtridwa}{\WW\myboldsymbol{(2,(2p)^{\times 3\times 2})}}
\newcommand{\bker}{\boldsymbol{\textup{ker}}\,}

\newcommand{\Xa}{X_{\mathrm{a}}}
\newcommand{\Xs}{X_{\mathrm{s}}}
\newcommand{\algUa}{{\mathbf{U}(\Xa})}
\newcommand{\algUs}{{\mathbf{U}(\Xs})}
\newcommand{\balgUa}{{\mathbf{U}\myboldsymbol{(X}_{\textbf{a}}\myboldsymbol{)}}}
\newcommand{\balgUs}{{\mathbf{U}\myboldsymbol{(X}_{\textbf{s}}\myboldsymbol{)}}}

\newcommand{\Ures}{\overline{\mathscr{U}}_{\q} s\ell(2)}
\newcommand{\Uresk}{\overline{\mathscr{U}}^{*}_{\q} s\ell(2)}

\newcommand{\algU}{{\mathbf{U}(X)}}
\newcommand{\algUU}{{\overline{\mathbf{U}}(X)}}
\newcommand{\balgU}{{\mathbf{U}\myboldsymbol{(X)}}}

\newcommand{\repZ}{\mathscr{Z}}

\newcommand{\Vsr}{\boldsymbol{V}\!}

\newcommand{\eqdef}{\mathrel{\mbox{\raisebox{1pt}{\footnotesize{:}}}{=}}}

\newcommand{\wh}{\widehat}
\newcommand{\dF}{\widehat{F}}

\newcommand{\Da}{\widehat{a}}
\newcommand{\Db}{\widehat{b}}
\newcommand{\DB}{\widehat{\Fi}}
\newcommand{\DF}{\widehat{\Fii}}

\newcommand{\bNich}{\boldsymbol{\mathfrak{B}}}
\newcommand{\Nich}{\mathfrak{B}}
\newcommand{\Nicha}{\Nich(X_{\asymm})}
\newcommand{\Nichs}{\Nich(X_{\symm})}

\newcommand{\symm}{\mathrm{s}}
\newcommand{\asymm}{\mathrm{a}}

\newcommand{\Fi}{B}
\newcommand{\Fii}{F}

\newcommand{\Ei}{C}
\newcommand{\Eii}{E}

\newcommand{\one}{\boldsymbol{1}}

\newcommand{\zS}{\boldsymbol{S}}

\newcommand{\AsymmF}{\mathscr{F}}
\newcommand{\SymmF}{\mathscr{G}}

\newcommand{\Tm}{T_{\text{m}}}
\newcommand{\cm}{c_{\text{m}}}

\newcommand{\Shuffle}{\mathop{\text{\mifody\sf Sh}}\nolimits}

\newcommand{\Bbin}[2]{{\Shuffle^{#1}_{#2}}}
\newcommand{\Bfac}[1]{\mathfrak{S}_{#1}}

\newcommand{\fobject}[1]{\object{\scriptstyle #1}}

\newcommand{\A}{\raisebox{.5pt}{\large$\mathpzc{S}\kern-1pt$}}

\newcommand{\leftact}{\kern1pt{\rightharpoonup}\kern1pt}

\newcommand{\cR}{\mathscr{R}}
\newcommand{\cU}{\mathscr{U}}
\newcommand{\cV}{\mathscr{V}}
\newcommand{\eval}[2]{\langle#1,\,#2\rangle\,}

\newcommand{\id}{\mathrm{id}}

\newcommand{\YDname}{\mathscr{Y\kern-1ptD}}
\newcommand{\GGyd}{{}\mbox{\small${}^{\Gamma}_{\Gamma}$}\YDname}
\newcommand{\BByd}{{}\mbox{\small${}^{\Nich(X)}_{\Nich(X)}$}\YDname}
\newcommand{\HHyd}{{}\mbox{\small${}^{H}_{H}$}\YDname}
\newcommand{\bHHyd}{{}\mbox{\small${}^{\myboldsymbol{H}}_{\myboldsymbol{H}}$}\boldsymbol{\YDname}}

\newcommand{\hA}{\mbox{\large$\mathpzc{s}\kern-.8pt$}}

\newcommand{\Smash}{\mathbin{\hash}}

\newcommand{\bup}[1]{^{{\scriptscriptstyle\!\underline{\,#1_{\vphantom{.}}\,}\!}}}

\newcommand{\modF}{\textup{\textsf{\textbf{F}}}} 
\newcommand{\modZ}{\textup{\textsf{\textbf{Z}}}}

\newcommand{\zero}{_{_{(0)}}}
\newcommand{\mone}{_{_{(-1)}}}

\newcommand{\bzero}{_{_{\,\underline{\,0\,}}}}
\newcommand{\bmone}{_{_{\underline{\!{-1}\!}}}}

\newcommand{\Azero}{f}
\newcommand{\Dvarphi}{\partial\Azero}

\newcommand{\Aone}{\varphi_1}
\newcommand{\Atwo}{\varphi_2}
\newcommand{\Aii}{\varphi}
\newcommand{\Athree}{\chi}
\newcommand{\Afour}{a}

\def\verma{\bullet\rule[-2pt]{0pt}{9pt}
  \ar@{-}@*{[|(2.5)]}[0,0]+<0pt,0pt>;[0,0]-<10pt,0pt>
  \ar@{-}@*{[|(2.5)]}[0,0]-<1pt,-1pt>;[0,0]+<7pt,-7pt>
  \ar@{}[0,0]+<0pt,0pt>;}

\def\vermaiv{\bullet\rule[-2pt]{0pt}{9pt}
  \ar@{-}@*{[|(4)]}[0,0]+<0pt,0pt>;[0,0]-<10pt,0pt>
  \ar@{-}@*{[|(4)]}[0,0]-<1pt,-1pt>;[0,0]+<7pt,-7pt>
  \ar@{}[0,0]+<0pt,0pt>;}

\def\Cverma{\circ\rule[-2pt]{0pt}{9pt}
  \ar@{-}@*{[|(2.5)]}[0,0]+<-3pt,0pt>;[0,0]-<10pt,0pt>
  \ar@{-}@*{[|(2.5)]}[0,0]+<2pt,-2pt>;[0,0]+<6pt,-7pt>
  \ar@{}[0,0]+<0pt,0pt>;}

\def\iverma{\bullet\rule[-2pt]{0pt}{9pt}
  \ar@{-}@*{[|(2.5)]}[]+<0pt,0pt>;[0,0]+<10pt,0pt>
  \ar@{-}@*{[|(2.5)]}[]+<1pt,1pt>;[0,0]+<-7pt,-7pt>
  \ar@{}[0,0]+<0pt,0pt>;}

\def\ivermaiv{\bullet\rule[-2pt]{0pt}{9pt}
  \ar@{-}@*{[|(4)]}[]+<0pt,0pt>;[0,0]+<10pt,0pt>
  \ar@{-}@*{[|(4)]}[]+<1pt,1pt>;[0,0]+<-7pt,-7pt>
  \ar@{}[0,0]+<0pt,0pt>;}

\def\Civerma{\circ\rule[-2pt]{0pt}{9pt}
  \ar@{-}@*{[|(2.5)]}[]+<3pt,0pt>;[0,0]+<10pt,0pt>
  \ar@{-}@*{[|(2.5)]}[]-<2pt,2pt>;[0,0]+<-6pt,-7pt>
  \ar@{}[0,0]+<0pt,0pt>;}

\def\nothing{*{\mbox{}}\ar@{};[0,0];}

\def\Bullet{*{\rule[-1pt]{0pt}{7pt}\bullet}\ar@{};[0,0];}
\def\Circ{*{\rule[-1pt]{0pt}{7pt}\circ}\ar@{};[0,0];}

\newcommand{\gr}{\mathop{\mathrm{gr}}\nolimits}
\newcommand{\rep}{\mathscr}
\newcommand{\ket}[1]{\mathchoice{%
    {\left|\smash[t]{#1}\right\rangle}}{|{#1}\rangle}{|{#1}\rangle}{|{#1}\rangle}}
\newcommand{\charSL}[2]{\chi_{%
    {\phantom{h}\kern-3pt #2}}^{\phantom{y}\kern-3pt #1}}
\newcommand{\hVerma}{\rep{M}}
\newcommand{\MFFplus}[1]{\mathbf{s}^+(#1)}
\newcommand{\MFFminus}[1]{\mathbf{s}^-(#1)}
\newcommand{\MFFname}{\mathbf{s}}
\newcommand{\jp}{\lambda^+}
\newcommand{\jm}{\lambda^-}

\newcommand{\bref}[1]{\textbf{\ref{#1}}}
\newcommand{\adj}{%
  \mathchoice{\mathbin{\blacktriangleright}}%
  {\mathbin{\mbox{\small${\blacktriangleright}$}}}%
  {\mathbin{{\blacktriangleright}}}%
  {\mathbin{{\blacktriangleright}}}}

\newcommand{\leftii}{%
  \mathchoice%
  {\mathbin{\mbox{\footnotesize${\triangleright}$}}}%
  {\mathbin{\scriptstyle\triangleright}}%
  {\mathbin{{\triangleright}}}%
  {\mathbin{{\triangleright}}}}

\newcommand{\thisfont}{\textsc}
\newcommand{\ab}[1]{\thisfont{a\kern-1pt b}_{#1}}
\newcommand{\ba}[1]{\thisfont{b\kern-1pt a}_{#1}}
\newcommand{\aba}[1]{\thisfont{a\kern-1pt b\kern-1pt a}_{#1}}
\newcommand{\bab}[1]{\thisfont{b\kern-1pt a\kern-1pt b}_{#1}}

\newcommand{\F}[1]{\thisfont{f\kern-1pt}_{#1}}
\newcommand{\Fb}[1]{\thisfont{f\kern-1pt b\kern-.5pt}_{#1}}
\newcommand{\bF}[1]{\thisfont{b\kern-1pt f\kern-1pt}_{#1}}
\newcommand{\bX}[1]{\thisfont{x\kern-1pt b\kern-1pt}_{#1}}
\newcommand{\bFb}[1]{\thisfont{b\kern-1pt f\kern-1pt b\kern-.5pt}_{#1}}

\newcommand{\FF}{\F}
\newcommand{\FB}{\Fb}
\newcommand{\BX}{\bX}

\newcommand{\BFB}{\bFb}

\newcommand{\Aint}[1]{\langle#1\rangle}
\newcommand{\Afac}[1]{\langle#1\rangle!\,}
\newcommand{\Abin}[2]{\mathchoice%
  {\Abinm{#1}{#2}\,}{\Abinmm{#1}{#2}\,}%
  {\Abinmm{#1}{#2}}{\Abinmm{#1}{#2}}}
\newcommand{\Abinm}[2]{\mbox{\footnotesize$\displaystyle
    \genfrac{\langle}{\rangle}{0pt}{}{#1}{#2}$}}
\newcommand{\Abinmm}[2]{\genfrac{\langle}{\rangle}{0pt}{}{#1}{#2}}

\newcommand{\Vertex}[1]{V^{\{#1\}}}

\newcommand{\tensor}{\mathbin{\otimes}}

\newcommand{\Dynkin}[3]{\xymatrix@C28pt@1{\circ\ar@{-}[r]^(.55){#2}\ar@{}^{#1}[]&\circ\ar@{}^{#3}[]}}

\newcommand{\mfrac}[2]{\raisebox{.8pt}{\mbox{\small$\displaystyle\frac{#1}{#2}$}}}
\newcommand{\ffrac}[2]{\raisebox{.5pt}{\mbox{\footnotesize$\displaystyle\frac{#1}{#2}$}}}
\newcommand{\fffrac}[2]{\raisebox{.9pt}{\mbox{\scriptsize$\displaystyle\frac{#1}{#2}$}}}
\newcommand{\half}{%
  \mathchoice{\ffrac{1}{2}}{\frac{1}{2}}{\frac{1}{2}}{\frac{1}{2}}}

\newcommand{\q}{\mathfrak{q}}
\newcommand{\xx}{\xi}
\newcommand{\hSL}[1]{\widehat{s\ell}(#1)}

\newcommand{\Co}[1]{{}#1{}\tensor{}}

\newcommand{\adjoint}{%
  \mathchoice{\mathbin{\blacktriangleright}}%
  {\mathbin{\mbox{\small${\blacktriangleright}$}}}%
  {\mathbin{{\blacktriangleright}}}%
  {\mathbin{{\blacktriangleright}}}}

\newcommand{\Ebeta}{\mathscr{E}}

\newcommand{\Lstate}{\mathscr{L}}
\newcommand{\Rstate}{\mathscr{R}}
\newcommand{\Ktop}{\mathscr{K}}
\newcommand{\Ktopleft}{\overset{\leftarrow}{\Ktop}}
\newcommand{\Ktopright}{\overset{\rightarrow}{\Ktop}}

\newcommand{\jplus}{j^+}
\newcommand{\jminus}{j^-}
\newcommand{\Jplus}{J^+}
\newcommand{\Jminus}{J^-}
\newcommand{\Jnaught}{J^0}

\newcommand{\oC}{\mathbb{C}}
\newcommand{\oZ}{\mathbb{Z}}
\newcommand{\oN}{\mathbb{N}}

\swapnumbers
\newtheorem{Thm}[subsection]{Theorem}

\theoremstyle{definition}

\begin{document}

\title{Logarithmic $\hSL2$ CFT{} models from Nichols algebras. 1}

\author{A.\,M.\;Semikhatov}
\author{I.\,Yu.\;Tipunin}

\address{Lebedev Physics Institute\hfill\mbox{}\linebreak
  \texttt{asemikha@gmail.com}, \ \texttt{tipunin@pli.ru}}

\begin{abstract}
  We construct chiral algebras that centralize rank-two Nichols
  algebras with at least one fermionic generator.  This gives
  ``logarithmic'' $W$-algebra extensions of a fractional-level $\hSL2$
  algebra.  We discuss crucial aspects of the emerging general
  relation between Nichols algebras and logarithmic CFT{} models:
  (i)~the extra input, beyond the Nichols algebra proper, needed to
  uniquely specify a conformal model; (ii)~a relation between the
  CFT{} counterparts of Nichols algebras connected by Weyl groupoid
  maps; and (iii)~the common double bosonization $\algU$ of such
  Nichols algebras.  For an extended chiral algebra, candidates for
  its simple modules that are counterparts of the $\algU$ simple
  modules are proposed, as a first step toward a functorial relation
  between $\algU$ and $W$-algebra representation categories.
\end{abstract}

\maketitle
\thispagestyle{empty}

\section{\textbf{Introduction}}
Logarithmic models of two-dimensional conformal field theory with
$\hSL2$ symmetry have been addressed in
\cite{[Gab-frac+],[N-SU2],[LMRS],[QQ-sl2],[arche],[Rid-0810+],
  [Rid-Cre]}.  Here, we approach logarithmic $\hSL2$ models from the
theory of Nichols
algebras~\cite{[Nich],[Wor],[Lu-intro],[Rosso-inv],[AG],[AS-pointed],
  [Andr-remarks],[AS-onthe],[Heck-Weyl],[Heck-class],[AHS],[ARS],
  [Ag-0804-standard],[Ag-1008-presentation]}.  The major advantages
are that (i)~we then have a recipe for seeking the \textit{extended
  chiral algebra}: as (the maximal local algebra in) the kernel of
Nichols algebra generators represented by screening
operators~\cite{[FHST]}, and (ii)~very strong hints about the chiral
algebra representation theory come from the Nichols-algebra side.

Our starting point is therefore brazenly algebraic---a braiding
matrix $(q_{i,k})_{\substack{1\leq i\leq\text{rank}\\
    1\leq k\leq\text{rank}}}$, which is just a collection of factors
defining the diagonal braiding of generators of a Nichols algebra.
Finite-dimensional Nichols algebras with diagonal braiding have been
impressively classified \cite{[Heck-class]}; in the spirit
of~\cite{[STbr],[c-charge]}, this calls for translating the available
structural results into the language of logarithmic CFT{} models.  We
here go beyond the case of rank-$1$ Nichols
algebras~\cite{[b-fusion]}, and in rank~$2$ select two braiding
matrices
\begin{equation}\label{qij-first}
  Q_{\asymm} = 
  \begin{pmatrix}
    -1 & (-1)^j \q^{-1} \\
    (-1)^j \q^{-1} & \q^2
  \end{pmatrix}
  \quad\text{and}\quad
  Q_{\symm} = 
  \begin{pmatrix}
    -1 & (-1)^j \q \\
    (-1)^j \q & -1
  \end{pmatrix},\quad
  \q = e^{\frac{i\pi}{p}}
\end{equation}
with an integer $p\geq 2$.

A \,characteristic \,result \,deduced \,from \eqref{qij-first} \,is
\,the \,triplet--triplet \,extended \,algebra $\WWtritri$ that
centralizes\pagebreak[3] the Nichols algebra (the one corresponding to
$Q_{\asymm}$ with $j=0$, for definiteness): it extends $\hSL2$ at the
level $k=\frac{1}{p} - 2$ and is generated by dimen\-sion-$2p$ fields
\begin{equation*}
  \W^+(z) = \AsymmF_1(z),
  \qquad
  \W^-(z) =\MFFminus{2,p}\AsymmF_{-1}(z),
\end{equation*}
where $\AsymmF_{\pm1}(z)$ are $\hSL2$ primaries of charge $\pm1$,
constructed in a free-field realization associated with the chosen
braiding matrix, and $\MFFminus{2,p}$ is an $\hSL2$ singular vector
operator.  The algebra is triplet--triplet because $\W^{\pm}(z)$,
together with a~$\W^{0}(z)$, make up a triplet, and in addition each
$\W^{+,0,-}(z)$ is part of a triplet under the zero-mode
subalgebra~of~$\hSL2$.

Other results include a similar construction for the second braiding
matrix, the ``triplet--multiplet'' $W$-algebras that extend some other
fractional-level $\hSL2$, and links between the extended algebras and
certain Hopf-algebraic counterparts.  The reader is also invited to
``\textit{Logarithmic $\hSL2$ CFT{} models from Nichols
  algebras. $N$}'' with $N>1$ for modular transformations of the
extended (pseudo)char\-acters, fusion of representations, etc.  Part
of our effort, in this paper especially, is given to grasping the
general principles governing how a given Nichols algebra can be
``mapped'' into a logarithmic~CFT{}.

\subsubsection*{``Non-Nichols-algebra'' data for logarithmic CFTs}
We refrain from claiming, except in private, that a collection of
roots of unity $q_{i,k}$ contains exactly as much information as is
needed for uniquely reconstructing a CFT{} model.  What precisely the
extra input consists of was a major part of the intrigue in writing
this paper.  One such input is $j$ in each case in~\eqref{qij-first}.
When realizing the $q_{i,k}$ as the braiding of screening operators,
the range of $j$ is promoted from $\oZ_2$ to $\oZ$ (hence an
arbitrariness), with drastic consequences for logarithmic theories
with different~$j$ (we actually consider $j=0,1,2,\dots$ for
$Q_{\asymm}$ and $j=-1, 0,1,2,\dots$ for~$Q_{\symm}$).  The extended
algebra is \textit{triplet--triplet} for a single value of $j$ in each
case, and \textit{triplet--multiplet} for other $j$.

\subsubsection*{More on the relevant integers: $\myboldsymbol{p}$,
  $\myboldsymbol{j}$, and $\myboldsymbol{p'}$} The theory of Nichols
algebras only requires that for each case in~\eqref{qij-first}, $\q^2$
be a primitive $p$th root of unity.  If
$\tilde{\q}=e^{i\pi\frac{\tilde p}{p}}$ (with an odd $\tilde p$
coprime with $p$) is chosen instead of $\q$ in~\eqref{qij-first}, then
the screening operators whose braiding reproduces the braiding matrix
depend on two full-fledged coprime integers, $p$ and $p' \eqdef j p +
\tilde p$.  This is how $(p,p')$-type logarithmic models appear in our
context.  To somewhat reduce technical details, we restrict ourself to
the case $p' = j p + 1$; the relevant $\hSL2$ machinery is then
structurally the same as for the general $p'$ coprime with $p$, but
slightly simpler to follow.

Thus, setting
\begin{equation}\label{the-q}
  \q = e^{\frac{i\pi}{p}}
\end{equation}
in \eqref{qij-first} and introducing the $j$ parameter, we effectively
deal with $2p$th roots of unity of the form $e^{i\pi\frac{j p +
    1}{p}}$, which is the source of $(p, j p + 1)$ logarithmic models
that occur in what follows.

\subsubsection*{Nichols algebras and fermionic screenings}
The Nichols algebras associated with braiding matrices
in~\eqref{qij-first} are selected from the list of about two dozen
rank-two Nichols algebras with diagonal braiding \cite{[Heck-1+1]}
based on two properties:
\begin{enumerate}
\item these are not isolated but ``serial'' algebras, existing for
  each integer $p\geq 2$ (to become the $p$ in the labeling of
  logarithmic CFT{} models), and

\item at least one generator of the Nichols algebra has the
  self-braiding factor $q_{i,i}=-1$, which is standardly rephrased by
  saying that this generator is a fermion. 
\end{enumerate}

Because Nichols algebra generators $F_j$ are identified with screening
operators, the $F_i$ with $q_{i,i}=-1$ are generally referred to as
\textit{fermionic screenings}.  Fermionic screenings occur in CFT{}
models rather frequently, but are nevertheless slightly disquieting
from a ``conceptual'' standpoint because of the lack of an obvious
relation to root systems: while some other ``good'' systems of
screenings are reproducible by rescaling classic root systems,
typically by factors like $1/\sqrt{p}$, those with fermionic
screenings are not.  With the power of Nichols algebras, however, this
complication passes unnoticed; moreover, there is a procedure to
extract a (generalized) Cartan matrix $A=(a_{i,j})$ from a braiding
matrix (dependent on several ``$p$-type'' parameters in general).
Such a Cartan matrix is an important part of the connection between
Nichols algebras and conformal models, as we see shortly.

The classes to which braiding matrices~\eqref{qij-first} belong are
defined by the conditions~\cite{[Heck-1+1]}
\begin{equation}\label{eq:three-Q}
  \begin{aligned}
    Q_{\asymm}:\qquad&q_{11}=-1,\quad q_{12}q_{21}q_{22}=1,\quad 
    q_{12}q_{21} \in \mathsf{R}_p,
    \\
    Q_{\symm}:\qquad&q_{11}=-1,\quad q_{12}q_{21}\in \mathsf{R}_p,
    \quad q_{22}=-1,
    \\
    \overline{Q}_{\asymm}:\qquad&q_{12}q_{21}q_{11}=1,\quad 
    q_{12}q_{21} \in \mathsf{R}_p,\quad q_{22}=-1,
  \end{aligned}
\end{equation}
where $\mathsf{R}_p$ is the set of primitive $p$th roots of unity,
$p\geq 2$, and where we temporarily (until~\eqref{Weyl-on-R}) add the
third case that is merely a $1\leftrightarrow 2$ relabeling of the
first.  In each case, the associated Cartan matrix is the $A_2$ one,
$(a_{i,j}) =\left(\begin{smallmatrix}
    2&-1\\
    -1&2
  \end{smallmatrix}\right)$.  The Nichols algebra has dimension $4p$
in each case.

\subsubsection*{From braiding to screenings}
The point of contact with CFT{}, already mentioned in the foregoing,
is the interpretation of Nichols algebra generators as screenings,
\begin{equation}\label{Fi}
  F_i = \oint e^{\alpha_i . \varphi},\quad 1\leq i \leq
  \text{rank}\equiv \theta,
\end{equation}
where $\varphi(z)$ is a $\theta$-plet of scalar fields with standard
normalization, the dot denotes Euclidean scalar product, and
$\alpha_i\in\oC^{\theta}$ are called the momenta of the screenings.
The relation between the screening momenta and the braiding matrix is
postulated~\cite{[c-charge]} in the form of $\half\theta(\theta+1)$
equations
\begin{equation}\label{q-to-alpha}
  e^{i\pi \alpha_j.\alpha_j} = q_{j,j},\qquad
  e^{2i\pi \alpha_j.\alpha_k} = q_{j,k}q_{k,j}
\end{equation}
and the $\theta^2 - \theta$ logical-``or'' conditions
\begin{equation}\label{cartan-test-log}
  a_{i,j} \alpha_i.\alpha_i = 2\alpha_i.\alpha_j
  \quad\mbox{\small$\bigvee$}\quad
 (1-a_{i,j})\alpha_i.\alpha_i = 2
\end{equation}
imposed for each pair $i\neq j$ and involving the Cartan matrix
$a_{i,j}$ associated with the given braiding
matrix.

For the braiding matrices defined by each line in~\eqref{eq:three-Q},
we solve Eqs.~\eqref{q-to-alpha}--\eqref{cartan-test-log} (with
$\theta=2$) for the scalar products of the momenta, with $\q$ as
in~\eqref{qij-first}. \ The respective solutions are
\begin{equation}\label{three-R}
  \begin{aligned}
    R_{\asymm}(j):&\qquad
    \alpha_1.\alpha_1 = 1,
    \quad \alpha_1.\alpha_2 =  -\fffrac{j p + 1}{p},
    \quad \alpha_2.\alpha_2 = 2 \fffrac{j p + 1}{p},
    \\
    R_{\symm}(j):&\qquad
    \alpha_1.\alpha_1 = 1,
    \quad\alpha_1.\alpha_2 = \fffrac{j p + 1}{p},
    \quad\alpha_2.\alpha_2 = 1,
    \\
    \overline{R}_{\asymm}(j):&\qquad
    \alpha_1.\alpha_1 = 2 \fffrac{j p + 1}{p},
    \quad \alpha_1.\alpha_2 =  -\fffrac{j p + 1}{p},
    \quad \alpha_2.\alpha_2 = 1,
  \end{aligned}
\end{equation}
with an arbitrary $j\in\oZ$ in each case.  These three systems of
$\oC^2$ vectors are related by Weyl-groupoid pseudoreflections, which
map as (`$1$' with respect to $\alpha_1$ and `$2$' with respect to
$\alpha_2$)
\begin{equation}\label{Weyl-on-R}
  \xymatrix@R=12pt@C=12pt{
    &&R_{\symm}(j-1)
    \ar@{{<}{-}{>}}_2[lld]
    \ar@{{<}{-}{>}}^1[rrd]&&&
    \\
    \overline{R}_{\asymm}(j)
    \ar@{{<}{-}{>}}@(ld,rd)[]_{1}&&&&
    R_{\asymm}(j)
    \ar@{{<}{-}{>}}@(ld,rd)[]_{2}
    }
\end{equation}
(in our case of a classic, $A_2$, Cartan matrix, it is the same for
the entire orbit, and therefore the Weyl groupoid becomes the
corresponding Weyl group).  

The third, ``opposite asymmetric'' case does not have to be considered
separately.  In what follows, we therefore speak about
two---``asymmetric'' and ``symmetric''---systems of the screening
momenta, braiding matrices, braided spaces $\Xa$ and $\Xs$, Nichols
algebras $\Nich(\Xa)$ and $\Nich(\Xs)$, etc.; in the narrow, literal
sense, the \textit{symmetry} is with respect to the anti-diagonal of
$Q_{\symm}$ in~\eqref{qij-first}, but it also shows up on the CFT{}
side, as we describe below.

\subsubsection*{$\myboldsymbol{2\neq 3}$: From parafermions to
  $\myboldsymbol{\hSL2}$}
Taking $\text{rank}\!=\!2$ means that the kernel of two screenings is
defined more or less uniquely in a two-boson space; there, the kernel
contains parafermionic fields $j^+(z)$ and $j^-(z)$---generators of
the $\hSL2/u(1)$ coset theory, which are mutually nonlocal.
Apparently, one can proceed within the theory of para\-fermions and
their ``representations,'' but this is a challenge we are not up to,
in this paper at least. Instead, we introduce a third, ``auxiliary''
scalar field representing the $u(1)$ in the above coset, with the
OPE\footnote{We use noncanonical normalizations.  The canonical ones
  are of course easy to restore, but at the expense of factors like
  $\sqrt{\pm k}$, entailing either a choice for the sign of $k$ or
  unnecessary imaginary units.}
\begin{equation}\label{localize-2}
  \Athree(z)\, \Athree(w) = 2k\log(z-w),
\end{equation}
where $k=\frac{1}{p} + j - 2$ in the asymmetric case and
$k=\frac{1}{p} + j - 1$ in the symmetric case, and use it to convert
the parafermions into level-$k$ $\hSL2$ currents
\begin{equation*}
  J^{\pm}(z) = j^{\pm}(z) e^{\pm\frac{1}{k}\Athree(z)}.
\end{equation*}
We can then use the full power of the $\hSL2$ representation theory,
but the price is the occurrence of twisted (spectral-flow transformed)
$\hSL2$ representations, because the locality requirement alone leaves
a size-$\oZ$ arbitrariness in how the third scalar enters the relevant
operators.

Eliminating $\Athree(z)$ in order to return to a ``clean''
Nichols-algebra-motivated two-boson conformal model is not difficult
at the level of formulas (for generators of the extended algebra,
etc.), and can be quite interesting at the level of characters
(cf.~\cite{[ABC]}), but the meaning of ``algebras'' generated by
mutually nonlocal fields would then have to be worked out first,
before the ``three-boson'' results in this paper can be restated in
that language.\footnote{Operator product~\eqref{localize-2} conceals
  yet another size-$\oZ$ arbitrariness in mapping from a given Nichols
  algebra to CFT{}. \ Instead of~\eqref{localize-2}, we could have
  chosen the OPE
\begin{equation*}
  \Athree(z)\, \Athree(w) = (2 + k n)k\log(z-w),\quad n\in\oZ.
\end{equation*}
The parafermions would then be ``dressed'' not into $J^{\pm}(z)$ but
into local fields $\mathscr{J}^{\pm}(z)$ with the OPEs
$\mathscr{J}^{\pm}(z)\mathscr{J}^{\pm}(w)\propto (z-w)^n$ and
$\mathscr{J}^{+}(z)\mathscr{J}^{-}(w)=
\mathrm{const}\cdot{}\linebreak[1](z-w)^{-n-2}+\dots$. \ For $n=1$, in
particular, $\mathscr{J}^+(z)$ and $\mathscr{J}^-(z)$ generate the
$N=2$ super-Virasoro algebra.  The ``parafermionic core'' of all these
theories with different $n$ is the same.  The value $n=0$ chosen
in~\eqref{localize-2} and corresponding to $\hSL2$ is the lowest in
the sense that taking $n\leq -1$ leads to unconventional theories
where $\mathscr{J}^{+}(z)$ (and $\mathscr{J}^{-}(z)$) has a
nonvanishing OPE with itself.  We do not return to this point here and
stay with~\eqref{localize-2}.}

\subsubsection*{From $\myboldsymbol{\hSL2_k}$ to $\WW$ algebras}
The $\hSL2$ algebra is used as a seed to grow extended algebras living
in the kernel of the screenings inside $\hSL2$ representations.  The
relevant representations are Wakimoto-type modules~\cite{[W],[FF]}.
These modules are related by intertwining maps, which can also be
associated with screenings~\cite{[BMP]}, and the intersection of the
screening kernels we are interested in is the socle of such modules
(cf.~\cite{[QQ-sl2]}).

In the kernel, we select mutually local fields, and among them,
$\hSL2$ primaries with the minimal Sugawara dimension as generators of
the \textit{extended algebra}---a $W$-algebra of mutually local fields
in the kernel.
The generators come in a triplet $\W^+(z)$, $\W^0(z)$, and $\W^-(z)$
and, moreover, each of these is also a multiplet with respect to the
zero-mode $s\ell(2)$ subalgebra of~$\hSL2$. \ We construct the
$W$-algebra generators as explicitly as the celebrated
construction~\cite{[MFF]} of $\hSL2$ singular vectors allows.  The
$W$-algebras constructed in the asymmetric and symmetric realizations
are isomorphic and are in fact related as two ``rebosonizations'' of
the Wakimoto representation, as we discuss
in~\bref{sec:Wak}.\enlargethispage{\baselineskip}

For the asymmetric realization (the top row in~\eqref{three-R}) for
definiteness, with $j\geq 0$, the three generators have the (Sugawara)
dimension $2 p (2 j p+1)$ and each is part of a $(4 j p + 3)$-plet
with respect to the zero-mode $s\ell(2)$ algebra.  The cases $j=0$
(dimension $2p$ and zero-mode \textit{triplets})\pagebreak[3] and
$j\geq1$ are essentially different.  For $j=0$, the $\W^+(z)$ field is
a pure exponential, and there exists a ``long'' screening\footnote{The
  screenings that are identified with Nichols algebra generators are
  sometimes conventionally called ``short'' screenings, as opposed to
  ``long'' screenings, which centralize the nonextended algebra
  (Virasoro in~\cite{[FHST]}, $W_3$ in~\cite{[log-W3]}, and $\hSL2$ in
  this paper).} mapping as $\W^+(z)\stackrel{\Ebeta}{\mapsto}\W^0(z)
\stackrel{\Ebeta}{\mapsto}\W^-(z)$; for $j\geq 1$, by contrast, all
three fields (denoted as $\jW^+(z)$, $\jW^0(z)$, and $\jW^-(z)$) have
the form of a differential polynomial times an exponential and the
long screening does not map that way.  Another characteristic
illustration of the difference is the result of Hamiltonian reduction
to single-boson models: to the $\pone$ triplet algebra
\cite{[Kausch],[GK+],[FHST],[AM-3+latt]} for $j=0$ and to the triplet
algebras of $(p,p')$ logarithmic models introduced in~\cite{[FGST3]}
(also see~\cite{[AM-3+latt],[AM-pp']}) for $j\geq 1$ (with $p'=p j +
1$ due to our choice of $\q$ in~\eqref{the-q}, as discussed
above).\footnote{This already shows that when it comes to comparing
  appropriate categories, the extra data used explicitly or implicitly
  in ``reconstructing'' a logarithmic model from a Nichols algebra can
  play a decisive role: the categories compare reasonably well for
  $\pone$ models \cite{[FGST2],[NT],[TW1]} and somewhat worse for
  $(p,p')$ models~\cite{[FGST3],[FGST4],[GRW],[RGW],[TW2]}.  Modular
  transformations are remarkably robust, however: they are the same on
  the CFT{} and Hopf-algebraic sides of $(p,p')$
  models~\cite{[FGST3],[FGST4]}.}

Our ``$W$-style'' notation for the triplet--triplet algebra
(corresponding to $j=0$ in the asymmetric realization) is $\WWtritri$,
where $2$ is the dimension of the energy--momentum tensor and $2p$ is
the dimension of all the $\W^{+,0,-}(z)$ fields.  The first
multiplicity of $3$ is for the superscript of the~$\W^{+,0,-}(z)$
fields, and the second, for the fact that each of these fields is a
triplet under the zero-mode~$s\ell(2)$.

\subsubsection*{Choosing an $\bHHyd$ category}
The relation between $\Nich(X)$ and the extended algebra starts
working at its full strength when lifted to the level of a
correspondence between the categories of their representations.  On
the Nichols-algebra side, the relevant representations are
Yetter--Drinfeld $\Nich(X)$ modules.  The class of models that we wish
to have in CFT{} should be ``rational in a reasonably broad sense''
(certainly with finitely many simple modules and possibly with the
$C_2$ cofiniteness property~\cite{[CF]}, etc.; ultimately, the
rigorous framework of~\cite{[HLZ]} must be kept in mind).  To avoid
``strongly nonrational'' effects, the scalar fields must take values
in a torus; in other words, a particular lattice vertex-operator
algebra should be selected.\enlargethispage{\baselineskip}

On the Nichols-algebra side, the corresponding
finiteness$/$discreteness is controlled by choosing a nonbraided Hopf
algebra $H$ used for ``YDinization,'' i.e., representing all relevant
braided spaces (starting with $X$) as objects in~$\HHyd$ (the category
of Yetter--Drinfeld $H$-modules).  Such an $H$ is by no means
unique~\cite{[T-survey]}, and its choice is another input beyond the
Nichols algebra itself, because $\Nich(X)$ is independent of an $H$
such that $X\in\HHyd$, but the resulting representation categories
acquire such a dependence.  For diagonal braiding, we have $H=k\Gamma$
for an Abelian group~$\Gamma$.\footnote{Our base field is $\oC$, but
  we write $k\Gamma$ for (perhaps misinterpreted) esthetic reasons.} \
The choice of $\Gamma$ is in turn related to the representation of
screenings in terms of free bosons in~\eqref{Fi}: once scalar fields
are introduced,\pagebreak[3] there is a zero-mode operator $\varphi_0
^{(i)}$ for each scalar $\varphi^{(i)}(z)= \overline{\varphi}^{(i)} +
\varphi_0^{(i)}\log z + \sum_{n\in\oZ\setminus 0}\frac{1}{n}
\varphi_n^{(i)} z^{-n}$, and operators
$e^{i\pi\kappa\varphi_0^{(i)}}$, with $\kappa$ determined by the
chosen lattice vertex-operator algebra, can be considered generators
of $\Gamma$.

We do not quite study $\BByd$ for $X\in\HHyd$ in this paper, however,
the reason being that, the results in~\cite{[b-fusion]}
notwithstanding, we made a much faster progress within a more
``old-fashioned'' setting of $\algU$-modules for a double bosonization
$\algU = \Nich(X^*)\tensor\Nich(X)\tensor H$ of our Nichols algebras.

\subsubsection*{Double bosonization: $\balgU$}
Double bosonization of braided Hopf algebras was discussed previously
in more general contexts in~\cite{[Majid-double],[Sommerh-deformed]};
it is to produce a nonbraided Hopf algebra from a braided Hopf algebra
$\cR$ in $\HHyd$, its dual $\cR^*$, and the Hopf algebra $H$ itself,
under a number of conditions satisfied by all the ingredients.  The
construction ``doubles'' the standard Radford bosonization (byproduct
formula)~\cite{[Radford-bos]} whereby the smash product $\cR\Smash H$
is endowed with the structure of a Hopf algebra.  For a Nichols
algebra $\Nich(X)$ with diagonal braiding and for a commutative
cocommutative $H=k\Gamma$, the list of double-bosonization conditions
is somewhat reduced, the most essential remaining one being the
symmetricity of the braiding in the standard sense $q_{i,j}=q_{j,i}$.
Double bosonization then yields a (nonbraided) Hopf algebra
\begin{equation}\label{algU}
  \algU = \Nich(X^*)\tensor\Nich(X)\tensor k\Gamma,
\end{equation}
which contains the ``single'' bosonizations $\Nich(X)\Smash k\Gamma$
and $\Nich(X^*)\mathbin{{\Smash}'} k\Gamma$ as Hopf subalgebras, where
the prime indicates that the $k\Gamma$ action and coaction are changed
by composing each with the antipode (they remain a left action and a
left coaction for $H=k\Gamma$).

This double bosonization previously appeared
in~\cite{[HY],[ARS],[HS-double]}.  We arrive at~\eqref{algU} in a way
that may be interesting in and of itself, the key observation being
that even without the assumption that $X$ is in $\HHyd$ (but with
diagonal braiding), there is an associative algebra structure on
$\algUU = \Nich(X^*)\tensor\Nich(X)\tensor k\bGamma$ for an Abelian
group $\bGamma$ that is read off from the monodromy $\Psi^2:X\tensor
X\to X\tensor X$.  This associative algebra is only ``half the way''
to~\eqref{algU} because, not surprisingly, $\bGamma$ is ``too coarse''
to endow $(X,\Psi)$ with the structure of a Yetter--Drinfeld module,
and for essentially the same reason, $\algUU$ is not a Hopf algebra.
In fact, $\bGamma$ is generated by the squares of the generators of a
``minimal'' $\Gamma$ featuring in~\eqref{algU}, and $\algUU$ can be
considered a subalgebra in~$\algU$.

Applying~\eqref{algU} to the two (nonisomorphic) rank-two Nichols
algebras studied here, the ``asymmetric'' $\Nich(\Xa)$ and
``symmetric'' $\Nich(\Xs)$, yields two $64p^4$-dimensional Hopf
algebras $\algUa$ and $\algUs$.  These turn out to be isomorphic as
associative algebras:
\begin{equation}\label{ssigma}
  \ssigma:\algUs\stackrel{\sim}{\to}\algUa
\end{equation}
Moreover, the two coproducts, $\Delta_{\asymm}$ and $\Delta_{\symm}$,
are related by a similarity transformation: there exists an invertible
element $\Phi\in\algUa\tensor\algUa$ satisfying an appropriate cocycle
condition such that
\begin{equation}\label{PHI}
  \Phi^{-1}\Delta_{\asymm}(\ssigma(x))\Phi =
  (\ssigma\tensor\ssigma)\ccirc \Delta_{\symm}(x)
\end{equation}
for all $x\in\algUs$. \ Hence, in particular, the representation
categories of $\algUa$ and $\algUs$ are equivalent as monoidal
categories.


The representation theory of $\algU$ is expected, at the very least,
to produce hints as to the representation theory on the CFT{} side;
but the hope---and the subject of future research---is that the
category of $\algU$-modules is equivalent to the representation
category of a chiral algebra.  Here is a subtlety, however, which has
no clear-cut analogue in the $\Wpone$$/$Virasoro case.

\subsubsection*{The $\bWWtridwa$ algebra} 
Only a subalgebra in $\WWtritri$ commutes with the generators of the
Abelian group $\Gamma$ that fixes the underlying Yetter--Drinfeld
category or, equivalently, allows constructing the double
bosonization.  This subalgebra does not contain $\Jplus(z)$ and
$\Jminus(z)$ currents, but contains their squares (normal-ordered
products) $(\Jplus(z))^2$ and $(\Jminus(z))^2$.  For consistency, the
fields sitting in the middle of each triplet with respect to the
zero-mode $s\ell(2)$ also have to be dropped.  The triplets then
become ``doublets'' with respect to $((J^{\pm}(z))^2)_0$, and the
notation for this algebra is~$\WWtridwa$.

The representation theory of $\algU$ (developed to some depth
in~\cite{[nich-sl2-2]}) and the representation theory of $\WWtridwa$
(rudimentary) show as much agreement as the ``rudimentary'' one
allows.  In particular, there are natural candidates for
irreducible $\WWtridwa$ modules.  The irreducible modules are divided
into spectral flow orbits, which are in a $1:1$ correspondence (not
only in number, but also as $\Gamma$-modules) with the $\algU$ simple
modules.

It is for the representation category of $\WWtridwa$ that we expect a
number of ``good'' properties resembling those of
$\Wpone$-models~\cite{[FGST]}.  Although the existing data is scarce,
it is very tempting to conjecture that much information about the
representation category of $\WWtridwa$ (including fusion and,
possibly, modular transformations) can be obtained from investigating
the $\algU$ category.

\subsubsection*{Notation and conventions}
The braiding matrix elements are defined such that
\begin{equation}\label{qij}
  \Psi: F_i\tensor F_k\mapsto q_{i,k} F_k\tensor F_i,\quad
  1\leq i,k\leq\text{rank}.
\end{equation}

The notations for conformal fields, such as $\Aone(z)$ and $\Atwo(z)$,
are ``local'' to each case, asymmetric (Sec.~\ref{CFTasymm}) and
symmetric (Sec.~\ref{CFTsymm}): they are for \textit{different} scalar
fields, whose OPEs are related to the respective system of scalar
products in~\eqref{three-R}.\enlargethispage{\baselineskip}

Our $\hSL2$ conventions are given in~\bref{app:sl2-conv}.

We use $q$-integers, factorials, and binomials defined as
\begin{equation*}
  \Aint{n}=\ffrac{q^{2n}-1}{q^2-1},\qquad
  \Afac{n}=\Aint{1}\dots\Aint{n},\qquad
  \Abin{m}{n}=\ffrac{\Afac{m}}{\Afac{n}\Afac{m-n}} \quad(m\geq n\geq0),
\end{equation*}
all of which are assumed specialized to $q=\q$ in~\eqref{the-q}.


\section{\textbf{Nichols algebra: the ``asymmetric''
    case}}\label{sec:Nich-asymm}
We explicitly describe the Nichols algebra $\Nich(X)$ corresponding to
the first line in~\eqref{eq:three-Q}: we take $X=\Xa$ to be a braided
vector space with basis $\Fi$, $\Fii$ and with the diagonal braiding
in this basis specified by the braiding matrix
\begin{equation}\label{qij-asymm}
  (q_{ij})=
  \begin{pmatrix}
    -1&\xx^{-1}\q^{-1}\\
    \xx \q^{-1}&\q^2
  \end{pmatrix}
\end{equation}
(the conventions are such that, e.g.,
$\Psi(\Fi\tensor\Fii)=\xx^{-1}\q^{-1}\Fii\tensor\Fi$).\footnote{The
  symbol $\Fi$ chosen for a fermionic generator may not suggest the
  best mnemonics; the association, possibly somewhat far-fetched, is
  with the physicists' notation $b,c$ for a pair of fermionic
  (``ghost'') fields.  Indeed, a $\Ei$ appears as a pair to $B$ in
  what follows.}  \ Here, $\xx$ is any $2p$th root of unity,
$\xx^{2p}=1$; this is an ``inessential'' variable, which can be
eliminated (``set equal to~$1$'') by a twist map associated with a
certain 2-cocycle~\cite{[AS-onthe]}.  In the applications in what
follows, we restrict to the case where $\xx^2=1$ anyway, but we
nevertheless keep $\xx$ in this section, to add some flavor to the
explicit formulas.

\subsection{$\boldsymbol{\Nich(X)}$: a
  presentation}\label{asymm-quotient}
For the braiding matrix in~\eqref{qij-asymm}, the Nichols algebra
$\Nich(X)$ is the quotient~\cite{[Ag-0804-standard]} (also
see~\cite{[Hel]})
\begin{equation}\label{eq:asymm-quotient}
  \Nich(X)=T(X)/ \bigl([\Fii,[\Fii,\Fi]],\ \Fi^2,\ \Fii^{p}\bigr)
\end{equation}
if $p\geq 3$ and
\begin{equation}\label{eq:asymm-quotient-p=2}
  \Nich(X)=T(X)/ \bigl(\Fi^2,\ [\Fi,\Fii]^2,\ \Fii^{2}\bigr)
\end{equation}
if $p=2$, with $\dim\Nich(X)=4p$ in all cases; here, square brackets
denote $q$-commutators, $[\Fii,\Fi]=\Fii \Fi - \xi\q^{-1} \Fi \Fii $
and $[\Fii,[\Fii,\Fi]] = \Fii [\Fii,\Fi] - \xi\q [\Fii ,\Fi]\Fii $.
The double commutator in~\eqref{eq:asymm-quotient} therefore
represents the relation
\begin{equation}\label{eq:asymm-Serre}
  \xx^2 \Fi \Fii^2 - \xx (\q + \q^{-1})  \Fii \Fi \Fii + 
  \Fii^2  \Fi = 0
\end{equation}
This implies that
\begin{equation}\label{FBFB}
  \Fi \Fii \Fi \Fii - \xx^{-2} \Fii \Fi \Fii \Fi = 0.
\end{equation}
For $p=2$, this last relation is a rewriting of $[\Fi,\Fii]^2=0$ 
under the condition that $\Fi^2 = \Fii^{2} = 0$.

The multiplication on $T(X)$ and $\Nich(X)$ understood
in~\eqref{eq:asymm-quotient} and the subsequent formulas is by
``concatenation'' $X^{\otimes m}\tensor X^{\otimes n}\to
X^{\otimes(m+n)}$, which simply maps
$(x_1,\dots,x_m)\tensor(y_1,\dots,y_n)\mapsto
(x_1,\dots,x_m,y_1,$\linebreak[0]$\dots,y_n)$; comultiplication is
then by ``deshuffling.''

\subsection{$\boldsymbol{\Nich(X)}$ as a subalgebra}
Another description of any Nichols algebra $\Nich(X)$, which is in
fact the description\pagebreak[3] in terms of screening operators, is
not as a quotient of but as a subspace in $T(X)$---as the algebra
generated by basis elements of $X$ under the shuffle
product~\cite{[Rosso-inv],[STbr]}.  The two realizations of $\Nich(X)$
are related by the total braided symmetrizer map $\Bfac{n}:X^{\otimes
  n}\to X^{\otimes n}$ in each grade.  It is an isomorphism precisely
because in the presentation
$\Nich(X)^{(n)}=T(X)^{(n)}/\mathscr{I}^{(n)}$, the ideal
$\mathscr{I}^{(n)}$ is the kernel of $\Bfac{n}$.

\subsubsection{}\label{sec:a-basis}
In this language, the $\Nich(X)$ in~\bref{asymm-quotient} is the
linear span of the $4p$ PBW elements
\begin{align*}
  \F{n} &= \fffrac{1}{\Afac{n}}\underbrace{\Fii *\Fii
    *{}\dots{}*\Fii}_{n},\quad 0\leq n\leq p-1
  \\[-6pt]
  \intertext{(with $\F{0}=\one$),}
  \Fb{n} &= \fffrac{1}{\Afac{n-1}}\underbrace{\Fii *\Fii *{}\dots{}*\Fii}_{n-1}{}*\Fi,\quad 1\leq n\leq p,
  \\
  \bX{n} &= \fffrac{1}{\Afac{n - 2}} \bigl(
  \underbrace{\Fii *\Fii *{}\dots{}*\Fii}_{n-2}{}*\Fi *\Fii - \xx^{-1} \q
  \underbrace{\Fii *\Fii *{}\dots{}*\Fii}_{n-1}{}*\Fi\bigr),
  \quad 2\leq n\leq p+1,
  \\
  \bFb{n} &= \fffrac{1}{\Afac{n-3}}\underbrace{\Fii *\Fii *{}\dots{}*\Fii}_{n-3}{}*\Fi *\Fii *\Fi,
  \quad 3\leq n\leq p+2,
\end{align*}
where $*$ denotes the shuffle product associated with the given
braiding,
\begin{equation*}
  (x_1,\dots,x_m)\tensor(y_1,\dots,y_n)\mapsto
  \Bbin{}{m,n}(x_1,\dots,x_m,y_1,\dots,y_n)
\end{equation*}
(our conventions are just as in~\cite{[STbr]}, except that $*$ was not
used for the shuffle product there).  In terms of the concatenation
product, we have
\begin{align*}
  \F{n}&=\Fii ^n,
  \qquad
  \Fb{n} = \sum_{i=1}^{n}\xx^{n - i} \q^{-n + i} \Fii ^{i - 1} \Fi \Fii ^{n - i},
  \\[-4pt]
  \bX{n} &= \sum_{i=2}^{n}
  \xx^{n - i - 1} \q^{n - i - 1} (1 - \q^2) \Aint{i - 1}
  \Fii ^{i - 1} \Fi \Fii ^{n - i},
\end{align*}
and some lower-degree elements are given by
\begin{alignat*}{3}
  \Fb{1} &= \Fi,\quad&
  \Fb{2} &= \xx \q^{-1} \Fi \Fii + \Fii \Fi,\quad&
  \Fb{3} &=
  \xx^2 \q^{-2} \Fi \Fii \Fii + \xx \q^{-1} \Fii \Fi \Fii +  \Fii \Fii \Fi,
  \\
  &&\bX{2} &= \q^{-1} \xx^{-1} (1 - \q^2) \Fii \Fi,\quad
  &\bX{3} &=(1 - \q^2) \Fii \Fi \Fii + q^{-1} \xx^{-1}(1 - \q^4) \Fii \Fii \Fi,\\
  &&&&\bFb{3} &= (1 - \q^{-2}) \Fi \Fii \Fi.
\end{alignat*}

The shuffle multiplication table in the above basis is evaluated as
\begin{alignat*}{2}
  \FF{n} * \FF{m} &= \Abin{n + m}{n} \FF{n + m} ,\quad&
  \FF{n} * \Fb{m} &= \Abin{n + m - 1}{n} \Fb{n + m} ,\\
\FF{n} * \bX{m} &= \Abin{n + m - 2}{n} \BX{n + m},
  \quad&
  \FF{n} * \bFb{m} &= \Abin{n + m - 3}{n} \bFb{n + m} ,\\
  \Fb{n} * \FF{m} &=
  \q^{1 - m} \xx^{1 - m} \Abin{n + m - 2}{n - 1} \BX{n + m} + 
  \q^m \xx^{-m} \Abin{n + m - 1}{n - 1} \FB{n + m}
  \kern-160pt\\  
  \Fb{n} * \Fb{m} &=
  \xx^{2 - m}\q^{2 - m} \Abin{n + m - 3}{n - 1} \bFb{n + m}
  \quad\text{(in particular, $\Fb{n} * \Fb{1}=0$),}\kern-150pt\\
  \Fb{n} * \bX{m} &= -\q^{m - 1} \xx^{1 - m} \Abin{n\!+\!m\!-\!3}{n - 1}
  \BFB{n + m},
  \quad&
  \Fb{n} * \bFb{m} &= 0,\\
  \bX{n} * \FF{m} &= \q^{-m} \xx^{-m} \Abin{n + m - 2}{n - 2} \BX{n + m},
  \\
  \bX{n} * \Fb{m} &= \q^{1 - m} \xx^{1 - m} \Abin{n + m - 3}{n - 2}
  \BFB{n + m},
  \quad&
  \bX{n} * \bX{m} &= 0,\\
  \bX{n} * \bFb{m} &= 0,\quad&
  \bFb{n} * \FF{m} &= \xx^{-2 m} \Abin{n\!+\!m\!-\!3}{n - 3}\bFb{n + m},
\end{alignat*}
and all other products with $\bFb{n}$ vanish.

Although this is obvious from the general theory, we note explicitly
that the relations by which the quotient is taken
in~\eqref{eq:asymm-quotient} are now ``resolved''---hold
identically---due to the properties of the shuffle product; in
particular, the $*$-form of ``Serre relation''~\eqref{eq:asymm-Serre}
holds identically:
\begin{equation*}
  \xx^2 \Fi * \Fii * \Fii  - \xx (\q + \q^{-1}) \Fii * \Fi * \Fii
  + \Fii * \Fii * \Fi = 0.
\end{equation*}

\subsubsection{}\label{a:coproduct}With multiplication given by the
shuffle product, comultiplication is by deconcatenation (see
\cite{[STbr]}). \ For the above basis elements, their deconcatenation
can be calculated from the definition, for example,
\begin{multline*}
  \Delta\BFB{3}
  =(1 - \q^{-2})\bigl(\Co{\Fi \Fii \Fi}\one +  \Co{\Fi \Fii} 
   \Fi +  \Co{\Fi} \Fii \Fi +  \Co{\one} \Fi \Fii \Fi\bigr) \\
   = \Co{\BFB{3}}\one
   + \Co{\one} \BFB{3}
   + \Co{\BX{2}} \FB{1}
   - \xi\q^{-1}\Co{\FB{1}} \BX{2})
   + \xi^{-1}(\q - \q^{-1}) \Co{\FB{2}} \FB{1}.
 \end{multline*}
The result is
\begin{align*}
  \Delta\F{n} &=
  \sum_{i=0}^{n}\Co{\FF{i}} \F{n - i},
  \\
  \Delta\Fb{n} &= 
  \sum_{i=0}^{n - 1} \Co{\FF{i}}\Fb{n - i} + 
  \sum_{i=1}^{n} 
  \xx^{n - i} \q^{-n + i} \Co{\Fb{i}} \FF{n - i},
  \\
  \Delta\bX{n} &=  
  + \sum_{i=0}^{n - 2}\Co{\FF{i}}\BX{n - i}
  + \sum_{i=2}^{n} \xx^{n - i} \q^{n - i} \Co{\BX{i}}\FF{n - i}
  + \sum_{i=1}^{n - 1}\xx^{-1}\q^{2 n - 3} (\q^{-2 i} - 1)\Co{\FF{i}}\FB{n - i}
  ,\\
  \Delta\bFb{n} &=
  \sum_{i=0}^{n - 3}\Co{\FF{i}} \BFB{n - i}
  + \sum_{i=3}^{n}\xx^{2 n - 2 i} \Co{\BFB{i}} \FF{n - i}
  + \sum_{i=2}^{n - 1}\q^{n - i - 1} \xx^{n - i - 1} \Co{\BX{i}} \FB{n - i}\\
  &\quad{}
  - \sum_{i=1}^{n - 2}\xx^{n - i - 1} \q^{i - n + 1} \Co{\FB{i}} \BX{n - i}
  + \sum_{i=2}^{n - 1}\q^{n - i - 2} \xx^{n - i - 2} (\q^2 - 1)
  \Aint{i - 1} \Co{\FB{i}} \FB{n - i}.
\end{align*}

\subsubsection{}\label{a:antipode} The antipode is given by
``half-twist,'' the Matsumoto lift of the longest element in the
symmetric group~\cite{[STbr]}, which is evaluated for the above basis
elements as
\begin{align*}
  S(\FF{n}) &= (-1)^n \q^{n (n - 1)} \FF{n},\\
  S(\Fi ) &= -\Fi 
  ,\\
  S(\Fb{n}) &= (-1)^n \q^{(n - 4) (n - 1)} \Fb{n}
  + (-1)^n \xx \q^{(n - 4) (n - 1) + 1} \bF{n},\quad 2\leq n\leq p
  ,\\
  S(\bX{n}) &=
  (-1)^{n - 1} \q^{(n - 2) (n - 1)} \BX{n}
  + (-1)^n \xx^{-1}\q^{(n - 2) (n - 1) - 1} (1 - \q^2) \Aint{n-1}
  \FB{n}
  ,\\
  S(\bFb{n}) &= (-1)^{n + 1} \q^{(n-5) (n-2)} \bFb{n}.
\end{align*}

\subsection{Vertex operators and Yetter--Drinfeld $\Nich(X)$-modules}
For any pair of integers $\bar{r}$ and $\bar{s}$, we introduce a
one-dimensional braided vector space $Y^{\{\bar{r},\bar{s}\}}$, with a
fixed basis vector $V^{\{\bar{r},\bar{s}\}}$ (called a vertex
operator) such that
\begin{align*}
  \Psi: \Fi\tensor V^{\{\bar{r},\bar{s}\}}&\mapsto
  \q^{b}\,V^{\{\bar{r},\bar{s}\}}\tensor\Fi,\quad
  V^{\{\bar{r},\bar{s}\}}\tensor \Fi\mapsto
  \q^{b}\, \Fi\tensor V^{\{\bar{r},\bar{s}\}},
  \\
  \Psi: \Fii\tensor V^{\{\bar{r},\bar{s}\}}&\mapsto
  \q^{f}\,V^{\{\bar{r},\bar{s}\}}\tensor\Fii,\quad
  V^{\{\bar{r},\bar{s}\}}\tensor \Fii\mapsto
  \q^{f}\, \Fii\tensor V^{\{\bar{r},\bar{s}\}}.
\end{align*}
Each space $\Nich(X)\tensor
Y^{\{\bar{r}_1,\bar{s}_1\}}\tensor\Nich(X)\tensor
Y^{\{\bar{r}_2,\bar{s}_2\}}\tensor\dots\tensor \Nich(X)\tensor
Y^{\{\bar{r}_N,\bar{s}_N\}}$ is a Yetter--Drinfeld $\Nich(X)$-module
(an ``$N$-vertex'' Yetter--Drinfeld module in~\cite{[STbr]}) under the
left adjoint action and coaction given by deconcatenation up to the
first $Y$ space.  We here consider only one-vertex modules; the
$\Nich(X)$ action and coaction are then those in~\eqref{1-vertex};
coaction therefore reduces to just the comultiplication
in~\bref{a:coproduct}, and the adjoint action can be calculated using
the formulas in~\bref{sec:a-basis}--\bref{a:antipode}. \ The result is
\begin{align*}
  \Fi \adj \FF{n} \Vertex{\bar{r}, \bar{s}} &=
  \xx^{1 - n}\q^{1 - n}\BX{n + 1} \Vertex{\bar{r}, \bar{s}}
  + (1 - \q^{2 \bar{r} - 2 n}) \xx^{-n}\q^{n}\FB{n + 1} \Vertex{\bar{r}, \bar{s}}
  ,\\
  \Fi \adj \FB{1} \Vertex{\bar{r}, \bar{s}} &= 0,
  \\
  \Fi \adj \FB{n} \Vertex{\bar{r}, \bar{s}} &=
  \xx^{2 - n} \q^{2 - n} \BFB{n + 1} \Vertex{\bar{r}, \bar{s}},
  \quad 2\leq n\leq p,
  \\
  \Fi \adj \BX{n} \Vertex{\bar{r}, \bar{s}} &=
  -\xx^{1 - n}\q^{n-1} (1 - \q^{2 \bar{r} - 2 n + 2}) 
  \BFB{n + 1} \Vertex{\bar{r}, \bar{s}}
  ,\\
  \Fi \adj \BFB{n} \Vertex{\bar{r}, \bar{s}} &= 0
  \\
  \intertext{and}
  \Fii\adj \FF{n} \Vertex{\bar{r}, \bar{s}} &=
  \Aint{n+1} (1 - \q^{2 (n + \bar{s})}) \FF{n + 1} \Vertex{\bar{r}, \bar{s}},
  \\
  \Fii\adj \FB{n} \Vertex{\bar{r}, \bar{s}} &=
  \Aint{n} (1 - \q^{2 (n + \bar{s}- 1)}) \FB{n + 1} \Vertex{\bar{r}, \bar{s}}
  - \xx \q^{2 n + 2 \bar{s} - 3} \BX{n + 1} \Vertex{\bar{r}, \bar{s}}
  ,\\
  \Fii\adj \BX{n} \Vertex{\bar{r}, \bar{s}} &=
  \Aint{n-1}(1 - \q^{2 (n + \bar{s} - 2)}) \BX{n + 1} \Vertex{\bar{r}, \bar{s}}
  ,\\
  \Fii\adj \BFB{n} \Vertex{\bar{r}, \bar{s}} &=
  \Aint{n - 2} (1 - \q^{2 (n + \bar{s} - 3)}) \BFB{n + 1}\Vertex{\bar{r}, \bar{s}},
\end{align*}
where the range of $n$, unless specified explicitly, is in each case
as indicated in~\bref{sec:a-basis}.  We see that in the action on
$1$-vertex modules, $\bar{r}$ and $\bar{s}$ enter only modulo~$p$.  We
encounter this again in~\bref{V[n,m]}.

\section{\textbf{Nichols algebra: the ``symmetric'' case}}
We explicitly describe $\Nich(X)$ in the case corresponding to the
second line in~\eqref{eq:three-Q}; we take $X=\Xs$ to be a braided
vector space with basis $F_1$, $F_2$ and with diagonal braiding in
this basis specified by the braiding matrix
\begin{equation}\label{qij-symm}
  (q_{ij})=
  \begin{pmatrix}
    -1 & -\xx^{-1} \q \\
    -\xx\q & -1
  \end{pmatrix},
\end{equation}
where we allow $\xx$ to be any $2p$th root of unity, $\xx^{2p}=1$ (an
``inessential'' variable, eliminated by a twist
map~\cite{[AS-onthe]}).

\subsection{$\boldsymbol{\Nich(X)}$: a
  presentation}\label{symm-quotient}
The Nichols algebra $\Nich(X)$ associated with~\eqref{qij-symm} is the
quotient~\cite{[Ag-0804-standard]} (also see~\cite{[Hel]})
\begin{equation*}
  \Nich(X)=T(X)/ \bigl(F_1^{2},\
  (F_1 F_2)^p -  \xx^{-p}(F_2 F_1)^p,\
  F_2^{2}\bigr),
  \qquad \dim\Nich(X)=4p.
\end{equation*}

\subsection{$\boldsymbol{\Nich(X)}$ as a subalgebra}
Another description of $\Nich(X)$ in~\bref{symm-quotient}, which is
in fact the description in terms of screening operators, is not as a
quotient of but as a subspace in $T(X)$.  The product and coproduct in
$\Nich(X)$ are then the shuffle product and deconcatenation
(see~\cite{[Rosso-inv],[STbr]}).  We now describe these in some detail.

\subsubsection{The total symmetrizer map} Let $a$ and $b$ denote two
elements with respectively the same braiding as $F_1$ and $F_2$, but
without any algebraic constraints. Then the map $x\mapsto
\Bfac{\bullet}x$ by the total braided symmetrizer in each graded
component is as follows (see~\cite{[STbr]} for our conventions
on~$\Bfac{}$). \ Any noncommutative monomial containing $a^2$ or $b^2$
is mapped to zero, and the ``alternating'' monomials
\begin{equation*}
  \ab{n}=(ab)^n,\quad \ba{n}=(ba)^n,\quad \aba{n}=(ab)^n a,
  \quad \bab{n}=(ba)^n b  
\end{equation*}
 map as
\begin{align*}
  \ab{n}&\mapsto (1 - \q^2)^{n - 1} \Afac{n - 1}\, \bigl(\ab{n}
  - \xx^{-n} \q^n  \ba{n}\bigr),\\
  \ba{n}&\mapsto (1 - \q^2)^{n - 1} \Afac{n - 1}\, \bigl(\ba{n}
  - \xx^n \q^n  \ab{n}\bigr),\\
  \aba{n}&\mapsto (1 - \q^2)^{n} \Afac{n} \aba{n},\\
  \bab{n}&\mapsto (1 - \q^2)^{n} \Afac{n} \bab{n}.
\end{align*}
Hence, $\Bfac{2r+1}\aba{r}$ and $\Bfac{2r+1}\bab{r}$ are nonzero only
for $0\leq r\leq p-1$, and $\Bfac{2r}\ab{r}$ and $\Bfac{2r}\ba{r}$ are
nonzero only for $1\leq r\leq p$.  The image of $\Bfac{}$ contains
both $\ab{r}$ and $\ba{r}$ for $1\leq r\leq p-1$, but only the linear
combination $\ab{p} + \xx^{-p} \ba{p}$ for $r=p$.  
The algebra is therefore a linear span of the $4p$ elements
\begin{equation}\label{eq:ab-basis}
  \begin{alignedat}{2}
    &\one,
    \\
    &\ab{r},\quad 1\leq r\leq p-1,&\qquad& \ba{r},\quad 1\leq r\leq p-1,
    \\
    &\ab{p} + \xx^{-p} \ba{p},
    \\
    &\aba{r},\quad 0\leq r\leq p-1,&\qquad&\bab{r},\quad 0\leq r\leq p-1.
  \end{alignedat}
\end{equation}

\subsubsection{Shuffle product} 
The shuffle multiplication table of the above basis elements is as
follows:
\begin{alignat*}{2}
  \ab{r} * \ab{s} &= \Abin{r + s}{s} \ab{r + s},\\
  \ab{r} * \ba{s} &= \xx^s \q^s \Abin{r + s - 1}{s} \ab{r + s}
  + \xx^{-r} \q^r \Abin{r + s - 1}{s - 1} \ba{r + s},\kern-200pt\\
  \ab{r} * \aba{s} &=
  \Abin{r + s}{r} \aba{r + s},
  &\qquad
  \ab{r} * \bab{s} &= \xx^{-r} \q^r \Abin{r + s}{r} \bab{r + s},
  \\
  \ba{r} * 
   \ab{s} &= \xx^{-s} \q^s \Abin{r + s - 1}{s} \ba{r + s}
   + \xx^r \q^r \Abin{r + s - 1}{s - 1} \ab{r + s},\kern-200pt\\
  \ba{r} * \ba{s} &= \Abin{r + s}{s} \ba{r + s},
  &\qquad
  \ba{r} * \aba{s} &= \xx^r \q^r \Abin{r + s}{s} \aba{r + s},
  \\
  \ba{r} * \bab{s} &= \Abin{r + s}{s} \bab{r + s},
  &\qquad
  \aba{r} * \ab{s} &= \xx^{-s} \q^s \Abin{r + s}{s} \aba{r + s},
  \\
  \aba{r} * \ba{s} &= \Abin{r + s}{s} \aba{r + s},
  &\qquad
  \aba{r} * \aba{s} & = 0,
  \\
  \aba{r} * \bab{s} &= \Abin{r + s}{s} \ab{r + s + 1}
  - \xx^{-r - s - 1} \q^{r + s + 1} \Abin{r + s}{s} \ba{r + s + 1},\kern-200pt\\
  \bab{r} * \ab{s} &= \Abin{r + s}{s} \bab{r + s},
  &\qquad
  \bab{r} * \ba{s} &= \xx^s \q^s \Abin{r + s}{s} \bab{r + s},\kern-200pt\\
  \bab{r} * \aba{s} &= \Abin{r + s}{s} \ba{r + s + 1}
  - \xx^{r + s +1} \q^{r + s + 1} \Abin{r + s}{s} \ab{r + s + 1},\kern-200pt\\
  \bab{r} * \bab{s} &= 0.
\end{alignat*}
Due to the binomial coefficient vanishing, no elements outside the
ranges specified in~\eqref{eq:ab-basis} occur in the right-hand sides.
For $\ab{p} + \xx^{-p} \ba{p}$, strictly speaking, all products must
be listed separately, which is easy because all of them are zero;
these zero products are also reproduced by taking the appropriate
linear combinations of the above formulas, with due care; for example,
it follows that $(\ab{r} + \xx^{-r}
\ba{r})*\aba{s}=(1+\q^r)\Abin{r+s}{s}\aba{r+s}$, which vanishes at
$r=p$ for all $s\geq0$ already because $1+\q^p = 0$.

\subsubsection{}\label{sec:sDelta}The deconcatenation coproduct is
given by quite evident formulas,
\begin{align*}
  \Delta\ab{r}&=\one\tensor\ab{r}
   + a\tensor\bab{r-1}
   + \ab{1}\tensor\ab{r-1}
   + \dots
   + \aba{r-1}\tensor b + \ab{r}\tensor\one,
   \\
   \Delta\aba{r}&=\one\tensor\aba{r}
   + a\tensor\ba{r}
   + \ab{1}\tensor\aba{r-1}
   + \dots
   + \ab{r}\tensor a + \aba{r}\tensor\one,
\end{align*}
and similarly for $\ba{r}$ and $\bab{r}$.  The formula for
$\Delta(\ab{p} + \xx^{-p} \ba{p})$, once again, follows by extending
$\Delta\ab{r}$ and $\Delta\ba{r}$ to $r=p$ and combining them
appropriately.

\subsubsection{}
The antipode, given by the ``half-twist''~\cite{[STbr]} (the Matsumoto
lift of the longest element in the symmetric group), acts on the basis
monomials as
\begin{alignat*}{2}
  S(\ab{r}) &= (-1)^r \xx^{-r} \q^{r^2} \ba{r},\quad&
  S(\ba{r}) &= (-1)^r \xx^{r} \q^{r^2} \ab{r},
  \\
  S(\aba{r}) &= (-1)^{r + 1} \q^{r (r + 1)} \aba{r},\quad&
  S(\bab{r}) &= (-1)^{r + 1} \q^{r (r + 1)} \bab{r}
\end{alignat*}
(and $S(\ab{p} + \xx^{-p} \ba{p}) = \ab{p} + \xx^{-p} \ba{p}$).

\subsection{Vertex operators and Yetter--Drinfeld $\Nich(X)$-modules}
For any pair of integers $\bar{r}_1$ and $\bar{r}_2$, we introduce a
one-dimensional braided vector space $Y^{\{\bar{r}_1,\bar{r}_2\}}$, with a
fixed basis vector $V^{\{\bar{r}_1,\bar{r}_2\}}$ such that
\begin{align*}
  \Psi: a\tensor V^{\{\bar{r}_1,\bar{r}_2\}}&\mapsto
  \q^{\bar{r}_1}\,V^{\{\bar{r}_1,\bar{r}_2\}}\tensor a,\quad
  V^{\{\bar{r}_1,\bar{r}_2\}}\tensor a \mapsto
  \q^{\bar{r}_1}\, a\tensor V^{\{\bar{r}_1,\bar{r}_2\}},
  \\
  \Psi: b\tensor V^{\{\bar{r}_1,\bar{r}_2\}}&\mapsto
  \q^{\bar{r}_2}\,V^{\{\bar{r}_1,\bar{r}_2\}}\tensor b,\quad
  V^{\{\bar{r}_1,\bar{r}_2\}}\tensor b \mapsto
  \q^{\bar{r}_2}\, b\tensor V^{\{\bar{r}_1,\bar{r}_2\}}.
\end{align*}
Each space $\Nich(X)\tensor
Y^{\{{\bar{r}}_1^1,{\bar{r}}_2^1\}}\tensor\Nich(X)\tensor
Y^{\{{\bar{r}}_1^2,{\bar{r}}_2^2\}}\tensor\dots\tensor \Nich(X)\tensor
Y^{\{{\bar{r}}_1^N,{\bar{r}}_2^N\}}$ is a Yetter--Drinfeld
$\Nich(X)$-module (an ``$N$-vertex'' Yetter--Drinfeld
module)~\cite{[STbr]}.  We here consider only one-vertex modules; the
$\Nich(X)$ action and coaction are then those in~\eqref{1-vertex}.
The coaction is therefore literally the same as in~\bref{sec:sDelta},
and the (adjoint) action evaluates as
\begin{align*}
  a \adj \ab{n} V^{\{\bar{r}_1, \bar{r}_2\}}&=\xx^{-n} \q^n (1 - \q^{2 \bar{r}_1}) \aba{n} V^{\{\bar{r}_1, \bar{r}_2\}},\\
  a \adj\ba{n} V^{\{\bar{r}_1, \bar{r}_2\}}&=(1 - \q^{2 (\bar{r}_1 + n)}) \aba{n} V^{\{\bar{r}_1, \bar{r}_2\}},\\
  a \adj \aba{n} V^{\{\bar{r}_1, \bar{r}_2\}}&=0,
  \\
  a \adj \bab{n} V^{\{\bar{r}_1, \bar{r}_2\}}&=(1 - \q^{2 (\bar{r}_1 + n + 1)}) \ab{n + 1} V^{\{\bar{r}_1, \bar{r}_2\}} + 
  \xx^{-1 - n} \q^{n + 1} (\q^{2 \bar{r}_1} - 1) 
  \ba{n + 1} V^{\{\bar{r}_1, \bar{r}_2\}},
  \\
  b \adj \ab{n} V^{\{\bar{r}_1, \bar{r}_2\}} &=(1 - \q^{2 (\bar{r}_2 + n)}) \bab{n} V^{\{\bar{r}_1, \bar{r}_2\}},
  \\
  b \adj \ba{n} V^{\{\bar{r}_1, \bar{r}_2\}} &=\xx^n \q^n(1 - \q^{2 \bar{r}_2}) \bab{n} V^{\{\bar{r}_1, \bar{r}_2\}}
  \\
  b \adj \aba{n} V^{\{\bar{r}_1, \bar{r}_2\}} &= \xx^{n+1} \q^{n+1}(\q^{2 \bar{r}_2} - 1) 
  \ab{n + 1} 
  V^{\{\bar{r}_1, \bar{r}_2\}} + (1 - \q^{2 (\bar{r}_2 + n + 1)}) \ba{n + 1} 
  V^{\{\bar{r}_1, \bar{r}_2\}},
  \\
  b \adj \bab{n} V^{\{\bar{r}_1, \bar{r}_2\}} &= 0.
\end{align*}

\section{\textbf{From $\Nich(X)$ to extended chiral algebras: asymmetric
    case}}\label{CFTasymm}
We take two screening operators $\Fi$ and $\Fii$ with the respective
momenta $\alpha_1$ and $\alpha_2$ defined by the first line
in~\eqref{three-R}.  Then the corresponding braiding matrix is
$Q_{\asymm} = Q_{\asymm}(j,p)$ in~\eqref{qij-first}.
We then seek the kernel $\ker\Fi\bbigcap\ker\Fii$. \ In~\bref{aVir}
and~\bref{a-para}, we fix our conventions for scalar fields and, for
the convenience of the reader, recall the relevant points
from~\cite{[c-charge]}: a Virasoro algebra and parafermionic fields in
the kernel of the screenings; the nonlocal parafermionic fields are
converted into decent $\hSL2$ currents by introducing an ``auxiliary''
third scalar.  In~\bref{a-more} and~\bref{a-even-more}, we find more
fields in the kernel by looking at certain representations of this
$\hSL2$ and finally use the locality requirement and propose the
triplet--triplet and triplet--multiplet algebras, $\WWtritri$ and
$\WWtrimult$, which logarithmically extend the $\hSL2$ algebra.
In~\bref{a:hamred}, we briefly consider the Hamiltonian reduction of
the obtained $\WW$ algebras to the $\pone$ and $(p,p')$ triplet
algebras, in fact reproducing the constructions of the latter given
in~\cite{[FHST]} and~\cite{[FGST3]}.  In~\bref{a:W-modules}, we
outline the construction of some $\WW$-modules.

\subsection{Scalar fields and a Virasoro algebra}\label{aVir}
We introduce two scalar fields $\Aone(z)$ and $\Atwo(z)$ with the OPEs
determined by scalar products in the first line in~\eqref{three-R}:
\begin{alignat*}{2}
  \Aone(z)\,\Aone(w)&=\log(z-w),\qquad
  &\Aone(z)\,\Atwo(w)&=\bigl(-\ffrac{1}{p} - j\bigr)\log(z-w),
  \\
  \Atwo(z)\,\Atwo(w)&=\bigl(\ffrac{2}{p} + 2 j\bigr)\log(z-w).
\end{alignat*}
It is readily verified that with these OPEs, the kernel of our two
screenings
\begin{equation*}
  \Fi=\oint e^{\Aone}
  \quad
  \text{and}
  \quad
  \Fii=\oint e^{\Atwo}
\end{equation*}
contains the energy--momentum tensor\footnote{We tend to write $AB(z)$
  for the normal-ordered product of two fields $A(z)$ and $B(z)$, and
  $ABC(z)$, etc., for nested normal-ordered products $A(BC)(z)$, etc.
  The convention is not always obeyed, however; a notable case where
  it is violated is in expressions involving exponentials: for
  example, we write $(X(z)+\partial a \partial b(z)) e^{c(z)}$,
  whereas nested normal products are in fact understood in all cases,
  such as $(\partial a(\partial b\, e^{c}))(z)$.  Insisting on the
  rigorous writing would make many formulas incomprehensible.}
\begin{equation*}
  T(z)=
  -\ffrac{1}{k}\partial\Aone \partial\Aone(z)
  - \ffrac{1}{k} \partial\Aone \partial\Atwo(z)
  - \ffrac{1}{2 k (k +
    2)} \partial\Atwo\partial\Atwo(z)
  \\
  {}- \partial^2\Aone(z)
  - \fffrac{1}{2 (k + 2)}\partial^2\Atwo(z),
\end{equation*}
where we set
\begin{equation}\label{a-k2p}
  k \eqdef \ffrac{1}{p} + j - 2.
\end{equation}
This energy--momentum tensor represents the Virasoro algebra with the
central charge
\begin{equation}\label{eq:cc-asymm}
  c = \ffrac{3 k}{k + 2} - 1.
\end{equation}

It is useful to introduce $\omega_1,\omega_2\in\oC^2$ as ``fundamental
weights'' ($\omega_i\cdot\alpha_j = \delta_{i,j}$) with respect to
$\alpha_1$ and $\alpha_2$ defined by scalar products in the first line
in~\eqref{three-R}:
\begin{equation*}
  \omega_1 =
  \ffrac{1}{k} (-2 \alpha_1 - \alpha_2),\qquad
  \omega_2 =
  \ffrac{1}{k} \bigl(-\alpha_1 - \ffrac{1}{k + 2} \alpha_2\bigr).
\end{equation*}
It follows that
\begin{equation*}
  \omega_1.\omega_1=
  -\ffrac{2}{k},
  \quad
  \omega_1.\omega_2=
  -\ffrac{1}{k},
  \quad
  \omega_2.\omega_2=
  -\ffrac{1}{k(k + 2)}.
\end{equation*}
Anticipating the appearance of a third scalar, we also consider a
three-dimensional space spanned by $\alpha_1$, $\alpha_2$, and
$\alpha_3$, where $\alpha_3$ has zero scalar products with $\alpha_1$
and $\alpha_2$ and is normalized by the condition read off
from~\eqref{localize-2}, $\alpha_3.\alpha_3 = 2 k$.  Then the
corresponding ``fundamental weight'' is $\omega_3=
\fffrac{1}{2 k} \alpha_3$.

We slightly abuse the notation and write $\omega_i(z)$ for the
appropriate linear combinations of the fields,
\begin{equation*}
  \omega_1(z) =
  \ffrac{1}{k} (-2 \Aone(z) - \Atwo(z)),\quad
  \omega_2(z) =
  \ffrac{1}{k} \bigl(-\Aone(z) - \ffrac{1}{k + 2} \Atwo(z)\bigr),
  \quad
  \omega_3(z)= \ffrac{1}{2 k}\Athree(z).
\end{equation*}

\subsection{Parafermionic fields in
  $\bker\myboldsymbol{\Fi{}\bigcap{}}\bker\myboldsymbol{\Fii}$ and the
  $\myboldsymbol{\hSL2_k}$ algebra}\label{a-para} The form of the
central charge in~\eqref{eq:cc-asymm} immediately suggests that the
kernel must contain generators of the $\hSL2/u(1)$ coset; indeed,
these are given by
\begin{equation}\label{para:asymm}
  \begin{aligned}
    \jplus(z)&=
    e^{\omega_1(z)}
    \\
    \jminus(z)&=
    -\bigl(\partial\Aone \partial\Aone(z) 
    + \partial\Aone \partial\Atwo(z)
    +(k + 1) \partial^2\Aone(z)\bigr)
    e^{-\omega_1(z)}
  \end{aligned}
\end{equation}

\subsubsection{A long screening} The fields $T(z)$ and $j^{\pm}(z)$
are in the kernel not only of the two ``short'' screenings $\Fi$ and
$\Fii$ but also of the ``long'' screening
\begin{equation*}
  \Ebeta=\oint e^{-\frac{1}{k + 2}\Atwo}
  =\oint e^{-\frac{p}{j p + 1}\Atwo}.
\end{equation*}
The long screening is readily verified to commute with both short
screenings:
\begin{equation*}
  [\Ebeta,\Fi]=0,\qquad
  [\Ebeta,\Fii]=0.
\end{equation*}

\subsubsection{The $\myboldsymbol{\hSL2}$ currents}
The fields $j^{\pm}(z)$ are nonlocal with respect to one another
(their OPEs contain noninteger powers $(z-w)^{\pm 2 p/((j-2) p +
  1)}$), and in fact represent a coset theory $\hSL2/u(1)$.  To deal
with local fields, we introduce a third, auxiliary scalar $\Athree(z)$
with OPE~\eqref{localize-2} and construct the $\hSL2_k$ currents (with
our conventions given in~\bref{app:sl2-conv})
\begin{align}
  \Jplus(z)&=
  \jplus(z) e^{\frac{1}{k}\Athree(z)},
  \notag
  \\
  \label{sl2-asymm}
  \Jnaught(z)&=\half\partial\Athree(z),
  \\
  \Jminus(z)&=
  \jminus(z) e^{-\frac{1}{k}\Athree(z)}. \notag
\end{align}
The associated Sugawara energy--momentum tensor~\eqref{eq:Sug} is then
evaluated as
\begin{equation*}
  T_{\text{Sug}}(z)
  = T(z) +\ffrac{1}{4k}\partial\Athree\partial\Athree(z).
\end{equation*}

\subsubsection{Orthogonalizing the scalar fields}\label{ortho} It is
useful to pass to a triplet of pairwise OPE-orthog\-onal scalar
fields, $(\Aone(z),\Atwo(z),\Athree(z)) \to
(\Afour(z),\Aii(z),\Athree(z))$, where
\begin{align*}
  \Afour(z) &= -2\Aone(z) - \Atwo(z),
  \\
  \Aii(z)&=\Atwo(z) 
\end{align*} %
are mutually orthogonal (have a regular OPE) and $\Afour(z)$ is
normalized as
\begin{equation*}
  \Afour(z)\,\Afour(w) = -2 k \log(z-w).
\end{equation*}
We keep the notation $\omega_i(z)$ introduced in~\bref{aVir}; we then
have the OPEs
\begin{equation*}
  \Afour(z)\,\omega_1(w) = -2\log(z-w),
  \qquad
  \Afour(z)\,\omega_2(w) = -\log(z-w).
\end{equation*}
Whenever the context requires this, we understand $\omega_1(z)$ and
$\omega_2(z)$ to be reexpressed in terms of the new fields, as
\begin{equation*}
  \omega_1(z) = \ffrac{1}{k}\Afour(z),\qquad
  \omega_2(z) = \ffrac{1}{2k}\Afour(z)+\ffrac{1}{2(k+2)}\Aii(z).
\end{equation*}

Then the $\hSL2_k$ currents in~\eqref{sl2-asymm} take the
form~\cite{[S-inv]}
\begin{equation}\label{through-matter}
  \begin{aligned}
    \Jplus(z) &=
    e^{\omega_1(z) + 2\omega_3(z)}
    ,\\
    \Jnaught(z) &= \half\partial\Athree(z),\\
    \Jminus(z) &= \bigl((k + 2) \Tm(z) - 
    \ffrac{1}{4}\partial\Afour\partial\Afour(z)
    + \ffrac{k + 1}{2} \partial^2\Afour(z)\bigr)
    e^{-\omega_1(z) - 2\omega_3(z)},
  \end{aligned}
\end{equation}
where
\begin{equation}\label{Tm-expressed}
  \Tm(z) \eqdef \ffrac{1}{4 (k + 2)}\partial\Aii\partial\Aii(z) + 
  \ffrac{k + 1}{2 (k + 2)} \partial^2\Aii(z).
\end{equation}
This $\Tm(z)$ is an energy--momentum tensor with central charge
\begin{equation}\label{mattercc}
  \cm = 13 - 6(k + 2) - \ffrac{6}{k + 2},
\end{equation}
i.e., it satisfies the OPE
\begin{equation*}
  \Tm(z)\,\Tm(w) = \ffrac{\cm/2}{(z-w)^4}
  + \ffrac{2\Tm(w)}{(z-w)^2}
  + \ffrac{\partial\Tm(w)}{z-w}
\end{equation*}
encoding a Virasoro algebra.  We note for the future use that primary
fields of this Virasoro algebra and their conformal dimensions
are
\begin{equation}\label{matter-fields-dim}
  V^{\text{m}}_{r,s}(z)
  =e^{\left(\half(1 - r) + \frac{s - 1}{2 (k + 2)}\right)\Aii(z)},\quad
  \Deltamatter_{r,s}=
  \fffrac{r^2 - 1}{4}(k + 2) + \fffrac{s^2 - 1}{4(k + 2)}
  + \fffrac{1}{2}(1 - s r).
\end{equation}

We also note that the Sugawara energy--momentum tensor is now
reexpressed as a sum of three energy--momentum tensors,
\begin{equation*}
  T_{\text{Sug}}(z) = \Tm(z)
  + \bigl(-\ffrac{1}{4 k}\partial\Afour\partial\Afour(z)
  + \half \partial^2\Afour(z)\bigr)
  + \ffrac{1}{4 k}\partial\Athree\partial\Athree(z),
\end{equation*}
with the respective central charges $\cm$, \ $6 k + 1$, and $1$.

\subsubsection{}\label{two-standpoints} The fact that the $\Aii(z)$
field enters the $\hSL2$ currents in~\eqref{through-matter} only
through $\Tm(z)$ \cite{[S-inv]} can be viewed, depending on one's
taste, as either i)~a technicality or ii)~an important structural
piece; we switch between the two standpoints at will.
\begin{itemize}

\item[i)] We can of course assume that $\Tm(z)$ is expressed through a
  free scalar, as in~\eqref{Tm-expressed}, and continue working with
  the three scalars $\Afour(z)$, $\Aii(z)$, and $\Athree(z)$; an
  enveloping algebra of $\hSL2_k$ is then selected from differential
  polynomials in $\partial\Afour(z)$, $\partial\Aii(z)$,
  $\partial\Athree(z)$, and $e^{\pm(\frac{1}{k} \Athree(z) +
    \frac{1}{k}\Afour(z))}$ by taking the kernel of both screenings
  \begin{equation*}
    \Fi = \oint e^{-\half\Aii - \half\Afour},\qquad
    \Fii = \oint e^{\Aii}.
  \end{equation*}

\item[ii)] Equivalently, and more interestingly, we can recall that
  the kernel of the single screening $\Fii$ selects differential
  polynomials in $\Tm(z)$ (the enveloping of the Virasoro algebra)
  from the algebra of differential polynomials in
  $\partial\Aii(z)$. We can therefore ``use up'' the $\Fii$ screening
  to forget about the $\Aii(z)$ scalar and deal with the $\hSL2_k$
  currents  in~\eqref{through-matter}\pagebreak[3] expressed in terms
  of two scalars $\Afour(z)$ and $\Athree(z)$ and an ``abstract''
  energy--momentum tensor---\textit{not} a free scalar---with central
  charge~\eqref{mattercc}.  The remaining screening,
  \begin{equation}\label{F2-matter}
    \Fi = \oint V^{\text{m}}_{2,1} e^{-\half\Afour},
  \end{equation}
  then serves to select the enveloping algebra of $\hSL2_k$ from the
  algebra of differential polynomials in $\partial\Afour(z)$,
  $\partial\Athree(z)$, $\Tm(z)$, and $e^{\pm(\frac{1}{k} \Athree(z) +
    \frac{1}{k}\Afour(z))}$. \ In~\eqref{F2-matter},
  $V^{\text{m}}_{2,1}(z)$ is the ``21'' primary field of the Virasoro
  algebra with central charge $\cm$: it is defined by the
  OPE
  \begin{equation*}
    \Tm(z)\,V^{\text{m}}_{2,1}(w)
    = \ffrac{\Deltamatter_{2,1} V^{\text{m}}_{2,1}(w)}{(z-w)^2}
    + \ffrac{\partial V^{\text{m}}_{2,1}(w)}{z-w},\qquad
    \Deltamatter_{2,1}= \ffrac{1}{4}(3k+4),
  \end{equation*}
  and the differential equation\footnote{This equation is needed, in
    particular, in verifying that \eqref{F2-matter} \textit{is} a
    screening for the $\hSL2$ algebra in~\eqref{through-matter}.}
  \begin{equation*}
    \partial^2 V^{\text{m}}_{2,1}(z) - (k + 2) \Tm
    V^{\text{m}}_{2,1}(z) = 0.
  \end{equation*}
\end{itemize}
In what follows, we speak of the $\Aii$ sector as the \textit{matter}
theory; the name is motivated by the relation to Hamiltonian
reduction, as we see below.

\subsection{More of
  $\bker\myboldsymbol{\Fi{}\bigcap{}}\bker\myboldsymbol{\Fii}$}
\label{a-more}
Another piece of the kernel is easy to find.  It contains the fields
\begin{equation}\label{AsymmF}
  \AsymmF_h(z)=
  e^{2 h (\omega_2(z) + \omega_3(z))}
\end{equation}
for any $h=\half,1,\frac{3}{2},2,\dots$.  Each $\AsymmF_h(z)$ is an
$\hSL2$ primary,
\begin{gather*}
  \Jminus_{1}\AsymmF_h(z) = 0,\quad \Jnaught_0 \AsymmF_h(z) = h
  \AsymmF_h(z),\quad \Jplus_0 \AsymmF_h(z) = 0,
\end{gather*}
and generates a ``horizontal'' $(2h+1)$-plet under the action of the
zero-mode $s\ell(2)$ algebra:
\begin{gather*}  
  (\Jminus_0)^{2h+1}\AsymmF_h(z) = 0,\qquad
  (\Jminus_0)^{2h}\AsymmF_h(z)\neq0.
\end{gather*}
Of course, the entire $\hSL2$ module generated from $\AsymmF_h(z)$ is
in the kernel.  

If the matter theory is singled out as explained above, we reexpress
$\AsymmF_h(z)$ as
\begin{equation*}
  \AsymmF_h(z) =
  V^{\text{m}}_{1,2h+1}(z)e^{\frac{h}{k}\Afour(z) + \frac{h}{k}\Athree(z)},
\end{equation*}
where $V^{\text{m}}_{1,s}(z)$ are $\Tm(z)$-primary fields of dimension
$\Deltamatter_{1,s}$ (see~\eqref{matter-fields-dim}).
The Sugawara dimension of $ \AsymmF_h(z)$ is of course
\begin{equation}\label{Delta-gen}
  \Delta_j(h)= \ffrac{h (h + 1)}{k + 2} = \ffrac{h (h + 1) p}{j p + 1}.
\end{equation}

\subsection{Even more of
  $\bker\myboldsymbol{\Fi{}\bigcap{}}\bker\myboldsymbol{\Fii}$ and the
  extended algebras}\label{a-even-more}
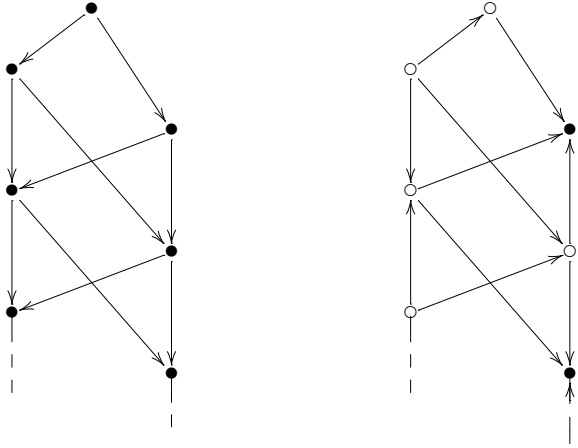
\begin{figure}[tbh]
  \centering
  \begin{equation*}
    \xymatrix@C=2pt@R=14pt{
      &&&\Bullet="Vtop"&&&
      &&&&&&&&&
      &&&\Circ="top"&&&
      \\
      \Bullet="Vcos1"&&&&&&
      &&&&&&&&&
      \Circ="cos1"&&&&&&
      \\
      &&&&&&\Bullet="Vsing1"
      &&&&&&&&&
      &&&&&&\Bullet="sing1"
      \\
      \Bullet="Vcos2"&&&&&&
      &&&&&&&&&
      \Circ="cos2"&&&&&&
      \\
      &&&&&&\Bullet="Vsing2"
      &&&&&&&&&
      &&&&&&\Circ="sing2"
      \\
      \Bullet="Vcos3"&&&&&&
      &&&&&&&&&
      \Circ="cos3"&&&&&&
      \\
      &&&&&&\Bullet="Vsing3"
      &&&&&&&&&
      &&&&&&\Bullet="sing3"
      \ar"cos1";"top"
      \ar"top";"sing1"
      \ar"cos1";"cos2"
      \ar"cos3";"cos2"
      \ar"sing2";"sing1"
      \ar"sing2";"sing3"
      \ar"Vtop";"Vcos1"
      \ar"Vtop";"Vsing1"
      \ar"Vcos1";"Vcos2"
      \ar"Vcos2";"Vcos3"
      \ar"Vsing1";"Vsing2"
      \ar"Vsing2";"Vsing3"
      \ar"cos1";"sing2"
      \ar"cos2";"sing3"
      \ar"cos2";"sing1"
      \ar"cos3";"sing2"
      \ar"Vcos1";"Vsing2"
      \ar"Vcos2";"Vsing3"
      \ar"Vsing1";"Vcos2"
      \ar"Vsing2";"Vcos3"
      \ar@{{}{--}{-}}"cos3"+<0pt,-30pt>;"cos3"
      \ar"sing3"+<0pt,-10pt>;"sing3"
      \ar@{{-}{--}{-}}"sing3"+<0pt,-24pt>;"sing3"
      \ar@{{}{--}{-}}"Vcos3"+<0pt,-30pt>;"Vcos3"
      \ar@{{}{--}{-}}"Vsing3"+<0pt,-20pt>;"Vsing3"
    }
  \end{equation*}
  \caption{\small\textsc{Left:} Embedding diagram of an $\hSL2$ Verma
    module.  The highest-weight vector is at the top.  Arrows are
    drawn toward singular vectors (submodules). \ \textsc{Right:} The
    corresponding Wakimoto module.  The same subquotients are
    ``glued'' to one another in a different manner.  Black dots show
    the socle of the module.}
  \label{fig:zzz}
\end{figure}

We temporarily fix a positive integer or half-integer $h$.  It is easy
to verify that $\AsymmF_{-h}(z)$ is in $\ker\Fi$, but not in
$\ker\Fii$.  The intersection of the two kernels is to be found deeper
in the Wakimoto-type free-field module~\cite{[W],[FF]} associated with
$\AsymmF_{-h}(z)$.\pagebreak[3] The actual picture depends on the
value of~$j$.  We consider the special case $j=0$ separately and then
discuss the cases $j>0$.  When we speak of nonvanishing and vanishing
singular vectors in what follows, we refer to the typical picture in
Fig.~\ref{fig:zzz}.
The Verma-module embedding pattern on the left changes in Wakimoto
modules to the one on the right.  In particular, the ``reversal'' of
an arrow leading from the top means that the corresponding singular
vector vanishes.

\subsubsection{The case $\myboldsymbol{j=0}$ and a triplet--triplet
  algebra}\label{a:j=0}
In the $\hSL2_k$ Verma module whose highest-weight vector has the same
charge (eigenvalue of $\Jnaught_0$) and dimension as those of
$\AsymmF_{-h}(z)$, there are two basic singular vectors, which happen
to lie on the same level (see Fig.~\ref{fig:a-TRIPLET}; our notation
for $\hSL2$ singular vectors is explained in
Appendix~\ref{app:sl2}):
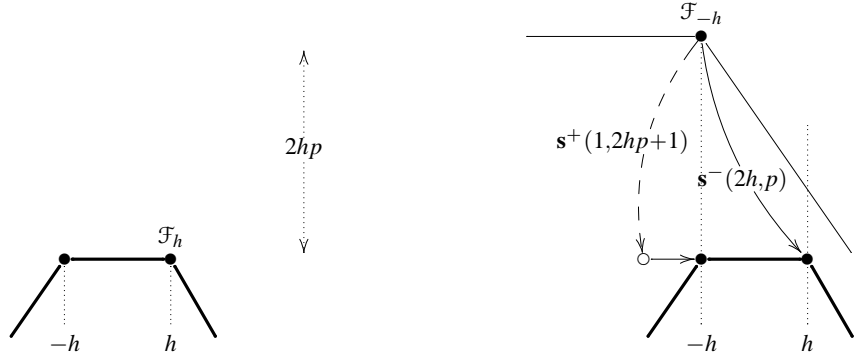
\begin{figure}[tb]
  \centering
  \vspace*{-1.5\baselineskip}
  \begin{equation*} 
    \xymatrix@C8pt@R20pt{
      &&&&&&&&&&&&&&&&&&&\\
      &&&&&&&\ar@{{<}{.}{>}}|{2 h p}[ddd]&&&&&&&&&
      *{\bullet}\ar@{-}[dddrrr]\ar@{-}[llll]\ar@{{.}{.}{.}}[dddd]
      \ar@{}[u]|(.3){\AsymmF_{-h}}
      \ar@{{-}{--}{>}}
      @/_10pt/|{\MFFplus{1,2hp+1}\quad\ \ }[dddl]
      \ar
      @/_8pt/|(.6){\MFFminus{2h,p}\ }[dddrr]&&&
      \\
      &&&&&&&&&&&&&&&&&&\ar@{{.}{.}{.}}[ddd]&\\
      &&&&&&&&&&&&&&&&&&&\\
      &&*{\bullet}\ar@*{[|(3)]}@{-}[rr]\ar@*{[|(3)]}@{-}[dl]&
      &*{\bullet}\ar@*{[|(3)]}@{-}[dr]\ar@{}[u]|(.3){\AsymmF_h}
      &&&&&&&&&&
      &*{\circ}\ar[r]&*{\bullet}\ar@*{[|(3)]}@{-}[rr]\ar@*{[|(3)]}@{-}[dl]&
      &*{\bullet}\ar@*{[|(3)]}@{-}[dr]&\\
      &&{\scriptstyle -h}\ar@{{.}{.}{.}}[u]&
      &{\scriptstyle h}\ar@{{.}{.}{.}}[u]&&&&&&&&&&&
      &{\scriptstyle-h}&&{\scriptstyle h}&\\
    }
  \end{equation*}
  \caption{\small\textsc{Left}: the $\hSL2$ primary $\AsymmF_h(z)$,
    Eq.~\eqref{AsymmF}, is an element of a $(2h+1)$-plet under the
    action of the zero-mode $s\ell(2)$ algebra.  \textsc{Right}:
    $\AsymmF_{-h}(z)$ for $h>0$ and some structures in the associated
    module.  The open dot shows a ``cosingular'' vector in the grade
    where the $\MFFplus{1,2hp+1}$ singular vector vanishes.  The
    horizontal arrow is a map by $\Jplus_0$.  The right black dot is a
    nonvanishing singular vector, from which the zero-mode $s\ell(2)$
    algebra generates a $(2h+1)$-plet, of the same Sugawara dimension
    as in the left diagram.  The dimensions of $\AsymmF_{h}$ and
    $\AsymmF_{-h}$ compare as $\frac{h(h+1)}{k+2} -
    \frac{-h(-h+1)}{k+2} =2hp$.}
  \label{fig:a-TRIPLET}
\end{figure}%
\begin{itemize}
\item $\MFFplus{1,2hp+1}$ with $(\text{charge},\text{level})$ relative
  to those of the highest-weight vector equal to $(-1,2hp)$, and

\item $\MFFminus{2h,p}$ with the relative charge and level $(2h,2hp)$.
\end{itemize}
In the free-field module that we actually have,
$\MFFplus{1,2hp+1}\AsymmF_{-h}(z)$ vanishes; this is shown in
Fig.~\ref{fig:a-TRIPLET} with a dashed line. On the other hand, the
$\MFFminus{2h,p}$ singular vector evaluated on $\AsymmF_{-h}(z)$ does
not vanish and \textit{is in the kernel of both screenings}.  The
zero-mode $s\ell(2)$ algebra produces a $(2h+1)$-plet from it,
terminating in the grade next to the one with the vanishing singular
vector.

The value $j=0$ is singled out by the fact that the long screening is
then a generator of a Lie algebra $s\ell(2)$ rather than of a quantum
$s\ell(2)$ group, as for other integer $j$.  Mapping $\AsymmF_{h}(z)$
by the long screening produces just the $\MFFminus{2h,p}$ singular
vector:
\begin{align}\label{Elong-maps}
  (\Ebeta)^{2h}\AsymmF_{h}(z) &= \MFFminus{2h,p}\AsymmF_{-h}(z).
\end{align}
Because the long screening is in the ``matter'' sector,
Eq.~\eqref{Elong-maps} can also be written as
\begin{equation*}
  (\Ebeta)^{2h}\AsymmF_{h}(z)
  =\mathscr{S}V^{\text{m}}_{2p-1,2h+1}(z)\,
  e^{h (\omega_1(z) + 2 \omega_3(z))},
\end{equation*}
where $\mathscr{S}V^{\text{m}}_{2p-1,2h+1}(z)$ is the corresponding
\textit{Virasoro} singular vector evaluated on (the field
corresponding to) the Virasoro primary state
$V^{\text{m}}_{2p-1,2h+1}(z)$, which occurs here because
\begin{equation*}
  \AsymmF_{-h}(z) =
  V^{\text{m}}_{2p-1,2h+1}(z)
  e^{-h (\omega_1(z) + 2 \omega_3(z))}.
\end{equation*}
We also note the \textit{matter} dimension of the fields in the two
$(2h+1)$-plets:
\begin{equation}\label{delta(h)}
  \Deltamatter_{1, 2 h + 1} = h (p(h+1) - 1).
\end{equation}

Continuing with the embedding diagrams of Wakimoto-type modules allows
describing all of the socle (the black dots in
Fig.~\ref{fig:zzz},~right, plus similar dots in other Wakimoto
modules, which increase in number as we down the diagram), but we stop
here because our main task now is to propose generators
(minimal-dimension fields) of the maximum local algebra in the kernel.

The mutual (non)locality of $\AsymmF_{h}(z)$ and $\AsymmF_{h'}(w)$ is
measured by the powers $(z-w)^{2 h h' p}$ in their operator product.
We choose $h$ such that these exponents be integer for all, integer
and half-integer, $h'$.  This means taking integer $h$, and the local
algebra generators are those with the smallest positive integer~$h=1$,
$\AsymmF_1(z)$ and $\MFFminus{2,p}\AsymmF_{-1}(z)$.  These fields are
then the leftmost and the rightmost operators in a triplet under the
action of $\Ebeta$.  In addition, each of these fields is the
rightmost element in a triplet with respect to the action of the
zero-mode~$s\ell(2)$.

For each integer $p\geq2$,
\begin{it}%
  we propose the algebra $\WWtritri$ generated by
  the fields
\begin{equation}\label{a:the-W-j=0}
  \begin{aligned}
    \W^+(z) &= \AsymmF_1(z),
    \\
    \W^0(z) &= \Ebeta\W^+(z) 
    \\
    \W^-(z) &=\Ebeta\W^0(z) = \MFFminus{2,p}\AsymmF_{-1}(z),
  \end{aligned}
\end{equation}
together with the corresponding $s\ell(2)$ triplets
$(\Jminus_0)^i\W^{\pm,0}(z)$, $i=0,1,2$, as a ``logarithmic''
extension of the $\hSL2$ algebra at the level $k=\frac{1}{p} - 2$.
\end{it}%
\ To be more explicit, we recall that the $\AsymmF_{\pm1}(z)$ are here
given by
\begin{equation*}
  \begin{aligned}
    \AsymmF_1(z)
    &=
    e^{2 \omega_2(z) + 2 \omega_3(z)}
    =V^{\text{m}}_{1,3}(z)
    e^{\omega_1(z) + 2 \omega_3(z)},
    \\
    \AsymmF_{-1}(z)
    &=
    e^{-2 \omega_2(z) - 2 \omega_3(z)}
    =V^{\text{m}}_{2p-1,3}(z)
    e^{-\omega_1(z) - 2 \omega_3(z)},
  \end{aligned}
  \qquad k=\fffrac{1}{p} - 2.
\end{equation*}
Conjecturally, \textit{$\WWtritri$ contains all mutually local fields
  in $\ker\Fi\bbigcap\ker\Fii$} (in particular, the $\AsymmF_{n}(z)$
with integer $n\geq 2$ and their images $\Ebeta^m\AsymmF_{n}(z)$,
$1\leq m\leq 2n$, under the long screening).

Each field $\W^a(z)=\W^{+,0,-}(z)$ also belongs to a triplet
$((\Jminus_0)^2\W^a(z),\Jminus_0\W^a(z),\W^a(z))$ under the zero-mode
$s\ell(2)$ algebra, as in Fig.~\ref{fig:a-TRIPLET} (where we now set
$h=1$).  In particular,
\begin{multline}\label{fully-op}
  (\Jminus_0)^2\W^{+}(z)=\\
  \Bigl(
  \ffrac{1}{2}  \partial\Aii\partial\Aii(z)
  -\partial^2\Aii(z)
  + \partial\Aii\partial\Afour(z)
  + \ffrac{1}{2}  \partial\Afour(z)\partial\Afour(z)
  - \partial^2\Afour(z) \Bigr)
  e^{-2 \omega_1(z) + 2 \omega_2(z) - 2\omega_3(z)}.
\end{multline}

The fields $\W^+(z)$, $\W^0(z)$, $\W^-(z)$ (and the entire
``zero-mode'' triplets) have the Sugawara dimension
(see~\eqref{Delta-gen})
\begin{equation*}
  \Delta_0(1)=2p.
\end{equation*}

An example of the triplet--triplet algebra generators is given
in~\bref{app:example-a}.\enlargethispage{\baselineskip}

\medskip

\noindent\textbf{Remark: a doublet.}\ \ 
We note that setting $h=\half$ instead of $h=1$ in~\eqref{Elong-maps}
yields a \textit{doublet} of dimension-$\frac{3p}{4}$ fields
$\AsymmF_{\half}(z)$ and $\Ebeta\AsymmF_{\half}(z) =
\MFFminus{1,p}\AsymmF_{-\half}(z)$ (each of which also belongs to a
zero-mode $s\ell(2)$ doublet); they can be regarded as generators of
an ``almost local'' \textit{doublet} algebra---an analogue of a pair
of fermions (derivatives of the symplectic fermions $\psi^{\pm}(z)$)
well known from the logarithmic $(p=2,1)$ ``matter''
models.\footnote{Hamiltonian reduction (see~\bref{HR:j=0} below) of
  $\AsymmF_{\half}(z)$ and $\Ebeta\AsymmF_{\half}(z)$ gives a doublet
  of ``matter'' fields of dimension $\frac{3p-2}{4}$, which
  \textit{are} the $\partial\psi^{\pm}(z)$ for $p=2$.}

\subsubsection{Cases $\myboldsymbol{j\geq 1}$ and triplet--multiplet
  algebras}\label{a:j>0}
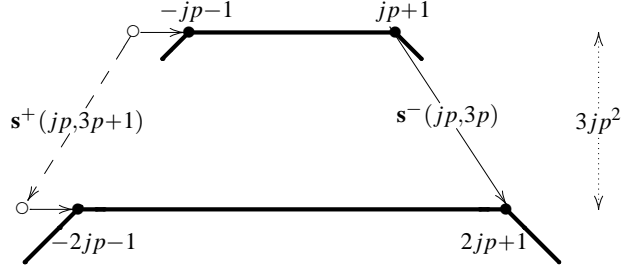
\begin{figure}[tb]
  \centering
  \begin{equation*}
    \xymatrix@C10pt@R50pt{
      &&&\Circ="cosing"&\iverma="L"&&&&\verma="R"&&&&\nothing="up"
      \\
      &\Circ="cosing2"&\ivermaiv="L2"&&&&&&&&\vermaiv="R2"&&\nothing="down"
      \ar@*{[|(4)]}@{-}"L"+<0pt,0pt>;"R"+<0pt,0pt>
      \ar@*{[|(4)]}@{-}"R"+<0pt,0pt>;"R"+<10pt,-10pt>
      \ar@*{[|(4)]}@{-}"L"+<0pt,0pt>;"L"-<10pt,10pt>
      \ar@*{[|(4)]}@{-}"L2";"R2"
      \ar@{{-}{--}{>}}|{\MFFplus{j p, 3 p + 1}}"cosing";"cosing2"
      \ar|{\MFFminus{j p, 3 p}}"R"+<-2pt,2pt>;"R2"+<0pt,2pt>
      \ar"cosing";"L"+<-2pt,0pt>
      \ar"cosing2";"L2"+<-2pt,0pt>
      \ar@{{<}{.}{>}}|{3 j p^2}"up";"down"
      \ar@*{[|(4)]}@{-}"R2";"R2"+<20pt,-20pt>
      \ar@*{[|(4)]}@{-}"L2";"L2"-<20pt,20pt>
      \ar@{}|{\scriptstyle-j p - 1}"L";"L"+<6pt,16pt>
      \ar@{}|{\scriptstyle j p + 1}"R";"R"+<6pt,16pt>
      \ar@{}|{\scriptstyle\ \  -2 j p - 1}"L2";"L2"+<6pt,-26pt>
      \ar@{}|{\scriptstyle 2 j p + 1\quad\ }"R2";"R2"+<6pt,-26pt>
    }
  \end{equation*}
  \caption{\small Relevant structures in the Wakimoto-type $\hSL2$
    module associated with the field $\AsymmF_{j p + 1}(z)$, which
    sits at the top right corner.  Vanishing singular vectors are
    shown with open dots.  The charges (eigenvalues of $\Jnaught_0$)
    of the states are indicated.  Two arrows show singular vectors,
    $\MFFname^{\pm}$, in the corresponding Verma module; the right
    singular vector is nonvanishing in our free-field realization.
    The dotted arrow shows the relative level (difference of Sugawara
    dimensions) of the two floors.}
  \label{fig:a:left:j>0}
\end{figure}

For $j\geq1$, we repeat the construction
in~\bref{a:j=0} mutatis mutandis,\pagebreak[3] noting from the start
that with $k = j - 2 + \fffrac{1}{p}$, mutual (non)locality of vertex
operators is measured by the powers $(z-w)^{2h h' p/(p j + 1)}$ in
their OPEs.  In seeking the local algebra in the kernel of the
screenings, we therefore start with the $\hSL2$ modules generated from
$\AsymmF_{\pm h}(z)$ with
\begin{equation}\label{eq:h(j)}
  h = j p + 1.
\end{equation}
In Fig.~\ref{fig:a:left:j>0}, we represent $\AsymmF_{j p + 1}(z)$ with
the top right corner.
It has a vanishing singular vector $\MFFplus{2 j p + 3, 1}$ (the top
left open circle) and the nonvanishing singular vector
\begin{align}\label{jWplus}
  \jW^+(z)&=\MFFminus{j p, 3 p}\AsymmF_{j p + 1}(z),
  \\
  \intertext{from which (due to another singular vector vanishing) the
    zero-mode $s\ell(2)$ algebra generates a $(4 j p + 3)$-plet.  It
    is easy to see that $\jW^+(z)$ has the structure}
  \jW^+(z)&=\mathscr{P}^+(z)\,
  e^{j p \omega_1(z) + 2 (j p + 1) \omega_2(z) + 2 (2 j p + 1)\omega_3(z)},
  \notag
\end{align}
where $\mathscr{P}^+(z)$ is a degree-$j p(3 p - 1)$ differential
polynomial in the fields (and the exponential has the Sugawara
dimension $j p^2 + j p + 2 p$).  All black dots in
Fig.~\ref{fig:a:left:j>0} are in $\ker\Fi\bbigcap\ker\Fii$, but we
select $\jW^+(z)$ as an extended algebra generator.

\begin{figure}[tb]
  \centering
  \begin{equation*}
    \xymatrix@C15pt@R25pt{
      \nothing="left"&&&&&&&\verma="R"&&&&\nothing="up"
      \\
      \nothing="left2"&&&&&&\Cverma="R2"&&&&
      \\
      &
      \\
      \Circ="cosing"&\ivermaiv="L3"&&&&&&&&\vermaiv="R3"&&\nothing="down"
      \ar@*{[|(2)]}@{-}"left";"R"
      \ar@*{[|(2)]}@{-}"left2";"R2"
      \ar@*{[|(2)]}@{-}"R";"R"+<50pt,-50pt>
      \ar@{{-}{--}{>}}|{\MFFplus{1, 2 p + 1}}"R"+<-1pt,0pt>;"R2"+<1pt,1pt>
      \ar|{\MFFminus{3 j p + 2,p}}"R"+<-1pt,0pt>;"R3"+<0pt,1pt>
      \ar@{{-}{--}{>}}|{\MFFplus{j p, 3 p + 1}}"R2"+<-3pt,-3pt>;"cosing"+<1pt,1pt>
      \ar@*{[|(4)]}@{-}"R3";"R3"+<20pt,-20pt>
      \ar@*{[|(4)]}@{-}"L3";"L3"-<20pt,20pt>
      \ar@*{[|(4)]}@{-}"R3";"R3"+<20pt,-20pt>
      \ar@*{[|(4)]}@{-}"L3";"R3"
      \ar"cosing";"L3"+<-2pt,0pt>
      \ar@{{<}{.}{>}}|{3 j p^2 + 2 p}"up";"down"
      \ar@{}|{\scriptstyle -j p - 1}"R";"R"+<6pt,16pt>
      \ar@{}|{\scriptstyle -j p - 2}"R2";"R2"+<-6pt,-26pt>
      \ar@{}|{\scriptstyle 2 j p + 1\quad\ }"R3";"R3"+<6pt,-26pt>
      \ar@{}|{\scriptstyle\ \  -2 j p - 1}"L3";"L3"+<6pt,-26pt>
    }
  \end{equation*}
  \caption{\small Relevant structures in the Wakimoto-type $\hSL2$
    module associated with the field $\AsymmF_{-j p - 1}(z)$, which
    sits in the top right corner.  Dashed lines show Verma-module
    singular vectors that vanish in the free-field realization.  Open
    dots are ``cosingular'' vectors at the grades where a singular
    vector vanishes.  The bottom right corner is a nonvanishing
    singular vector.  Charges of the states are
    indicated.}\label{fig:a:right:j>0}
\end{figure}
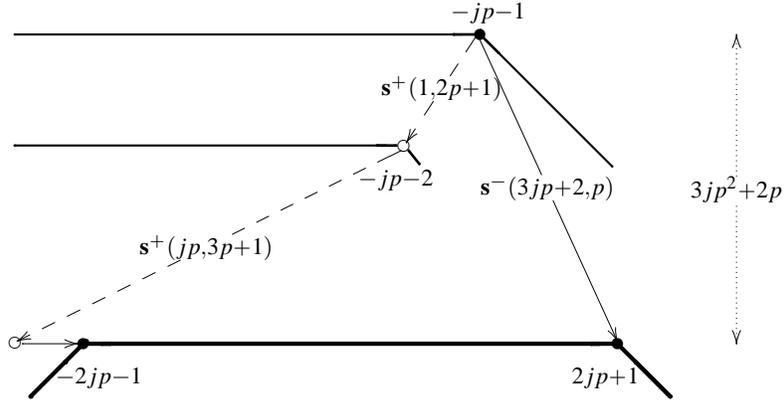

Figure~\ref{fig:a:left:j>0} is a ``refinement'' of the left part of
Fig.~\ref{fig:a-TRIPLET} for $h$ in~\eqref{eq:h(j)}, $j\geq 1$.
Instead of the right part of Fig.~\ref{fig:a-TRIPLET}, we then have
Fig.~\ref{fig:a:right:j>0}, where the top right corner represents
$\AsymmF_{-j p - 1}(z)$ and the bottom right corner is the
nonvanishing singular vector
\begin{align}\label{jWminus}
  \jW^-(z)&=\MFFminus{3 j p + 2,p}\AsymmF_{-j p - 1}(z),
  \\
  \intertext{which lies in $\ker\Fi\bbigcap\ker\Fii$ and has the
    structure}
  \jW^-(z)&=\mathscr{P}^-(z)\,
  e^{(3 j p + 2) \omega_1(z) - 2 (j p + 1) \omega_2(z) + 2 (2 j p + 1) \omega_3(z)},
  \notag
\end{align}
where $\mathscr{P}^-(z)$ is a differential polynomial in the fields of
the degree $(p - 1)(3 j p + 2)$ (and the exponential has the Sugawara
dimension $j p^2 + 3 j p + 2$).

The zero-mode $s\ell(2)$ algebra generates a $(4 j p + 3)$-plet from
$\jW^-(z)$ as well as from $\jW^+(z)$.  Because the Sugawara
dimensions of the top right corners in Figs.~\ref{fig:a:left:j>0}
and~\ref{fig:a:right:j>0} are
\begin{equation*}
  \dim\AsymmF_{j p + 1}(z)= j p^2 + 2 p
  \quad\text{and}\quad
  \dim\AsymmF_{-j p - 1}(z)= j p^2,
\end{equation*}
the $(4 j p + 3)$-plets in Figs.~\ref{fig:a:left:j>0} and
\ref{fig:a:right:j>0} have the same Sugawara dimension $4 j p^2 + 2
p$.

Without drawing another picture, we briefly describe the relevant
structure of the $\hSL2$-module associated with $\AsymmF_{0}(z)=\one$
(the unit operator).  There, the singular vector
\begin{equation}\label{jWnaught}
  \jW^0(z)=\MFFminus{2 j p + 1, 2 p}\one(z)
\end{equation}
is nonvanishing in the free-field realization and is also the
rightmost element of a zero-mode $(4 j p + 3)$-plet located at the
same Sugawara dimension as the two $(4 j p + 3)$-plets containing
$\jW^+(z)$ and $\jW^-(z)$.

We summarize our findings as the following conjecture on the extended
algebra.  For fixed integers $p\geq 2$ and $j\geq 1$, let
\begin{equation*}
  r_a = (2-a)j p + 1 - a,
  \quad\text{and}\quad
  s_a = (2+a)p
  \quad\text{for}\quad
  a = 1, 0, -1.
\end{equation*}%
\begin{it}%
  The three dimension-$(4 j p^2 + 2 p)$ fields \eqref{jWplus},
  \eqref{jWnaught}, \eqref{jWminus}, which can also be written~as
  \begin{equation}\label{a:the-W-j}
    \jW^a(z)=\MFFminus{r_a, s_a}\AsymmF_{a(j p + 1)}(z),\qquad
    a = 1, 0, -1,
  \end{equation}
  together with the entire $(4 j p + 3)$-plets $(\Jminus_0)^i
  \jW^a(z)$, $0\leq i\leq 4 j p + 2$, generate a $W$-algebra of
  mutually local fields in the kernel of the two screenings.
\end{it}

We call this $W$-algebra the triplet--multiplet algebra, $\WWtrimult$,
although its triplet structure is somewhat more elusive than that of
the $\WWtritri$ algebra in~\bref{a:j=0}: the long screening does not
map between elements of the triplet.

\subsection{Hamiltonian reduction}\label{a:hamred} The choice of
scalar fields made in~\bref{ortho} implies that the Hamiltonian
reduction of $J^{\pm,0}(z)$ and $\W^{\pm,0}(z)$ is obtained by simply
setting $\Afour(z)=\Athree(z)=0$, leaving us with only the matter
field $\Aii(z)$.

\subsubsection{$\myboldsymbol{j=0}$: the $\myboldsymbol{\pone}$
  triplet algebra}\label{HR:j=0}Setting $\Afour(z)=\Athree(z)=0$
reduces the three fields in \eqref{a:the-W-j=0} to three fields
generating \textit{the} triplet $\Wpone$ algebra, exactly as it was
obtained in~\cite{[FHST]}.  (In particular, formula~\eqref{delta(h)}
with $h=1$ gives the dimension $2p-1$ of the triplet algebra
generators).

Moreover, much as the $\hSL2$ currents were expressed in terms of the
``matter'' energy--momentum tensor $\Tm(z)$ and the two additional
scalars $\Athree(z)$ and $\Afour(z)$ in~\eqref{through-matter}, simple
analysis of the construction in~\eqref{a:the-W-j=0} readily shows how
to ``pack'' the generators \textit{and} certain fields of the $\Wpone$
algebra, properly dressed with the additional scalars, into the
$\W^{\pm,0}(z)$ generators of $\WWtritri$ (we do not need this
construction in this paper, however).

\subsubsection{$\myboldsymbol{j>0}$: the $\myboldsymbol{(p,p')}$
  algebras}
Setting $\Afour(z)=\Athree(z)=0$ in the expressions for
$\jW^{\pm,0}(z)$ in \eqref{a:the-W-j} reduces them to the
corresponding three fields generating the triplet $W_{p_+,p_-}$
algebra obtained in~\cite[\textbf{4.2.1.}]{[FGST3]}:
\begin{align*}
  W^+(z)&=\mathscr{P}_{j p(3 p - 1)}^+(z)\,e^{p\Aii(z)},
  \\
  W^0(z)&=\mathscr{P}_{(2 p - 1)(2 j p + 1)}^0(z),
  \\
  W^-(z)&=\mathscr{P}_{(p - 1)(3 j p + 2)}^-(z)\,e^{-p\Aii(z)},
\end{align*}
where $\mathscr{P}_m^{\pm,0}(z)$ are differential polynomials of the
indicated degree $m$ in $\partial^n\Aii(z)$, $n\geq1$.
In terms of $p_+ \eqdef p$ and $p_- \eqdef j p + 1$, the degrees of
these polynomials are $ j p(3p-1) = (p_- - 1)(3 p_+ - 1)$, $(2 p -
1)(2 j p + 1) = (2 p_+ -1)(2 p_- - 1)$, and $(p - 1)(3 j p + 2) = (p_+
- 1)(3 p_- - 1)$, which coincides with what we had in~\cite{[FGST3]}.
As regards the exponentials in the above formulas, they are also the
same as in~\cite{[FGST3]}, where the scalar field was normalized
canonically, giving rise to the factors $\pm\sqrt{p_+ p_-\mathstrut}
=\pm\sqrt{p\cdot p\cdot (\frac{2}{p} + 2 j)}$ in the exponents.

We also recall from~\cite{[FGST3]} that the OPE of $W^+(z)$ and
$W^-(z)$ is
\begin{equation}
  W^+(z)W^-(w)=\frac{S_{p_+,p_-}(T)}{(z-w)^{7p_+p_--3p_+-3p_-+1}}+
  \text{less singular terms},
\end{equation}
where $S_{p_+,p_-}(T)$ is the vacuum singular vector\,---\,the
polynomial of degree $\half
(p_+\,{-}\,1)\cdot$\linebreak[0]$(p_-\,{-}\,1)$ in $T$ and
$\partial^nT$, $n\geq1$, such that $S_{p_+,p_-}(T)=0$ is the
polynomial relation for the energy--momentum tensor in the $(p_+,p_-)$
Virasoro minimal model.  This degree is $\half j p(p-1)$ in our
current case.


\subsection{$\WW\myboldsymbol{(2,(2p)^{\times 3\times 3})}$
  highest-weight states}
\label{a:W-modules}
We return to the triplet--triplet algebra (in particular, we now have
$k = \frac{1}{p} - 2$) and construct a class of its highest-weight
states.

\subsubsection{}
For integer $r$, $s$, and $\vartheta$, we set
\begin{equation*}
  \Vsr_{s,r;\vartheta}(z)=
  e^{
    \frac{r}{p}\omega_1(z) - \frac{s - 1}{p} \omega_2(z)
    - \frac{s - 1 - 2 r + \vartheta - 2 p \vartheta}{p} \omega_3(z)
  }.
\end{equation*}
This is a twisted relaxed $\hSL2$ highest-weight state of twist
$\vartheta-1$: the top modes $J^{\pm}_n$ that do not
annihilate it act as
\begin{align*}
  \Jplus_{\vartheta-1}\Vsr_{s,r,\vartheta}(z) &= \Vsr_{s,r+p,\vartheta}(z),
  \\
  \Jminus_{1-\vartheta}\Vsr_{s,r,\vartheta}(z) &= -\ffrac{r (p + r - s)}{p^2}
  \Vsr_{s,r-p,\vartheta}(z).
\end{align*}
The Sugawara dimension of $\Vsr_{s,r,\vartheta}(z)$ is
\begin{equation}
  \Delta_{s,r;\vartheta}=
  \ffrac{(s+\vartheta -1)^2-4 r (\vartheta -1)}{4 p}
  +\half (1 - s - \vartheta^2).
\end{equation}

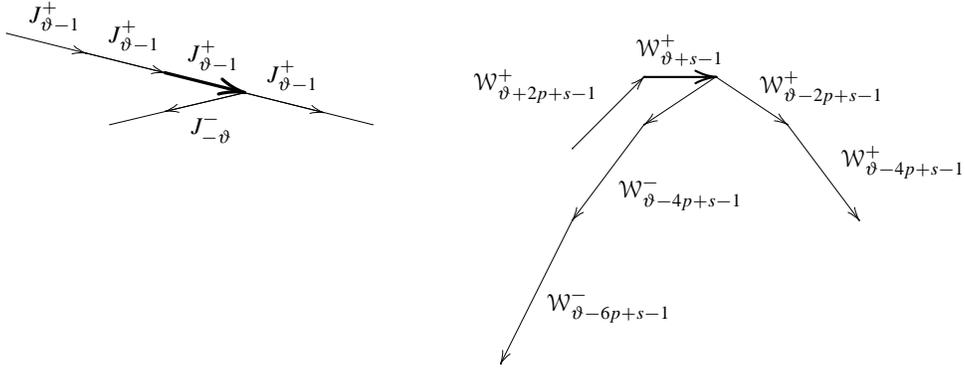
\begin{figure}[tb]
  \centering
\begin{equation*}
    \xymatrix@R=1pt@C=1pt{
      \nothing="left"&&&&&&&&&&&&&
      \\
      &&&&&&&&&&&&&
      \\
      &&&&&&&&&&&&&
      \\
      &\nothing="down"&&&&&&&&&&&\nothing="right"&
      \\
      &&&&&&&&&&&&&&
      \ar@{-}"left"+<0pt,0pt>;"right"+<0pt,0pt>
      \ar@*{[|(3)]}^{\Jplus_{\vartheta-1}}"left"+<30pt,-7.5pt>;"left"+<60pt,-15pt>="middle"
      \ar@{-}"middle"+<0pt,0pt>;"down"+<0pt,0pt>
      \ar^{\Jminus_{-\vartheta}}"middle"+<0pt,0pt>;"middle"+<-30pt,-7pt>
      \ar^{\Jplus_{\vartheta-1}}"middle"+<0pt,0pt>;"middle"+<30pt,-7.5pt>
      \ar^{\Jplus_{\vartheta-1}}"left"+<-30pt,7.5pt>;"left"+<0pt,0pt>="two"
      \ar^{\Jplus_{\vartheta-1}}"two"+<0pt,0pt>;"two"+<30pt,-7.5pt>
    }
    \quad
    \xymatrix@R=1pt@C=1pt{
      &&&&&&&&&&&&&&&&
      \\
      &&&&&&&\nothing="upup"&&&\nothing="up"&&&&&&&
      \\
      &&&&&&&&&&&&&&&&
      \\
      &&&&&&&\nothing="1"&&&&&&\nothing="right"&&&
      \\
      &&&&\nothing="upleft"&&&&&&&&&&&&
      \\
      &&&&&&&&&&&&&&&&
      \\
      &&&&&&&&&&&&&&&&
      \\
      &&&&\nothing="2"&&&&&&&&&&&&\nothing="right2"
      \\
      &&&&&&&&&&&&&&&&
      \\
      &&&&&&&&&&&&&&&&
      \\
      &&&&&&&&&&&&&&&&
      \\
      &&&&&&&&&&&&&&&&
      \\
      &&&&&&&&&&&&&&&&
      \\
      &\nothing="3"&&&&&&&&&&&&&&
      \ar^(.5){\W^{+}_{\vartheta+2p+s-1}}"upleft"+<0pt,0pt>;"upup"+<0pt,0pt>
      \ar@*{[|(2)]}^{\W^{+}_{\vartheta+s-1}}"upup"+<0pt,0pt>;"up"+<0pt,0pt>
      \ar^(.8){\W^{+}_{\vartheta-2p+s-1}}"up"+<0pt,0pt>;"right"+<0pt,0pt>
      \ar^(.6){\W^{+}_{\vartheta-4p+s-1}}"right"+<0pt,0pt>;"right2"+<0pt,0pt>
      \ar"up"+<0pt,0pt>;"1"+<0pt,0pt>
      \ar^{\W^{-}_{\vartheta-4p+s-1}}"1"+<0pt,0pt>;"2"+<0pt,0pt>
      \ar^{\W^{-}_{\vartheta-6p+s-1}}"2"+<0pt,0pt>;"3"+<0pt,0pt>
    }
  \end{equation*}
    \caption{\small\textsc{Left:} The modes $\Jplus_{\vartheta-1}$ map
    between the $\Vsr_{s,r,\vartheta}$ states; each
    $\Jplus_{\vartheta-1}$ arrow except the boldfaced one is
    invertible by the action of $\Jminus_{-\vartheta+1}$.  The state
    at the tip of the boldfaced arrow is annihilated by
    $\Jminus_{-\vartheta+1}$, and is mapped into a nonzero state by
    $\Jminus_{-\vartheta}$.  \ \ \textsc{Right:} The action of
    $\W^{\pm}$ modes on a $\Vsr_{s,r;\vartheta}^{0,\mu}[0,m](z)$
    state, which sits at the top vertex.  The arrow without a $\W$
    mode indicated represents $\W^{-}_{\vartheta-2p+s-1}$. The values
    of $n$ in $\Vsr_{s,r;\vartheta}^{0,\mu}[n,m](z)$ change from
    vertex to vertex in accordance with~\eqref{Wplus-act}
    and~\eqref{Wminus-act}.  For the submodule diagram, the value of
    $s$ governs its opening degree at the top.  The diagram lies in a
    plane intersecting the plane of the left diagram (and sharing the
    vertical direction with~it).}
  \label{fig:extremal}
\end{figure}

An extra vanishing condition,
$\Jminus_{1-\vartheta}\Vsr_{s,r;\vartheta}(z)=0$, occurs whenever
\begin{equation*}
  r = 0\quad\text{or}\quad r = s - p.
\end{equation*}
In these two cases, the corresponding operators
$\Vsr_{s,r;\vartheta}(z)$ are twisted highest-weight states:
\begin{align*}
  \Vsr_{s,0;\vartheta}(z)&= e^{\frac{1 - s}{p} \omega_2(z) + \frac{1 -
      s - \vartheta + 2 p \vartheta}{p} \omega_3(z)}
  \doteq\ket{\lambda^+(1,s);\vartheta},
  \\
  \Vsr_{s,s-p;\vartheta}(z)&= e^{\frac{s - p}{p} \omega_1(z) + \frac{1
      - s}{p} \omega_2(z) + \frac{1 - 2 p + s - \vartheta + 2 p
      \vartheta}{p} \omega_3(z) }
  \doteq\ket{\lambda^-(1,s+1);\vartheta}
\end{align*}
(with $\lambda^{\pm}(r,s)$ defined in~\bref{MFFthm}).  The occurrence
of a twisted highest-weight state is illustrated in the left part of
Fig.~\ref{fig:extremal}.  We also note that the \textit{matter}
dimension of each of these two states is $\Deltamatter_{s, 1}$
(see~\eqref{matter-fields-dim}).

\subsubsection{}\label{V[n,m]}
We reparameterize $\Vsr_{s,r;\vartheta}(z)$ by defining
\begin{align*}
  \Vsr_{s,r;\vartheta}^{\nu,\mu}[n,m](z) &=
  \Vsr_{s - p \nu - 2 p n, r + p \mu + 2 p m; \vartheta}(z)
  \\
  &=e^{(\frac{r}{p} + \mu + 2 m)\omega_1(z)
    + (\frac{1 - s}{p} + 2 n + \nu) \omega_2(z)
    + (\frac{1 + 2 r - s - \vartheta}{p} + 2 (n + 2 m + \vartheta + \mu) + \nu)
    \omega_3(z)
  },
\end{align*}
where now 
\begin{equation*}
  1\leq s\leq p,\quad 0\leq r\leq p-1,\quad \nu,\mu=0,1,\quad
  n,m,\vartheta\in\oZ.
\end{equation*}
This allows representing the $\WWtritri$ action very conveniently.  It
follows that $\Vsr_{s,r;\vartheta}^{\nu,\mu}[n,m](z)$ is always
annihilated by $\W^+_{i}$ with $i\geq \vartheta + s - p\nu - 2 p (n +
1)$ and by $\W^-_{i}$ with $i\geq \vartheta - s + p\nu + 2 p n$.  The
top modes of $\W^{+}$ and $\W^{-}$ that do not annihilate
$\Vsr_{s,r;\vartheta}^{\nu,\mu}[n,m](z)$ identically act~as
\begin{align}
  \label{Wplus-act}
  \W^+_{\vartheta + s - p\nu -2 p(n+1) - 1} \Vsr_{s,r;\vartheta}^{\nu,\mu}[n,m](z) &= \Vsr_{s,r;\vartheta}^{\nu,\mu}[n+1,m](z),
  \\
  \label{Wminus-act}
  \W^-_{\vartheta - s+ p\nu +2 p n - 1}
  \Vsr_{s,r;\vartheta}^{\nu,\mu}[n,m](z)&=\\
  &\kern-120pt\ffrac{(2 p - 1)!}{p^{4 p - 3}}
  \prod_{i=1}^{p - 1}(s - p\nu -2 pn-i)(s - p\nu -2 pn+i)\;
   \Vsr_{s,r;\vartheta}^{\nu,\mu}[n-1,m](z).
   \notag
\end{align}
The $n$ label changes as $n\mapsto n \pm 1$ in these formulas, and it
follows that a $\WWtritri$ module generically contains the direct sum
of subspaces spanned by the $\Vsr_{s,r;\vartheta}^{\nu,\mu}[n,m](z)$
over all~$n$.  We now detect cases where such a module is reducible
and then identify a submodule in~it.

\subsubsection{}For $n=0$ and $\nu=0$, annihilation conditions hold in
the form
\begin{align*}  
  \W^-_{i}
  \Vsr_{s,r;\vartheta}^{0,\mu}[0,m](z)
  &= 0,\quad i\geq \vartheta + s - 2p,
  \\
  \intertext{and for $n=0$ and $\nu=1$,}
  \W^-_{i}
  \Vsr_{s,r;\vartheta}^{1,\mu}[0,m](z) &= 0,\quad i\geq\vartheta.
\end{align*}
Hence, for each $s$, $1\leq s\leq p-1$, a submodule is generated from
$\Vsr_{s,r;\vartheta}^{\nu,\mu}[0,m](z)$.  For $\nu=0$, this submodule
is illustrated in the right part of Fig.~\ref{fig:extremal}.

We return to the consideration of these vertex operators in
Sec.~\ref{sec:7}, where we also discuss the relation to the $\algU$
algebra.

\section{\textbf{From $\Nich(X)$ to extended chiral algebras: symmetric
    case}}\label{CFTsymm}
We take two screening operators $F_1$ and $F_2$ with the momenta
defined by the second line in~\eqref{three-R}.  The corresponding
braiding matrix $Q_{\symm} = Q_{\symm}(j,p)$ is then the one
in~\eqref{qij-first}.  We seek the kernel $\ker
F_1\bbigcap{}$\linebreak[0]$\ker F_2$.  We fix our notation
in~\bref{sVir}, identify the parafermions in the kernel and convert
them into $\hSL2$ currents in~\bref{s-para}, and then use the $\hSL2$
representation theory to find other pieces of the kernel: the ``easy''
one in~\bref{s-more} and the ``difficult'' in~\bref{s-even-more},
where we identify the extended algebra generators.

\subsection{Scalar fields and a Virasoro algebra}\label{sVir}
We introduce two scalar fields $\Aone(z)$ and $\Atwo(z)$ with the OPEs
determined by scalar products in the second line in~\eqref{three-R}:
\begin{alignat*}{2}
\Aone(z)\, \Aone(w) &= \log(z - w),\qquad&
\Aone(z)\, \Atwo(w) &= \bigl(\ffrac{1}{p} + j\bigr)\log(z - w),\\
\Atwo(z)\, \Atwo(w) &= \log(z - w).
\end{alignat*}
It follows from these OPEs that the kernel of the two screenings
\begin{equation*}
  F_1=\oint e^{\Aone}\quad\text{and}\quad F_2=\oint e^{\Atwo}
\end{equation*}
contains the energy--momentum tensor
\begin{multline}
  T(z)= -\ffrac{1}{2 k (k + 2)} \partial\Aone\partial\Aone(z)
  + \ffrac{k + 1}{k (k + 2)} \partial\Aone\partial\Atwo(z)
  - \ffrac{1}{2 k (k + 2)}\partial\Atwo\partial\Atwo(z)
  \\
  {}- \ffrac{1}{2 (k + 2)}\partial^2\Aone(z)
  - \ffrac{1}{2 (k + 2)}\partial^2\Atwo(z),
\end{multline}
where we set
\begin{equation}\label{s-k2p}
  k \eqdef \ffrac{1}{p} + j - 1.
\end{equation}
In terms of this $k$, the central charge of $T(z)$ is expressed as
in~\eqref{eq:cc-asymm}.

It is convenient in what follows to introduce ``fundamental weights''
$\omega_1,\omega_2\in\oC^2$ for the vectors $\alpha_1$ and $\alpha_2$
in the second line in~\eqref{three-R}:
\begin{equation*}
  \omega_1 = \ffrac{1}{k (k + 2)} (-\alpha_1 + (k + 1) \alpha_2),\qquad
  \omega_2 = \ffrac{1}{k (k + 2)} ((k + 1) \alpha_1 - \alpha_2).
\end{equation*}
Then
\begin{equation*}
  \omega_1.\omega_1 = -\ffrac{1}{k (k + 2)},
  \qquad
  \omega_1.\omega_2 = \ffrac{k + 1}{k (2 + k)},
  \qquad
  \omega_2.\omega_2 = -\ffrac{1}{k (k + 2)}.
\end{equation*}
Anticipating the appearance of a third scalar field, we pass from
$\oC^2$ to $\oC^3$ by adding the vector $\alpha_3$, such that
$\alpha_3.\alpha_1=\alpha_3.\alpha_2=0$ and $\alpha_3.\alpha_3 = 2k$,
and the corresponding ``fundamental weight'' $\omega_3=\frac{1}{2 k}
\alpha_3$.  With a slight abuse of notation, we write $\omega_i(z)$
for the corresponding linear combinations of our scalar fields:
\begin{equation*}
  \omega_1(z) = \ffrac{1}{k (k + 2)} (-\Aone(z) + (k + 1) \Atwo(z)),\qquad
  \omega_2(z) = \ffrac{1}{k (k + 2)} ((k + 1) \Aone(z) - \Atwo(z)),
\end{equation*}
and $\omega_3(z)=\frac{1}{2 k} \Athree(z)$.

\subsubsection{A long screening} The above fields $T(z)$ and
$j^{\pm}(z)$ are also in the kernel of the ``long'' screening
\begin{equation*}
  \Ebeta=\oint\partial\Aone\,e^{-\frac{1}{k + 2}(\Aone + \Atwo)}.
\end{equation*}
(Up to a coefficient, $\Ebeta$ is equal to $\oint(a_1\partial\Aone +
a_2\partial\Atwo)\, e^{-\frac{1}{k + 2}(\Aone + \Atwo)}$ for any
$a_1\neq a_2$, because $\oint(\partial\Aone +\partial\Atwo)\,
e^{-\frac{1}{k + 2}(\Aone + \Atwo)}=0$.) \ It follows that
\begin{equation*}
  [\Ebeta,F_1]=0,\qquad
  [\Ebeta,F_2]=0.
\end{equation*}

\subsection{Parafermionic fields in
  $\bker\myboldsymbol{F_1\bigcap{}}\bker\myboldsymbol{F_2}$ and
  $\myboldsymbol{\hSL2_k}$}\label{s-para}
The kernel of the two screenings contains parafermionic fields
\begin{equation}\label{para:symm}
  \begin{aligned}
    \jplus(z)&=\partial\Aone(z)
    e^{-\omega_1(z) + \omega_2(z)},
    \\
    \jminus(z)&=\partial\Atwo(z)
    e^{\omega_1(z) - \omega_2(z)}
  \end{aligned} 
\end{equation} 
(we note that $\jplus(z)$ is $F_1$-exact and $\jminus(z)$ is
$F_2$-exact).  These nonlocal fields can be dressed into $\hSL2$
currents by introducing an auxiliary scalar with the OPE
in~\eqref{localize-2} (of course, with $k$ given
by~\eqref{s-k2p}):
\begin{equation*}
  \begin{aligned}
    \Jplus(z)&=
    \partial\Aone(z)e^{-\omega_1(z) + \omega_2(z) + 2 \omega_3(z)},
    \\
    \Jnaught(z)&=\half\Athree(z),
    \\
    \Jminus(z)&=
    \partial\Atwo(z)
    e^{\omega_1(z) - \omega_2(z) - 2 \omega_3(z)}.
  \end{aligned}
\end{equation*}

\subsection{More of $\bker\myboldsymbol{F_1\bigcap{}}\bker
  \myboldsymbol{F_2}$}\label{s-more} Another part of $\ker
F_1\bbigcap\ker F_2$ is easy to find.  For any $h=1,2,\dots$, the
field
\begin{equation*}
  \SymmF_{h}(z)
  = e^{h \omega_1(z) + h \omega_2(z)}
  = e^{\frac{h}{k+2}\Aone(z) + \frac{h}{k+2}\Atwo(z)}
\end{equation*}
is in the kernel of both screenings, is a relaxed highest-weight
state,
\begin{equation*}
  \Jminus_{1}\SymmF_{h}(z)=
  \Jnaught_{1}\SymmF_{h}(z)= \Jplus_{1}\SymmF_{h}(z)=0,
\end{equation*}
and is the central element in a $(2 h + 1)$-plet with respect to the
zero-mode $s\ell(2)$ algebra:
\begin{equation*}
  \begin{alignedat}{2}
    (\Jplus_0)^h\SymmF_{h}(z)&\neq 0,&\qquad
    (\Jplus_0)^{h+1}\SymmF_{h}(z)&= 0,\\ 
    (\Jminus_0)^h\SymmF_{h}(z)&\neq 0,&\qquad
    (\Jminus_0)^{h+1}\SymmF_{h}(z)&= 0,
  \end{alignedat}\qquad h\geq 1.
\end{equation*}
In particular, $(\Jplus_0)^h\SymmF_{h}(z)$ is a highest-weight state
and $(\Jminus_0)^h\SymmF_{h}(z)$ is a $\vartheta=1$ twisted
highest-weight state, which we represent as
\begin{equation}\label{s:triplet-easy}
  \xymatrix@C15pt@R24pt{
    \ivermaiv="fromR"&&*{\bullet}\ar@{}|{\SymmF_h}[];[]+<5pt,24pt>&&\vermaiv="fromL"
    \ar@*{[|(4)]}@{-}"fromL";"fromR"
    \ar@{}|{\scriptstyle (\Jminus_0)^h\SymmF_{h}}"fromR";"fromR"+<5pt,24pt>
    \ar@{}|{\scriptstyle(\Jplus_0)^h\SymmF_{h}}"fromL";"fromL"+<-5pt,24pt>
    \ar@*{[|(4)]}@{-}"fromR";"fromR"+<-20pt,-20pt>
    \ar@*{[|(4)]}@{-}"fromL";"fromL"+<20pt,-20pt>
  }
\end{equation}

\noindent
The Sugawara dimension of the fields in this $(2 h + 1)$-plet is
$\frac{h (h + 1)}{k + 2}=\frac{h (h + 1) p}{(j+1) p + 1}$.

\subsection{Even more of $\bker\myboldsymbol{F_1\bigcap{}}\bker
  \myboldsymbol{F_2}$ and the extended algebras}\label{s-even-more}
The field $\SymmF_{h}(z)$ with a negative integer $h$ is not in the
kernel of either $F_1$ or $F_2$, but is simply related to fields that
are in \textit{one} of these kernels, as we now describe.

In what follows, we have to consider $\SymmF_{h}(z)$ and
$\SymmF_{-h}(z)$ simultaneously; we therefore assume $h$ to be a
positive integer parameter from now on, and define
\begin{align*}
  \Lstate_{-h}(z)&=\frac{1}{\prod\limits_{i=0}^{h - 1}(-h - i)}
  (\Jminus_0)^{h} \SymmF_{-h}(z) =
  e^{-2 h \omega_2(z) - 2 h \omega_3(z)}
  \in\ker F_1,
  \\
  \Rstate_{-h}(z)&=\frac{1}{\prod\limits_{i=0}^{h - 1}(-h - i)}
  (\Jplus_0)^{h} \SymmF_{-h}(z) =
  e^{-2 h \omega_1(z) + 2 h \omega_3(z)}
  \in\ker F_2
\end{align*}
(with $\Lstate_{-h}(z)\notin\ker F_2$ and $\Rstate_{-h}(z)\notin\ker
F_1$).  These two states are identified with a highest-weight state
and a $\vartheta=1$ twisted highest-weight state~as
\begin{equation}\label{eq:LR}
  \Lstate_{-h}(z)\doteq\ket{-h},\qquad
  \Rstate_{-h}(z)\doteq\ket{h+\fffrac{k}{2};1},
\end{equation}
which we show with corners in the diagram
\begin{equation}\label{s:module-top}
  \xymatrix@C15pt@R24pt{
    &&&&\verma="L"&\Circ="V"&&&&\Circ="U"&\iverma="R"&&
    \ar@*{[|(2)]}@{-}"L";"L"-<80pt,0pt>
    \ar@*{[|(2)]}@{-}"L";"L"+<30pt,-30pt>
    \ar@*{[|(2)]}@{-}"R";"R"+<80pt,0pt>
    \ar@*{[|(2)]}@{-}"R";"R"+<-30pt,-30pt>
    \ar@*{[|(2)]}@{-}"V"+<0pt,1pt>;"U"+<0pt,1pt>
    \ar@*{[|(1)]}"V"+<0pt,1pt>;"L"
    \ar@*{[|(1)]}"U"+<0pt,1pt>;"R"
    \ar@{}|{\scriptstyle\Ktopleft}"V";"V"+<6pt,20pt>
    \ar@{}|{\scriptstyle\Ktopright}"U";"U"+<6pt,20pt>
    \ar@{}|{\ \ \scriptstyle\Rstate_{-h}}"R";"R"+<6pt,16pt>
    \ar@{}|{\scriptstyle\Lstate_{-h}\ \ }"L";"L"+<6pt,16pt>
    \ar@{}|{\scriptstyle -h\ \ }"L";"L"+<5pt,-24pt>
    \ar@{}|{\scriptstyle\ \  h}"R";"R"+<-5pt,-24pt>
  }
\end{equation}

\noindent
where $\protect\Ktopleft=(\Jminus_0)^{h-1} \SymmF_{-h}(z)$ and
$\protect\Ktopright=(\Jplus_0)^{h-1} \SymmF_{-h}(z)$, and $-h$ and $h$
indicate the charges (eigenvalues of $\Jnaught_0$) of the appropriate
states.

The intersection of the kernels is to be sought deeper in the module
whose top is shown above.  The actual picture depends on $j$, and we
here restrict ourself to $j\geq-1$; the case $j=-1$ is special.

\subsubsection{The case $\myboldsymbol{j=-1}$ and a triplet--triplet
  algebra}\label{s:j=-1}
\begin{figure}[tb]
  \centering
  \begin{equation*}
    \xymatrix@C18pt@R32pt{
      &&&&\verma="L"&\Circ="V"&&\Circ="U"&\iverma="R"&&\\
      &&&&&&&&&&&&\\
      &&&&&&&&&&&&\\
      &&&\Circ="iL"&\ivermaiv="fromR"&&&&\vermaiv="fromL"&\Circ="iR"&
      \ar@*{[|(2)]}@{-}"L";"L"-<80pt,0pt>
      \ar@*{[|(2)]}@{-}"L";"L"+<30pt,-30pt>
      \ar@*{[|(2)]}@{-}"R";"R"+<80pt,0pt>
      \ar@*{[|(2)]}@{-}"R";"R"+<-30pt,-30pt>
      \ar@*{[|(2)]}@{-}"V"+<0pt,1pt>;"U"+<0pt,1pt>
      \ar@*{[|(1)]}"V"+<0pt,1pt>;"L"
      \ar@*{[|(1)]}"U"+<0pt,1pt>;"R"
      \ar|(.8){\MFFminus{2 h, p}}"L"+<-1pt,0pt>;"fromL"+<-1pt,2pt>
      \ar@{{-}{--}{>}}|(.8){\MFFplus{1,2 h p + 1}}"L";"iL"
      \ar|(.7){\MFFplus{2 h, p;1}}"R"+<1pt,0pt>;"fromR"+<1pt,2pt>
      \ar@{{-}{--}{>}}|(.7){\MFFminus{1,2 h p + 1;1}}"R";"iR"
      \ar"iL";"fromR"+<-2pt,0pt>
      \ar"iR";"fromL"+<2pt,0pt>
      \ar@*{[|(4)]}@{-}"fromL";"fromR"
      \ar@*{[|(4)]}@{-}"fromL";"fromL"+<15pt,-15pt>
      \ar@*{[|(4)]}@{-}"fromR";"fromR"+<-15pt,-15pt>
      \ar@{}|{\scriptstyle\Ktopleft}"V";"V"+<6pt,20pt>
      \ar@{}|{\scriptstyle\Ktopright}"U";"U"+<6pt,20pt>
      \ar@{}|{\ \ \scriptstyle\Rstate_{-h}}"R";"R"+<6pt,16pt>
      \ar@{}|{\scriptstyle\Lstate_{-h}\ \ }"L";"L"+<6pt,16pt>
      \ar@{}|{\scriptstyle -h\ \ }"fromR";"fromR"+<5pt,-24pt>
      \ar@{}|{\scriptstyle\ \  h}"fromL";"fromL"+<-5pt,-24pt>
      \ar@{}|{\scriptstyle -h-1\ \ \ }"iL";"iL"+<5pt,-20pt>
      \ar@{}|{\scriptstyle\ \  h+1}"iR";"iR"+<5pt,-20pt>
    }
  \end{equation*}
  \caption{\small Left open circle shows a ``cosingular vector'' in
    the grade where $\MFFplus{1,2 h p + 1}\Lstate_{-h}(z)$ vanishes,
    and similarly for the right open circle.  The $(2h+1)$-plet in the
    lower part is in $\ker F_1\bbigcap\ker F_2$.}\label{fig:s-TRIPLET}
\end{figure}
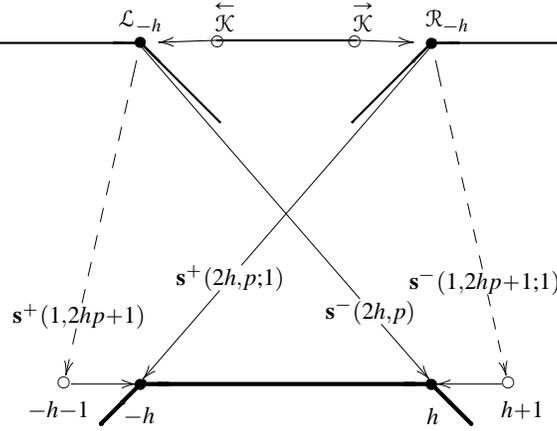

We recall that $h=1,2,\dots$.  For $j=-1$, with $k+2=\frac{1}{p}$, the
Verma-module highest-weight state $\ket{-h}$ has singular vectors
$\MFFplus{1,2 h p + 1}$ and $\MFFminus{2h,p}$, on the relative level
$2 h p$ both and at the respective charge grades $-h-1$ and $h$. \ In
our free-field realization, the first of these vanishes, which we show
with the left dashed arrow in Fig.~\ref{fig:s-TRIPLET},
but the second one is nonvanishing (the north-west--south-east arrow
in Fig.~\ref{fig:s-TRIPLET}).  Next, the Verma-module highest-weight
state $\ket{h+\frac{k}{2}}$ has singular vectors $\MFFplus{2 h,p}$ and
$\MFFminus{1,2 h p + 1}$, which lie at the respective charges $-h$ and
$h+1$.  After the twist by $\vartheta=1$, both $\MFFplus{2
  h,p;1}\ket{h+\frac{k}{2};1}$ and $\MFFminus{1,2 h p +
  1;1}\ket{h+\frac{k}{2};1}$ are on the relative level $2 h p$.  The
singular vector $\MFFminus{1,2 h p + 1;1}$ vanishes when evaluated on
$\Rstate_{-h}(z)$, which is shown by the right dashed arrow in
Fig.~\ref{fig:s-TRIPLET}; ``cosingular vectors'' in the corresponding
grades are shown with open circles.

The nonvanishing singular vectors $\MFFminus{2h,p}\Lstate_{-h}(z)$ and
$\MFFplus{2 h,p;1}\Rstate_{-h}(z)$ are endpoints of a $(2h+1)$-plet,
which is in $\ker F_1\bbigcap\ker F_2$.  The fields in this
$(2h+1)$-plet have the same Sugawara dimension, equal to $h(h+1)p$, as
the fields of the $(2h+1)$-plet
in~\eqref{s:triplet-easy}.\footnote{Figure~\ref{fig:s-TRIPLET} is a
  \textit{symmetric} counterpart of the right picture in
  Fig.~\ref{fig:a-TRIPLET}, but we must not forget that the screenings
  in the two cases are different.}

Continuing the embedding diagrams of Wakimoto modules, it is not
difficult to describe $\ker F_1\bbigcap$\linebreak[0]$\ker F_2$ quite
explicitly, but we here stop at the level in the module where the
fields presumably generating the extended algebra are located.  For
this, guided by mutual locality, we set $h=1$ in the above formulas;
the $(2h+1)$-plets are then triplets, which we now write in more
detail.

The ``plus'' triplet (Eq.~\eqref{s:triplet-easy}) is
\begin{equation}\label{s:wplus-x3}
  \begin{alignedat}{2}
    \W^{+}(z) &=\fffrac{1}{2}\Jplus_0\SymmF_1(z) &&=
    \fffrac{1}{2} e^{2 \omega_2(z) + 2 \omega_3(z)},
    \\
    \Jminus_0\W^{+}(z) &=\SymmF_1(z) &&=
    e^{\omega_1(z) + \omega_2(z)},
    \\
    (\Jminus_0)^2\W^{+}(z) &=\Jminus_0\SymmF_1(z)&&=
    e^{2 \omega_1(z) - 2 \omega_3(z)}.
  \end{alignedat}
\end{equation}
The action of the long screening generates the ``middle'' triplet
\begin{equation}\label{s:wnaught-x3}
  \begin{aligned}
    \W^{0}(z)&=\Ebeta\W^{+}(z),
    \\
    \Jminus_0\W^{0}(z)&=\Ebeta\Jminus_0 \W^{+}(z),
    \\
    (\Jminus_0)^2\W^{0}(z)&=\Ebeta(\Jminus_0)^2\W^{+}(z).
  \end{aligned}
\end{equation}
The grade of $\W^{0}(z)$ coincides with the grade of the level-$2p$
singular vector $\MFFminus{1,2p}$ in the $\hSL2$ Verma module
$\rep{M}_{0}$ with zero weight, and $(\Jminus_0)^2\W^{0}(z)$ is one
grade to the right of the level-$2p$ singular vector $\MFFplus{2,p+1}$
in $\rep{M}_{0}$; both these singular vectors vanish in our free-field
realization (and hence the ``middle'' triplet is not in the
Wakimoto-type module associated with~$\rep{M}_{0}$).

Acting on~\eqref{s:wnaught-x3} with the long screening generates the
``minus'' triplet (Fig.~\ref{fig:s-TRIPLET}):
\begin{alignat}{2}
  \notag
  \W^{-}(z) &=
  \Ebeta \W^{0}(z)
  &&
  = \MFFminus{2,p}
  e^{-2 \omega_2(z) - 2 \omega_3(z)},
  \\
  \label{s:wminus-x3}
  \Jminus_0\W^{-}(z)&=\Ebeta\Jminus_0\W^{0}(z),
  \\
  \notag
  (\Jminus_0)^2\W^{-}(z)
  &=\Ebeta(\Jminus_0)^2\W^{0}(z)
  &&
  = \MFFplus{2,p;1} e^{-2 \omega_1(z) + 2 \omega_3(z)}
\end{alignat}
(with the equalities to singular vectors holding up to nonzero
factors).  We propose these dimension-$2p$ fields, together with the
$\hSL2$ currents, as the \textit{triplet--triplet algebra generators}
in the ``symmetric'' realization.  Conjecturally, the algebra contains
all mutually local fields in $\ker F_1\bbigcap\ker F_2$.

An example of the triplet--triplet algebra in the ``symmetric''
realization is given in~\bref{app:example-s}.

\subsubsection{Cases $\myboldsymbol{j\geq 0}$ and triplet--multiplet
  algebras}\label{sec:s-j}
We briefly discuss the local algebra in the kernel for $j\geq 0$ in
the symmetric realization.  The relevant operator products contain
characteristic factors $(z-w)^{\frac{2 h h' p}{(j+1) p + 1}}$; for
locality, we therefore choose the smallest $h$ ensuring integers in
the exponent,
\begin{equation*}
  h  = (j+1) p + 1.
\end{equation*}
In~\eqref{s:triplet-easy} with this $h$, the state
$(\Jplus_0)^h\SymmF_{h}(z)$ (the top right corner) has a nonvanishing
singular vector
\begin{equation*}
  \jW^+(z)=\MFFminus{h-1,3p} (\Jplus_0)^h\SymmF_{h}(z)
  =   h!
  \MFFminus{h-1,3p}
  e^{2 h \omega_2(z) + 2 h \omega_3(z)},
\end{equation*}
which is in the grade with charge $2 h-1$ and at the level $3p^2(j+1)$
relative to the top.  Also, the left corner in
diagram~\eqref{s:triplet-easy} is a visualization of the fact that
$(\Jplus_0)^h\SymmF_{h}$ also has a vanishing singular vector
$\MFFplus{2 h+1,1}=(\Jminus_0)^{2p(j+1)+3}$; in the corresponding
Verma module, this state has an $\MFFplus{h-1,3p+1}$ singular vector,
located at the charge $-2 h$ and on the same level $3p^2(j+1)$
relative to the top as $\jW^+(z)$.  This same singular vector is
``seen'' from $\jW^+(z)$, and its vanishing makes $\jW^+(z)$ the
rightmost element in a $(4 h-1)$-plet under the zero-mode $s\ell(2)$.

We have thus found a $(4 h-1)$-plet at the absolute level (Sugawara
dimension) $\frac{h(h+1)}{k+2} + 3p^2(j+1) = 2 p (2 p(j + 1) + 1)$.

\begin{figure}[tb]
  \centering
  \begin{equation*}
    \xymatrix@C15pt@R24pt{
      &&&\verma="L"&\Circ="V"&&\Circ="U"&\iverma="R"&&&\nothing="up1"&\nothing="up2"\\
      {}\\
      &&\Circ="coL"&&&&&&\Circ="coR"&&\nothing="down1"&\\
      {}\\
      {}\\
      \Circ="iL"&\ivermaiv="fromR"&&&&&&&&\vermaiv="fromL"&\Circ="iR"&\nothing="down2"
      \ar@{{<}{.}{>}}|{2 p}"up1";"down1"
      \ar@{{<}{.}{>}}|{(3  h - 1)p}"up2";"down2"
      \ar@*{[|(2)]}@{-}"L";"L"-<80pt,0pt>
      \ar@*{[|(2)]}@{-}"R";"R"+<100pt,0pt>
      \ar@*{[|(2)]}@{-}"V"+<0pt,1pt>;"U"+<0pt,1pt>
      \ar@*{[|(1)]}"V"+<0pt,1pt>;"L"
      \ar@*{[|(1)]}"U"+<0pt,1pt>;"R"
      \ar|(.7){\MFFminus{3h-1,p}}"L";"fromL"+<-1pt,2pt>
      \ar|(.7){\MFFplus{3h-1,p;1}}"R";"fromR"+<1pt,2pt>
      \ar@*{[|(4)]}@{-}"fromR";"fromL"
      \ar@*{[|(1)]}"iL"+<0pt,1pt>;"fromR"
      \ar@*{[|(1)]}"iR"+<0pt,1pt>;"fromL"
      \ar@*{[|(2)]}@{-}"iL"+<0pt,1pt>;"iL"+<-20pt,1pt>
      \ar@*{[|(2)]}@{-}"iR"+<0pt,1pt>;"iR"+<20pt,1pt>
      \ar@{-}"coL"-<50pt,0pt>;"coR"+<70pt,0pt>
      \ar@{{-}{--}{>}}|{\MFFplus{1,2p+1}}"L"+<-1pt,1pt>;"coL"
      \ar@{{-}{--}{>}}|(.3){\MFFplus{h - 1, 3 p + 1}}"coL";"iL"
      \ar@{{-}{--}{>}}|{\MFFminus{1, 2 p + 1; 1}}"R"+<1pt,1pt>;"coR"
      \ar@{{-}{--}{>}}|(.3){\MFFminus{h - 1, 3 p + 1; 1}}"coR";"iR"
      \ar@{}|{\ \ \scriptstyle\Rstate_{-h}}"R";"R"+<6pt,16pt>
      \ar@{}|{\scriptstyle\Lstate_{-h}\ \ }"L";"L"+<6pt,16pt>
      \ar@{}|{\scriptstyle\jW^-}"fromL";"fromL"+<6pt,18pt>
      \ar@{}|{\scriptscriptstyle -h}"L";"L"+<5pt,-24pt>
      \ar@{}|{\scriptscriptstyle h}"R";"R"+<-5pt,-24pt>
      \ar@{}"fromR";"fromR"+<0pt,-30pt>|{-2 h + 1}
      \ar@{}"fromL";"fromL"-<0pt,30pt>|{2 h - 1}
      \ar@{}|{\scriptscriptstyle\qquad -h - 1}"coL";"coL"+<-9pt,-18pt>
      \ar@{}|{\scriptscriptstyle h + 1\quad\ }"coR";"coR"+<9pt,-18pt>
    }
  \end{equation*}
  \caption{\small $h = (j + 1) p + 1$ for $j\geq 0$.}
  \label{fig:s-MULTIPLET}
\end{figure}
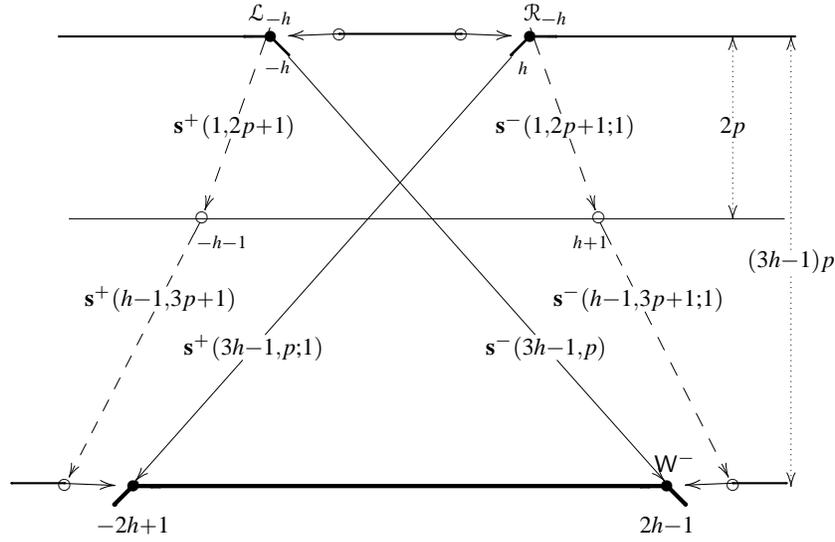

We next find a $(4 h-1)$-plet, at the same absolute level, in the
module whose top is shown in~\eqref{s:module-top}.  With the chosen
$h$, $\Lstate_{-h}$ has the vanishing singular vector
$\MFFplus{1,2p+1}$ and the nonvanishing one $\MFFminus{3 h-1,p}$.
These are shown in Fig.~\ref{fig:s-MULTIPLET}, where another crucial
piece is the
vanishing singular vector $\MFFplus{h-1,3p+1}$ and, importantly, the
``mirror images'' of all these singular vectors---twisted singular
vectors in the twisted module associated with $\Rstate_{-h}(z)$.  The
result is that
\begin{equation*}
  \jW^{-}(z) = \MFFminus{3 h-1, p}\Lstate_{-(j+1)p - 1}(z)
  \in\ker F_1\bbigcap\ker F_2
\end{equation*}
is the rightmost element of a $(4 h-1)$-plet, at the same
Sugawara dimension $2p(2p(j+1)+1)$ as $\jW^+(z)$.

In this case, there is also a nonvanishing singular vector in the
module associated with the unit operator:
\begin{equation*}
  \jW^0(z) = \MFFminus{2 h-1, 2p}\one(z)\in\ker F_1\bbigcap\ker F_2
\end{equation*}
(the arrangement of the relevant singular vectors is similar to the
one in Fig.~\ref{fig:s-MULTIPLET}, and we omit the details).

We repeat that each $\jW^{\pm,0}(z)$ is an $\hSL2$ primary state of
dimension $2p(2p(j+1)+1)$, is part of a $(4 (j+1) p+3)$-plet, and the
three of them, together with the $\hSL2$ currents, conjecturally
generate a $W$-algebra of local fields in $\ker F_1\bbigcap\ker F_2$.

\subsection{Relation between the symmetric and asymmetric
  realizations}\label{s:Wak}\addtocounter{subsubsection}{-1}
The symmetric and asymmetric realizations of extended algebras can be
mapped onto one another by a nonlocal field transformation.
Introducing it invokes the Wakimoto bosonization, and this deserves a
terminological comment.

\subsubsection{}
In CFT, representing anything as an exponential of free scalar(s) is
standardly called bosonization.  A typical example is a free-fermion
first-order system of fields $\eta(z)$ and $\xi(z)$ expressed as
\begin{equation*}
  \eta(z) = e^{-F(z)},\qquad
  \xi(z) = e^{F(z)},
\end{equation*}
where $F(z)$ is a scalar field with canonical normalization.  The
concept and the terminology have been extended to representing bosonic
first-order systems, with the OPE
\begin{equation*}
    \beta(z)\,\gamma(w)=\ffrac{-1}{z-w},
\end{equation*}
in terms of two scalars~\cite{[FMS]}, even though the procedure is a
map from bosons to bosons.  By the same token, the Wakimoto free-field
construction for $\hSL2$ currents is also commonly called bosonization
(although it involves by far not only exponentials).

\subsubsection{Wakimoto bosonization}\label{sec:Wak}
We recall the standard Wakimoto bosonization for $\hSL2_k$
currents \cite{[W],[FF]},\footnote{A standard notational atrocity
  committed in formulas~\eqref{JJJ-Wak} is the indiscriminate use of
  ``throughout'' symbols for $\hSL2$ currents; physicists tend to read
  such notation as the statement that the right-hand sides satisfy an
  $\hSL2$ algebra, of some level (in this case, $k$) deduced from the
  context.}
\begin{align}
  \notag
  \Jplus(z)&=-\beta(z),
  \\
  \label{JJJ-Wak}
  \Jnaught(z)&= \beta\gamma(z)
  +\ffrac{\sqrt{k+2}}{\sqrt{2}} \partial\Azero(z),
  \\
  \notag \Jminus(z)&= \beta\gamma\gamma(z) + k \partial\gamma(z)
  +\sqrt{2} \sqrt{k+2}\,\gamma\,\partial\Azero(z),
\end{align}
where $\beta(z)$ and $\gamma(z)$ constitute a first-order bosonic
system and $\Azero(z)$ is an independent canonically normalized free
boson,
\begin{equation*}
  \Azero(z)\,\Azero(w)=\log(z-w).
\end{equation*}

\subsubsection{Wakimoto $\myboldsymbol{\to}$ symmetric realization}
The $\beta$,\,$\gamma$,\,$\Dvarphi$ system of fields can be embedded
into the algebra of fields in~\bref{s-para} by a map $\maps{\cdot}$
such that
\begin{align*}
  \maps{\beta}(z)
  &=-\partial\Aone(z)e^{-\omega_1(z) + \omega_2(z) + 2 \omega_3(z)},
  \\
  \maps{\gamma}(z) &=e^{\omega_1(z) - \omega_2(z) - 2 \omega_3(z)},
  \\
  \maps{\partial\Azero}(z) &=
  \ffrac{\sqrt{2}}{\sqrt{k\!+\!2}}
  \Bigl(
  \ffrac{k\!+\!1}{k }\partial\Aone(z)
  -\ffrac{1}{k}\partial\Atwo(z)
  +\ffrac{k\!+\!2}{2 k}\partial\Athree(z)\Bigr)
  =\sqrt{2(k\!+\!2)} (\partial\omega_2(z) + \partial\omega_3(z)).
\end{align*}
It is worth noting that the $\beta\gamma$ current is then mapped as
\begin{equation*}
  \maps{\beta\gamma}(z) = -(k + 2)\partial\omega_2(z)
  - 2\partial\omega_3(z).
\end{equation*}

Under $\maps{\cdot}$, the currents~\eqref{JJJ-Wak} are mapped onto the
$\hSL2$ currents in~\bref{s-para}.  Moreover, the extended algebra
generators are also in the image of $\maps{\cdot}$.  Indeed, the
states at the left and right corners in~\eqref{s:triplet-easy} are
respectively given by
\begin{equation*}
  \maps{\gamma(z)^{2h}\,e^{h\frac{\sqrt{2}}{\sqrt{k + 2}} \Azero(z)}}
  \qquad\text{and}\qquad
  \maps{e^{h\frac{\sqrt{2}}{\sqrt{k + 2}} \Azero(z)}},\qquad h>0.
\end{equation*}
For $h=1$, these are respectively $(\Jminus_0)^2\W^{+}(z)$ and
$\W^{+}(z)$, Eq.~\eqref{s:wplus-x3}.  Next, the highest-weight state
in~\bref{s-even-more} is readily seen to be expressed as
\begin{equation*}
  \Lstate_{-h}(z) = \maps{e^{-h\frac{\sqrt{2}}{\sqrt{k + 2}} \Azero(z)}},
\end{equation*}
and hence the $\MFFminus{2,p}$ singular vector evaluates on this state
in the $\beta$,\,$\gamma$,\,$\Dvarphi$ theory.  The same applies to
$\jW^+(z)$, $\jW^0(z)$, and $\jW^-(z)$ in~\bref{sec:s-j}.  The twisted
highest-weight state $\Rstate_{-h}(z)$ is not expressible in terms of
the Wakimoto bosonization ingredients for $h>0$, and hence a ``half''
of the diagrams in Figs.~\ref{fig:s-TRIPLET} and~\ref{fig:s-MULTIPLET}
is not expressible in the Wakimoto bosonization, but the $\W^{-}(z)$
field (and, for $j\geq 0$, $\jW^{-}(z)$)~is.

To construct the inverse map, we have to express the three currents
$\partial\Aone(z)$, $\partial\Atwo(z)$, and
$\partial\Athree(z)$---which we temporarily denote as
$\partial\Aone^{\symm}(z)$, $\partial\Atwo^{\symm}(z)$, and
$\partial\Athree^{\symm}(z)$---in terms of three currents in the
$\beta$,\,$\gamma$,\,$\Dvarphi$ theory.  Two of these are
$\beta\gamma(z)$ and $\Dvarphi(z)$, and the third is the $\eta\xi(z)$
current ``hidden'' inside the $\beta$,\,$\gamma$ system~\cite{[FMS]}:
\begin{align}
  \notag
  \partial\Aone^{\symm}(z)&=\maps{\eta\xi(z)},
  \\
  \label{rho-s}
  \partial\Atwo^{\symm}(z)&=
  \maps{(k + 2) \beta\gamma(z) + \sqrt{2(k + 2)} \Dvarphi(z) + (k + 1)
    \eta\xi(z)},
  \\
  \partial\Athree^{\symm}(z)&=
  \maps{2 \beta\gamma(z) + \sqrt{2 (k + 2)} \Dvarphi(z)}.
  \notag
\end{align}

\subsubsection{Wakimoto $\myboldsymbol{\to}$ asymmetric realization}
The $\beta$,\,$\gamma$,\,$\Dvarphi$ fields are mapped into the algebra
of fields in~\bref{a-para} as
\begin{align*}
  \mapa{\beta}(z) &=-e^{\omega_1(z) + 2\omega_3(z)},
  \\
  \mapa{\gamma}(z) &=(\half\partial\Aii(z) + \half\partial\Afour(z))
  e^{-\omega_1(z) - 2\omega_3(z)},
  \\
  \mapa{\partial\Azero}(z) &=
  \ffrac{\sqrt{k+2}}{\sqrt{2}}\Bigl(
  \ffrac{1}{k+2}\,\partial\Aii(z)
  +\ffrac{1}{k}\,\partial\Athree(z)
  +\ffrac{1}{k}\,\partial\Afour(z)\Bigr)
  =\sqrt{2 (k + 2)} (\partial\omega_2(z) + \partial\omega_3(z)).
\end{align*}
Also,
\begin{equation*}
  \mapa{\beta\gamma}(z)
  = -\ffrac{1}{k}\,\partial\Athree(z)
  -\ffrac{k+2}{2 k}\,\partial\Afour(z)
  -\ffrac{1}{2}\,\partial\Aii(z)
  = -(k + 2) \partial\omega_2(z) - 2 \partial\omega_3(z).
\end{equation*}

We omit the maps between relevant vertices, and only give the
relations that allow inverting~$\mapa{\cdot}$, where we temporarily
write $\Aone^{\asymm}(z)$, $\Atwo^{\asymm}(z)$, and
$\Athree^{\asymm}(z)$ for the fields introduced in~\bref{aVir}.  Then
\begin{align}
  \notag
  \partial\Aone^{\asymm}(z)&=-\mapa{\eta\xi}(z),
  \\
  \label{rho-a}
  \partial\Atwo^{\asymm}(z)&=
  \mapa{(k + 2) \beta\gamma(z) + \sqrt{2 (k + 2)} \Dvarphi(z)
    + (k + 2) \eta\xi(z)},
  \\
  \partial\Athree^{\asymm}(z)&=
  \mapa{2 \beta\gamma(z)  + \sqrt{2 (k + 2)} \Dvarphi(z)}.
  \notag
\end{align}

\subsubsection{}
The map between the symmetric and asymmetric realizations is
\begin{equation*}
  \mapa{\cdot}\ccirc\mapsinv{\cdot},
\end{equation*}
which is a ``highly nonlocal change of variables'' in the two-boson
space (the $\Athree(z)$ are the same in~\eqref{rho-s}
and~\eqref{rho-a}).  We note that $k$ in~\eqref{rho-s}
and~\eqref{rho-a} and the related formulas is a common parameter, the
$\hSL2$ level.  It is noninteger in our setting, and therefore the
screenings $F=\oint e^{\Atwo^{\asymm}}$ and $F_2=\oint
e^{\Atwo^{\symm}}$ are undefined in intrinsic $\beta$,\;$\gamma$
terms.  If one starts with the Wakimoto realization of $\hSL2$ and
maps it by $\mapa{\cdot}$ or $\maps{\cdot}$, then the appearance of
these screenings is ``accidental.''


\section{\textbf{Double bosonization of $\bNich\myboldsymbol{(X)}$}}
There is a functorial, vector-space-preserving correspondence between
multivertex Yet\-ter--Drinfeld $\Nich(X)$ modules and modules over a
nonbraided Hopf algebra $\algU$, the double bosonization
of~$\Nich(X)$.  In this section, we show how $\algU$ can be
constructed in general (\bref{Nich-dual}--\bref{double-B}) and
evaluate it for our two Nichols algebras (\bref{sec:algUa}
and~\bref{sec:algUs}).  The two resulting $\algU$ compare nicely, as
we show in~\bref{sec:iso}.

\subsection{The dual Nichols algebra.}\label{Nich-dual}
For a Nichols algebra $\Nich(X)$, we let $\Nich(X)^*$ denote its
graded dual. The reader can consider $X$ an object in a
rigid braided category; we do not go into the details of (standard)
axioms and simply assume all the necessary structures to exist.  The
first of these is the evaluation $\eval{~}{}:\Nich(X)^*\tensor
\Nich(X)\to k$, which is diagrammatically denoted by
\raisebox{-4pt}{$\
\begin{tangles}{l}
  \ev
\end{tangles}
\ $}. \ By abuse of notation, similar diagrams are used to represent
the restriction of $\eval{~}{}$ to $X^*\tensor X$; it is extended to
tensor products as \raisebox{-4pt}{$\
\begin{tangles}{l}
  \ev\step[-1.5]\hev
\end{tangles}
\ \;$}, and so on.  The product, coproduct, and braiding in
$\Nich^*=\Nich(X)^*$ are defined by the respective rules
\begin{equation}\label{respective}
  \begin{tangles}{l}
    \fobject{\Nich^*}\step[2]\fobject{\Nich^*}\step[1]\fobject{\ \ \Nich}\\
    \cu\step[1]\id\\
    \step[1]\ev
  \end{tangles}
  \ \ = \ \ 
  \begin{tangles}{l}
    \id\step[1]\id\step[1]\hcd\\
    \Ev\hev
  \end{tangles}
  \qquad\qquad
  \begin{tangles}{l}
    \hcd\step[1]\id\step[1]\id\\
    \Ev\hev
  \end{tangles}
  \ \ = \ \
  \begin{tangles}{l}
    \id\step[1]\cu\\
    \ev
  \end{tangles}
  \qquad\qquad
  \begin{tangles}{l}
    \hx\step[1]\id\step[1]\id\\
    \Ev\hev
  \end{tangles}
  \ \ = \ \
  \begin{tangles}{l}
    \id\step[1]\id\step[1]\hx\\
    \Ev\hev
  \end{tangles}
\end{equation}

Every (left--left) Yetter--Drinfeld $\Nich(X)$-module carries a left
action of $\Nich(X)^*$ defined standardly using the
$\Nich(X)$-coaction, as $\alpha y=\eval{\alpha}{y\bmone}y\bzero$,
which for one-vertex modules (and, in fact, for all multivertex
Yetter--Drinfeld modules, cf.~\cite{[STbr]}) reduces to the coproduct
on~$\Nich(X)$:
\begin{equation}\label{dual-acts}
  \begin{tangles}{l}
    \vstr{200}
    \lu
    \step[1]\id\object{\kern-10pt\rule[18pt]{12pt}{4pt}}\step[.1]\id
  \end{tangles}\ \
  = \ \ 
  \begin{tangles}{l}
    \id\step[1]\hcd\step[1]\id\step[.1]\id\\
    \hev\step[1]\id\step[1]\id\step[.1]\id
  \end{tangles}
\end{equation}
On graded components, this gives the maps
\begin{equation*}
  \mathsf{e}^{(r)}_s : {X^*}^{\otimes r}\tensor X^{\otimes s}\tensor
  Y\to X^{\otimes(s-r)}\tensor Y,\qquad s\geq r,
\end{equation*}
such that
\begin{equation}\label{-r-act}
  \begin{tangles}{l}
    \fobject{X^*}\step[1]\fobject{\ \ X^*}\step[3]\fobject{X}\step[1]\fobject{X}\step[1]\fobject{\ X}
    \step[3]\fobject{\ X}\step\fobject{\ \ Y}
    \\[-1pt]
    \vstr{160}\id\step\id\step\raisebox{12pt}{$\tensor$}\step
    \id\step\id\step\id\step\raisebox{12pt}{\dots\ \ }\id\step\id\step[.1]\id
    \raisebox{12pt}{\ \ $\to$\ }
    \step\ev\step[-1.5]\hev\step\id\step\raisebox{12pt}{\dots\ \ }\id\step\id\step[.1]\id
  \end{tangles}
\end{equation}
for $r=2$ and totally similarly for all $r\leq s$.

We now discuss the nomenclature of tensor products ${X^*}^{\otimes
  r}\tensor X^{\otimes s}\tensor Y$.  \ We here label the tensor
factors as $(-r,\dots,-1,1,\dots,s;s+1)$ (the semicolon separates
$Y$).  The pairing such as in the right-hand side of~\eqref{-r-act} is
then conveniently written as $\rho_{(-r,r)}$; we slightly abuse the
notation by writing $\rho_{(-\ell,\ell)}: {X^*}^{\otimes r}\tensor
X^{\otimes s}\tensor Y\to {X^*}^{\otimes(r-\ell)}\tensor
X^{\otimes(s-\ell)}\tensor Y$ (with $r\geq\ell$ and $s\geq\ell$) also
for the map that should be more rigorously written as
$\id^{\otimes(r-\ell)}\tensor\rho_{(-\ell,\ell)}\tensor\id^{\otimes(s-\ell)}$.
The ``leg notation'' for elements of the braid group extends to
negative labels according to the pattern
\begin{alignat*}{2}
  \Psi_{-2}&=\id^{\otimes(r-3)}\tensor\Psi\tensor\id\tensor\id^{\otimes
    s} &&:{X^*}^{\otimes r}\tensor{X}^{\otimes s} \to{X^*}^{\otimes
    r}\tensor{X}^{\otimes s},
  \\
  \Psi_{-1}&=\id^{\otimes(r-2)}\tensor\Psi\tensor\id^{\otimes s}
  &&:{X^*}^{\otimes r}\tensor{X}^{\otimes s} \to{X^*}^{\otimes
    r}\tensor{X}^{\otimes s},
  \\
  \intertext{and, ``in the middle,''}
  \Psi_{0}&=\id^{\otimes(r-1)}\tensor\Psi\tensor\id^{\otimes(s-1)}
  &&:{X^*}^{\otimes r}\tensor{X}^{\otimes s}
  \to{X^*}^{\otimes(r-1)}\tensor X\tensor X^*\tensor{X}^{\otimes(s-1)}
\end{alignat*}
(and, of course, $\Psi_{1}=\id^{\otimes
  r}\tensor\Psi\tensor\id^{\otimes(s-2)}$, and so on).

We let $\mathsf{f}_s:X\tensor X^{\otimes s}\tensor Y\to
X^{\otimes(s+1)}\tensor Y$ denote the map such that for any $x\in X$
and $y\in X^{\otimes s}\tensor Y$, $\mathsf{f}_s(x,y) = x\adjoint y$,
the adjoint action by $x$.  We keep the same notation for the adjoint
action map also in the case where some $X^*$ factors precede the $X$
in the tensor product and the rigorous writing should be $\id^{\otimes
  r}\tensor\mathsf{f}_s:{X^*}^{\otimes r}\tensor X\tensor X^{\otimes
  s}\tensor Y\to {X^*}^{\otimes r}\tensor X^{\otimes(s+1)}\tensor Y$.

\subsection{``Commutator'' identities}
We now see how the left adjoint action of $\Nich(X)$ and the above
left action of $\Nich(X)^*$ commute with each other.

\subsubsection{}
We first recall a ``commutator'' identity for the maps $X^*\tensor
X^{\otimes(s+1)}\tensor Y\to X^{\otimes s}\tensor Y$ effected by
$\mathsf{e}^{(1)}$ and $\mathsf{f}$ (with appropriate subscripts, to
be restored momentarily)~\cite{[STbr]}.  We first show the identity in
graphic form, with $s=2$ for definiteness:
\begin{equation}\label{commutator-first}
  \begin{tangles}{l}
    \step[2.75]\id\\[-20pt]
    \vstr{200}
    \id\step[1]
    \lu\object{\raisebox{18pt}{\kern-4pt\tiny$\blacktriangleright$}}
    \step[.8]\object{\rule[18pt]{18pt}{4pt}}
    \step[.8]\id\step[.1]\id\\[-4.5pt]
    \ev\step[1.6]\id\step[.1]\id\\[-33.5pt]
    \vstr{167}\step[2.5]\id\step[.6]\id
  \end{tangles}\ \ \
  - \ \
  \begin{tangles}{l}
    \vstr{50}\id\step[1]\hx\step[1]\id\step[1.6]\id\step[.1]\id\\
    \vstr{120}
    \id\step[1]\id\step[1]
    \lu\object{\raisebox{10pt}{\kern-4pt\tiny$\blacktriangleright$}}
    \step[.8]\object{\rule[10pt]{18pt}{4pt}}
    \step[.8]\id\step[.1]\id\\[-12pt]
    \vstr{60}\step[3.75]\id\\[-.5pt]
    \hev\step[2]\id\step[.75]\id\step[.85]\id\step[.1]\id
  \end{tangles}\ \ \
    = \ \
    \begin{tangles}{l}
      \vstr{260}\ev\step[1]\id\step[.5]\id\step[.5]\id\step[.1]\id
    \end{tangles}\ \
    - \ \
    \begin{tangles}{l}
      \vstr{50}\id\step[1]\hx\step[1]\id\step[1]\id\step[.1]\id\\
      \vstr{50}\id\step[1]\id\step[1]\hx\step[1]\id\step[.1]\id\\
      \vstr{50}\id\step[1]\id\step[1]\id\step[1]\hx\\
      \vstr{50}\id\step[1]\id\step[1]\id\step[1]\hx\\
      \vstr{50}\id\step[1]\id\step[1]\hx\step[1]\id\step[.1]\id\\
      \vstr{50}\id\step[1]\hx\step[1]\id\step[1]\id\step[.1]\id\\
      \hev\\[-20pt]
      \vstr{33}\step[2]\id\step[1]\id\step[1]\id\step[.1]\id
    \end{tangles}
\end{equation}
Each diagram is a map $X^*\tensor X\tensor X^{\otimes 2}\tensor Y\to
X^{\otimes 2}\tensor Y$. \ In the first diagram in the left-hand side,
the adjoint action by $x\in X$ is applied first, and is followed by
the action of an element of~$X^*$; in the second diagram, the order is
reversed, at the expense of a braiding.  The second diagram in the
left-hand side can of course be rewritten as
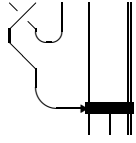
\begin{equation*}
  \begin{tangles}{l}
    \vstr{50}\hxx\step[1]\id\step[1]\id\step[1.5]\id\step[.1]\id\\
    \step[1]\hev\\[-20pt]
    \vstr{50}\dh\step[2.5]\id\step[1.5]\id\step[.1]\id\\
    \vstr{50}\step[.5]\hd\step[2]\id\step[1.5]\id\step[.1]\id\\
    \step[1]\lu[2]\object{\raisebox{8pt}{\kern-4pt\tiny$\blacktriangleright$}}\step[.8]\object{\rule[8pt]{18pt}{4pt}}\step[.7]\id\step[.1]\id\\[-10pt]
    \vstr{50}\step[3.8]\id
  \end{tangles}
\end{equation*}
where $\Psi_0^{-1}:X^*\tensor X\to X\tensor X^*$ appears.  

The general (and ``analytic'') form of~\eqref{commutator-first} for
maps $X^*\tensor X^{\otimes(s+1)}\tensor Y\to X^{\otimes s}\tensor Y$
is~\cite{[STbr]}
\begin{equation}\label{commutator-1}
  \mathsf{e}^{(1)}_{s+1}\ccirc\mathsf{f}_s -
  \mathsf{f}_{s-1}\ccirc
  \mathsf{e}^{(1)}_s\ccirc
  \Psi^{-1}_{0}
  = \rho_{(-1,1)}
  \tensor\id^{\tensor(s+1)}
  - \KK(s+1),
\end{equation}
where $\KK(s+1)$ is the \textit{monodromy operation}
\begin{equation}\label{K2}
  \!\!\!
  \KK(s+1)
  =(\rho_{(-1,1)}
  \tensor\id^{\otimes(s+1)})
  \ccirc(\Psi_1\dots\Psi_{s+1}\Psi_{s+1}\dots\Psi_1)
  \ : \
  \begin{tangles}{l}
    \fobject{X^*}\step[2]\fobject{X}\step[1]\fobject{X}\step\fobject{X}
    \step[3]\fobject{X}\step\fobject{Y}\\
    \vstr{62}\id\step[2]\hx\step[1]\id\step[3]\id\step[1]\id\step[.1]\id\\
    \vstr{62}\id\step[2]\id\step[1]\hx\step[3]\id\step[1]\id\step[.1]\id\\
    \vstr{100}\id\step[2]\id\step[1]\id\step[2]\object{\ \dots\ }\step[2]
    \id\step\id\step[.1]\id\\
    \vstr{62}\id\step[2]\id\step[1]\id\step[3]\hx\step[1]\id\step[.1]\id\\
    \vstr{62}\id\step[2]\id\step[1]\id\step[3]\id\step\hx\\
    \vstr{62}\id\step[2]\id\step[1]\id\step[3]\id\step\hx\\
    \vstr{62}\id\step[2]\id\step[1]\id\step[3]\hx\step\id\step[.1]\id\\
    \vstr{100}\id\step[2]\id\step[1]\id\step[2]\object{\ \dots\ }\step[2]
    \id\step\id\step[.1]\id\\
    \vstr{62}\id\step[2]\id\step\hx\step[3]\id\step\id\step[.1]\id\\
    \vstr{62}\id\step[2]\hx\step\id\step[3]\id\step\id\step[.1]\id\\
    \vstr{62}\ev\step\id\step\id\step[3]\id\step\id\step[.1]\id\\
    \step[3]\fobject{2}\step\fobject{3}\step[3]\fobject{s+1\quad}
    \step\fobject{\quad s+2}
  \end{tangles}
\end{equation}

\subsubsection{}\label{sec:-r-act}
Straightforward calculation shows that identity~\eqref{commutator-1}
generalizes to a ``commutator'' of $\mathsf{f}$ with
$\mathsf{e}^{(r)}$ as follows (both sides are maps ${X^*}^{\otimes
  r}\tensor X^{\otimes(s+1)}\tensor Y\to X^{\otimes(s+1-r)}\tensor
Y$):
\begin{multline}\label{commutator-rs}
  \mathsf{e}^{(r)}_{s+1}\ccirc\mathsf{f}_s -
  \mathsf{f}_{s-r}\ccirc
  \mathsf{e}^{(r)}_s\ccirc
  \Psi_{-r+1}^{-1}\dots\Psi_0^{-1}  
  \\
  \shoveleft{{}=(\rho_{(-r,r)}\tensor\id^{\otimes(s-r+2)})
    \ccirc
    \bigl(\id+\Psi_{-1}+\Psi_{-1}\Psi_{-2}+\dots
    +\Psi_{-1}\dots\Psi_{-(r-1)}\bigr)}
  \\
  -\KK(s-r+2)\ccirc(\id^{\otimes 2}\tensor\rho_{(-r+1,r-1)})
  \ccirc\Psi^{-1}_{-r+2}\dots\Psi^{-1}_{0}
  \\
  {}\ccirc\bigl(\id+\Psi_{-r+1} + \Psi_{-r+1}\Psi_{-r+2}
  +\dots+\Psi_{-r+1}\Psi_{-r+2}\dots\Psi_{-1}\bigr).
\end{multline}
Here,
$\Psi_{-r+1}^{-1}\dots\Psi_0^{-1}:{X^*}^{\otimes r}\tensor X\to
X\tensor {X^*}^{\otimes r}$; for $r=3$, for example, this is the map
$\ \begin{tangles}{l}
  \vstr{25}\hstr{50}\id\step[1]\id\step[1]\hxx\\
  \vstr{25}\hstr{50}\id\step[1]\hxx\step[1]\id\\
  \vstr{25}\hstr{50}\hxx\step[1]\id\step[1]\id
\end{tangles}\ $.  The left-hand side can also be written as
$\mathsf{e}^{(r)}_{s+1}\ccirc\mathsf{f}_s - \mathsf{f}_{s-r}\ccirc
\mathsf{e}^{(r)}_s\ccirc \Psi_r\dots\Psi_1$.

To continue with the $r=3$ example, we write the right-hand side
of~\eqref{commutator-rs} explicitly for $r=3$ and $s=3$:
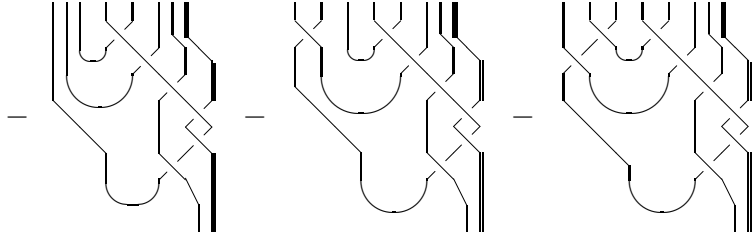
\begin{equation}\label{commutator-rs-graph}
  \begin{tangles}{l}
    \vstr{33}\id\step[.5]\id\step[.5]\id\step[1]\id\step[.5]\id\step[.5]\id\step[.5]\id\step[.5]\id\step[.1]\id\\
    \Ev\step[-.5]\ev\step[-1.5]\hev\step[1]\step[.5]\id\step[.5]\id\step[.1]\id\\
    \vstr{150}\step[3]\step[.5]\id\step[.5]\id\step[.1]\id
  \end{tangles}
  \ \ + \ \
  \begin{tangles}{l}
    \vstr{33}\id\step[.5]\id\step[.5]\id\step[1]\id\step[1]\id\step[.5]\id\step[.5]\id\step[.5]\id\step[.1]\id\\
    \vstr{50}\id\step[.5]\id\step[.5]\id\step[1]\hx\step[.5]\id\step[.5]\id\step[.5]\id\step[.1]\id\\
    \step[-.25]{\makeatletter\@ev{0,\hm@de}{20,\hm@detens}{35}b\makeatother}
    {\makeatletter\@ev{0,\hm@de}{20,\hm@detens}{25}b\makeatother}\step[1.25]\hev\step[2]\id\step[.5]\id\step[.1]\id\\
    \step[3.5]\step[.5]\id\step[.5]\id\step[.1]\id
  \end{tangles}
  \ \ + \ \
  \begin{tangles}{l}
    \vstr{33}\id\step[.5]\id\step[.5]\id\step[1]\id\step[1]\id\step[1]\id\step[.5]\id\step[.5]\id\step[.1]\id\\
    \vstr{50}\id\step[.5]\id\step[.5]\id\step[1]\hx\step[1]\id\step[.5]\id\step[.5]\id\step[.1]\id\\
    \vstr{50}\id\step[.5]\id\step[.5]\id\step[1]\id\step[1]\hx\step[.5]\id\step[.5]\id\step[.1]\id\\
    {\makeatletter\@ev{0,\hm@de}{20,\hm@detens}{40}b\makeatother}
    \step[-.25]{\makeatletter\@ev{0,\hm@de}{20,\hm@detens}{25}b\makeatother}\step[-.25]{\makeatletter\@ev{0,\hm@de}{20,\hm@detens}{10}b\makeatother}
    \step[5]\id\step[.5]\id\step[.1]\id\\
    \vstr{33}\step[4]\step[.5]\id\step[.5]\id\step[.1]\id
  \end{tangles}\kern200pt
\end{equation}
\begin{equation*}
\kern100pt\ \ \ - \ \ \
  \begin{tangles}{l}
    \vstr{33}\id\step[.5]\id\step[.5]\id\step[1]\id\step[1]\id\step[1]\id\step[.5]\id\step[.5]\id\step[.1]\id\\
    \vstr{33}\step[5]\id\step[.1]\id\\[-6.3pt]
    \vstr{50}\id\step[.5]\id\step[.5]\id\step[1]\hx\step[1]\id\step[.5]\dh\dh\\
    \vstr{50}\id\step[.5]\id\step[2.5]\hx\step[1]\id\step[.5]\hd\\[-10pt]
    \step[1]\hev\step[4.1]\vstr{70}\id\\[-10pt]
    \vstr{50}\id\step[4]\hx\step[1]\id\\[-10pt]
    \step[-.25]{\makeatletter\@ev{0,\hm@de}{20,\hm@detens}{25}b\makeatother}\\[-10pt]
    \vstr{50}\d\step[3]\id\step[1]\hx\\
    \vstr{50}\step[1]\d\step[2]\id\step[1]\hx\\
    \vstr{50}\step[2]\id\step[2]\hx\step[1]\id\step[.1]\id\\[-.5pt]
    \step[2]\ev\step[1]\hd\step[.5]\id\step[.1]\id
  \end{tangles}
  \ \ \ - \ \ \
  \begin{tangles}{l}
    \vstr{33}\id\step[1]\id\step[1]\id\step[1]\id\step[1]\id\step[1]\id\step[.5]\id\step[.5]\id\step[.1]\id\\
    \vstr{33}\step[6]\id\step[.1]\id\\[-6.3pt]
    \vstr{50}\hx\step[1]\id\step[1]\hx\step[1]\id\step[.5]\dh\dh\\
    \vstr{50}\id\step[1]\id\step[3]\hx\step\id\step[.5]\hd\\[-10pt]
    \step[2]\hev\step[4.1]\vstr{75}\id\\[-10pt]
    \vstr{50}\dh\step[4.5]\hx\step[1]\id\\[-10pt]
    \step[1]\Ev\\[-10pt]
    \vstr{50}\step[.5]\d\step[3.5]\id\step[1]\hx\\
    \vstr{50}\step[1.5]\d\step[2.5]\id\step[1]\hx\\
    \vstr{50}\step[2.5]\id\step[2.5]\hx\step[1]\id\step[.1]\id\\
    \step[1.75]{\makeatletter\@ev{0,\hm@de}{20,\hm@detens}{25}b\makeatother}\step[4.25]\hd\step[.5]\id\step[.1]\id
  \end{tangles}
  \ \ \ - \ \ \
  \begin{tangles}{l}
    \vstr{33}\id\step[1]\id\step[1]\id\step[1]\id\step[1]\id\step[1]\id\step[.5]\id\step[.5]\id\step[.1]\id\\
    \vstr{33}\step[6]\id\step[.1]\id\\[-6.3pt]
    \vstr{50}\id\step[1]\hx\step[1]\hx\step[1]\id\step[.5]\dh\dh\\
    \vstr{50}\hx\step[3]\hx\step[1]\id\step[.5]\hd\\[-10pt]
    \step[2]\hev\step[4.1]\vstr{70}\id\\[-10pt]
    \vstr{50}\id\step[5]\hx\step[1]\id\\[-10pt]
    \step[.5]{\makeatletter\@ev{0,\hm@de}{20,\hm@detens}{30}b\makeatother}\\[-10pt]
    \vstr{50}\d\step[4]\id\step[1]\hx\\
    \vstr{50}\step[1]\d\step[3]\id\step[1]\hx\\
    \vstr{50}\step[2]\hd\step[2.5]\hx\step[1]\id\step[.1]\id\\
    \step[1.75]{\makeatletter\@ev{0,\hm@de}{20,\hm@detens}{25}b\makeatother}\step[4.25]\hd\step[.5]\id\step[.1]\id
  \end{tangles}
\end{equation*}

\noindent
with the two lines corresponding to the two terms in the right-hand
side of~\eqref{commutator-rs}.

\subsubsection{Diagonal braiding and a
  $\myboldsymbol{T(X^*)\tensor}\bNich\myboldsymbol{(X)\tensor
    k\bGamma}$ algebra}
We next specify~\eqref{commutator-rs} to the case of diagonal
braiding, defined by~\eqref{qij} on $\Nich(X)$ generators.  First, for
the basis $(\dF_i)$ in $X^*$ such that
$\eval{\dF_i}{F_j}=\delta_{i,j}$, we apply the right
diagram in~\eqref{respective} to deduce that the $q_{\cdot,\cdot}$ in
$\Psi:\dF_i\tensor F_j\mapsto q_{\wh i,j} F_j\tensor\dF_i$ and
$\Psi:\dF_i\tensor\dF_j\mapsto q_{\wh i,\wh j}\dF_j\tensor\dF_i$ are
given by
\begin{equation}\label{q-hat}
  q_{\wh i,j} = q_{j,i}^{-1}\quad\text{and}\quad q_{\wh i,\wh j} = q_{i,j}.
\end{equation}
Then the monodromy operation evaluates on generators as
\begin{multline}
  \KK(s+1):
  \dF_i\tensor F_{j_1}\tensor F_{j_2}\tensor\dots\tensor F_{j_s}\tensor y
  \\
  \mapsto
  \eval{\dF_i}{F_{j_1}} (q_{j_1,j_2}q_{j_2,j_1})\dots(q_{j_1,j_s}q_{j_s,j_1})
  (q_{j_1,y}q_{y,j_1})
  F_{j_2}\tensor\dots\tensor F_{j_s}\tensor y
\end{multline}
for a homogeneous $y\in Y$.
  
Let $\bGamma$ be an Abelian group with generators $\bK_j$,
$j=1,\dots,\theta$, such that their action on $\Nich(X)\tensor Y$,
interpreted as adjoint action, produces the monodromies as in the last
formula:
\begin{equation}\label{KF}
  \bK_i F_j \bK_i^{-1} = q_{i,j} q_{j,i} F_j
\end{equation}
(and $\bK_i y \bK_i^{-1} = q_{i,y} q_{y,i} y$).  Commutator identity
\eqref{commutator-1} then becomes
\begin{equation}\label{dFF}
  \dF_i F_j - q_{j,i} F_j \dF_i =
  \delta_{i,j}(1 - \bK_j).
\end{equation}
Also setting
\begin{equation}\label{KdF}
  \bK_i \dF_j \bK_i^{-1} = q_{i,j}^{-1}q_{j,i}^{-1} \dF_j,
\end{equation}
we obtain an associative algebra on $T(X^*)\tensor\Nich(X)\tensor
k\bGamma$ with relations given by those in $\Nich(X)$
and~\eqref{KF}--\eqref{KdF}.

Formula \eqref{commutator-rs} becomes the statement that for each
$F_j$, the map $T(X^*)\to T(X^*)\tensor k\bGamma$ defined as
\begin{equation}
  (\dF_{i_r}\dots\dF_{i_1})\rightder_j
  = \dF_{i_r}\dots\dF_{i_1} F_j -
  q_{j,i_r}\dots q_{j,i_1} F_j \dF_{i_r}\dots\dF_{i_1}
\end{equation}
(in particular, $\dF_{i}\rightder_j = \delta_{i,j}(1 - \bK_j)$) is a
braided right derivation:
\begin{equation}
  (\dF_{i_r}\dots\dF_{i_1})\rightder_j
  = \dF_{i_r}\dots\dF_{i_2}(\dF_{i_1}\rightder_j)
  + q_{j,i_1} \bigl((\dF_{i_r}\dots\dF_{i_2})\rightder_j\bigr)\dF_{i_1}.
\end{equation}

\subsubsection{The ``half-way algebra''
  $\bNich\myboldsymbol{(X)^*\tensor{}}\bNich\myboldsymbol{(X)\tensor
    k\bGamma}$}
We verify in the particular cases studied in this paper that the
associative algebra structure on $T(X^*)\tensor\Nich(X)\tensor
k\bGamma$ pushes forward to $\Nich(X^*)\tensor\Nich(X)\tensor
k\bGamma$, i.e., to the quotient $\Nich(X^*)=T(X^*)/\mathscr{I}^*$ in
the first factor (the general proof must be possible by borrowing some
relevant steps from~\cite{[HS-double]}).  This algebra is only half
the way from $\Nich(X)$ to its double bosonization because the Abelian
group $\bGamma$ read off from monodromy, not braiding, is too coarse
to yield a Hopf algebra.

\subsection{The case of an $\bHHyd$ braiding}
We next assume that the braiding $\Psi$ is the one in the
Yetter--Drinfeld category $\HHyd$ for some Hopf algebra $H$
(see~\bref{app:HHyd}).  Monodromy---double braiding---then evaluates
as
\begin{equation*}
  \Psi^2:x\tensor y\mapsto
  x\mone{}' \,y\mone \hA(x\mone{}''')\leftii x\zero
  \tensor x\mone{}''\leftii y\zero.
\end{equation*}

From now on, we assume $H$ to be commutative and cocommutative.  Then
the last formula is simplified to
\begin{equation*}
  \Psi^2:x\tensor y\mapsto
  y\mone\leftii x\zero\tensor x\mone\leftii y\zero.
\end{equation*}
This is to be used in the definition of $\KK$, Eq.~\eqref{K2}.  We ask
when the right-hand side can be interpreted ``nonsymmetrically'' with
respect to $x$ and $y$, as an operation acting on $y$, not involving
any action on $x$.  This is the case if Sommerh\"auser's
condition~\cite{[Sommerh-deformed]}
\begin{equation}\label{Sommerh}
  u\mone\leftii v\tensor u\zero = v\zero\tensor v\mone\leftii u
\end{equation}
is imposed for all relevant $u$ and $v$ (we actually need the cases
where $u\tensor v$ is in $\Nich(X)\tensor\Nich(X)$, $\Nich(X)\tensor
Y$, and $Y\tensor\Nich(X)$.  With~\eqref{Sommerh}, indeed, the double
braiding takes the form
\begin{equation}\label{becomes}
  \Psi^2:x\tensor y\mapsto{}
  x\zero{}\zero\tensor x\mone x\zero{}\mone\leftii y
  =  x\zero\tensor x\mone{}'\, x\mone{}''\leftii y.
\end{equation}

We assume~\eqref{Sommerh} to hold from now on.  Before proceeding, we
make two brief remarks:
\begin{enumerate}
\item It suffices to impose~\eqref{Sommerh} on the generators, because
  if this condition holds for pairs $(y,x)$ and $(z,x)$, then it also
  holds for $(y\tensor z,x)$:
\begin{align*}
  (y\tensor z)\mone\leftii x\tensor (y\tensor z)\zero &=
  y\mone z\mone\leftii x\tensor y\zero\tensor z\zero
  \\
  &=y\mone\leftii x\zero\tensor y\zero\tensor x\mone\leftii z
  \\
  &=x\zero{}\zero\tensor x\zero{}\mone\leftii y\tensor x\mone\leftii z
  \\
  &=x\zero\tensor x\mone{}''\leftii y\tensor x\mone{}'\leftii z
  \\
  &=x\zero\tensor x\mone{}\leftii(y\tensor z),
\end{align*}
where the last equality is of course valid because $H$ is
cocommutative.

\item In terms of the braiding matrix entries,
  condition~\eqref{Sommerh} implies that $q_{i,j}=q_{j,i}$. \ In
  particular, for braiding matrices~\eqref{qij-asymm}
  and~\eqref{qij-symm}, this means that $\xx^2=1$, which brings us
  back to the braiding matrices in~\eqref{qij-first}.
\end{enumerate}

For a commutative and cocommutative $H$, \ $H=k\Gamma$
(see~\bref{app:commcocomm}), double braiding~\eqref{becomes}, with
$x=F_i$, becomes
\begin{equation*}
  \Psi^2:F_i\tensor y\mapsto{}
  F_i\tensor (g_i)^2\leftii y.
\end{equation*}
Also interpreting the $\Gamma$ action as an adjoint action (i.e.,
writing $g_i\leftii F_j=g_i F_j g_i^{-1}$), we see that $\bK_i$
in~\eqref{KF}--\eqref{KdF} are given by
\begin{equation*}
  \bK_i = (g_i)^2.
\end{equation*}
Hence, there is an associative algebra structure on $\algU =
\Nich(X^*)\tensor\Nich(X)\tensor k\Gamma$ with
$\Nich(X^*)\tensor\Nich(X)\tensor k\bGamma$ a subalgebra in it.

We summarize our findings as follows.

\subsection{Double bosonization: the Hopf algebra
  $\balgU$}\label{double-B}
\begin{it}%
  For a braided vector space $X\in\GGyd$ with a chosen basis $F_i$
  \textup{(}and the dual basis $\dF_i$ in $X^*$\textup{)} and a
  symmetric braiding matrix $(q_{i,j})$ in this basis, the algebra
  $\algU$ on generators $F_i$, $\dF_i$, $g_i$ $(i=1,\dots,\theta)$
  contains $\Nich(X^*)$ and $\Nich(X)$ as subalgebras and, in
  addition, has the relations
  \begin{align}
    \notag
    g_i F_j g_i^{-1} &= q_{i,j} F_j,
    \\
    \label{cross-commutator}
    \dF_i F_j - q_{j,i} F_j \dF_i &= \delta_{i,j}(1 - (g_j)^2),
    \\
    g_i \dF_j g_i^{-1} &= q_{i,j}^{-1} \dF_j,\qquad\qquad
    i,j=1,\dots,\theta.
    \notag
  \end{align}%
  Moreover, $\algU$ is a Hopf algebra, with the coproducts such that
  all $g_i\in\Gamma$, $i=1,\dots,\theta$, are group-like and
  \begin{equation*}
    \Delta:F_j\mapsto g_j\tensor F_j + F_j\tensor 1
    \quad  
    \Delta:\dF_j\mapsto g_j\tensor \dF_j + \dF_j\tensor 1,
  \end{equation*}
  and with the antipode
  \begin{gather*}
    S(F_i)=-g_i^{-1} F_i,\qquad
    S(\dF_i)=-g_i^{-1}\dF_i.
  \end{gather*}%
\end{it}%

The formula for $\Delta(F_j)=\Delta(F_j\Smash 1)$ is nothing but the
Radford formula for $\Nich(X)\Smash k\Gamma$.  Hence, in particular,
the relations in $\Nich(X)$ are compatible with the coproduct (the
corresponding ideal is a Hopf ideal).  The formula for
$\Delta(\dF_j)$, similarly, is the Radford formula for
$\Nich(X)^*\Smash k\Gamma$ with the $k\Gamma$ action and coaction
changed by composing each with the antipode; for a commutative
cocommutative Hopf algebra, this still gives a left action and a left
coaction.  It therefore only remains to verify that
cross-commutator~\eqref{cross-commutator} is compatible with the above
coproduct.  This is straightforward.

We also note that the $q$-commutator in~\eqref{cross-commutator} can be
conveniently ``straightened out'' by defining $\phi_i = g_i^{-1}
\dF_i$. \ Then
\begin{gather*}
  F_j \phi_i - \phi_i F_j = \delta_{i,j}(g_i - g_i^{-1}),
  \\[-4pt]
  \intertext{and we also have}
  \Delta(\phi_i)= 1\tensor\phi_i + \phi_i\tensor g_i^{-1},
  \qquad
  S(\phi_i) = -\phi_i g_i.
\end{gather*}

\subsection{The $\balgUa$ algebra}\label{sec:algUa}
\subsubsection{}
We consider the graded dual $\Nich(X)^*$ of the Nichols algebra in
Sec.~\ref{sec:Nich-asymm} and define elements $\DB,\DF\in X^*$ by
requiring that the only nonzero evaluations that they have with the
PBW basis in $\Nich(X)$ be $\eval{\DB}{\Fi}=1$ and
$\eval{\DF}{\Fii}=1$ (recall that using the notation for the PBW basis
in~\bref{sec:a-basis}, $B=\FB{1}$ and $F=\F{1}$).

The quotient by the kernel of a bilinear form is known to be another
characterization of a Nichols algebra.  Also, from~\eqref{respective},
the braiding matrix in the basis $(\DB,\DF)$ coincides
with~\eqref{qij-asymm}.  Therefore, $\Nich(X^*)$ is isomorphic to
$\Nich(X)$ and is the quotient of $T(X^*)$ by the ideal generated
by\footnote{That the elements in \eqref{dual-ideal} are in the kernel
  of the form is in fact obvious for $\DB^2$ and $\DF^p$; for the last
  element in \eqref{dual-ideal}, it is easy to verify its vanishing on
  all degree-three elements of the PBW basis in~$\Nich(X)$
  (see~\bref{sec:a-basis}): $\Fb{3} = \xx^2 \q^{-2} \Fi \Fii \Fii +
  \xx \q^{-1} \Fii \Fi \Fii + \Fii \Fii \Fi$, \ $\bX{3} =(1 - \q^2)
  \Fii \Fi \Fii + q^{-1} \xx^{-1}(1 - \q^4) \Fii \Fii \Fi$, and
  $\bFb{3} = (1 - \q^{-2}) \Fi \Fii \Fi$.}
\begin{equation}\label{dual-ideal}
   \DB^2, \ \ \DF^p, \ \ \text{and} \ \ \xx^2
   \DB\DF\DF - \xx(\q+\q^{-1})\DF\DB\DF + \DF\DF\DB.
\end{equation}

\subsubsection{}\label{dFKF-relations}
The dual algebra $\Nich(X)^*$ acts on each multivertex
$\Nich(X)$-module in accordance with~\eqref{dual-acts}.  Claiming this
requires verifying that the action by elements \eqref{dual-ideal}
commutes with the left adjoint action of $\Nich(X)$.  We show this.

For the last element in~\eqref{dual-ideal}, we take the general
``commutator'' formula~\eqref{commutator-rs} with $r=3$ (in which case
it takes the graphic form that differs
from~\eqref{commutator-rs-graph} only in the number of strands
``inside the loop'').  The two brackets in the right-hand side of
\eqref{commutator-rs} then become
$(\id+\Psi_{-1}+\Psi_{-1}\Psi_{-2}\bigr)$ and $\bigl(\id+\Psi_{-2} +
\Psi_{-2}\Psi_{-1}\bigr)$ (using the conventions set
in~\bref{Nich-dual}), and it is straightforward to verify that both
vanish when applied to $\xi^2\DB\tensor\DF\tensor\DF -
\xi(\q+\q^{-1})\DF\tensor\DB\tensor\DF +\DF\tensor\DF\tensor\DB$.

For $\DF^p$, similarly, the two elements of the braid group algebra
that occur in applying \eqref{commutator-rs} are $(\id+\Psi_{-1} +
\Psi_{-1}\Psi_{-2} + \dots + \Psi_{-1}\Psi_{-2}\dots\Psi_{-p+1}\bigr)$
\,and \,$\bigl(\id+\Psi_{-p+1} +{}$\\$
\Psi_{-p+1}\Psi_{-p+2} + \dots + \Psi_{-p+1}\dots\Psi_{-1}\bigr)$. \
Both are immediately seen to vanish when applied to $\DF^{\otimes p}$.

For $\DB^2$, totally similarly but even simpler, everything reduces to
the ``basic property of a fermion'' $(\id + \Psi)\DB\tensor\DB=0$.

\subsubsection{}\label{a-prebosonize}
The ``half-way algebra'' $\Nich(X)^*\tensor\Nich(X)\tensor k\bGamma$
in~\eqref{KF}--\eqref{KdF} is therefore an associative algebra on
generators $\Fi$, $\Fii$, $\bK_{\Fi}$, $\bK_{\Fii}$, $\DB$, $\DF$ with
the relations
\begin{alignat*}{2}
  \bK_{\Fi}\Fi \bK_{\Fi}^{-1}&=\Fi,\quad&
  \bK_{\Fi}\Fii \bK_{\Fi}^{-1}&=\xi^{-2}\q^{-2} \Fii,
  \\
  \bK_{\Fii} \Fi \bK_{\Fii}^{-1}&=\xi^2\q^{-2}\Fi,\qquad &
  \bK_{\Fii}\Fii \bK_{\Fii}^{-1}&=\q^4\Fii,
  \\
  \DB \Fi + \Fi \DB&=\one-\bK_{\Fi},
  \qquad& \DB \Fii - \xi\q^{-1}\Fii \DB&=0,
  \\  
  \DF \Fi - \xi^{-1}\q^{-1}\Fi \DF&=0,
  \quad& \DF \Fii - \q^2\Fii \DF&=\one - \bK_{\Fii},
  \\
  \bK_{\Fi}\DB \bK_{\Fi}^{-1}&=\DB,\quad&
  \bK_{\Fi}\DF \bK_{\Fi}^{-1}&=\xi^{2}\q^{2} \DF,
  \\
  \bK_{\Fii} \DB \bK_{\Fii}^{-1}&=\xi^{-2}\q^{2}\DB,\qquad &
  \bK_{\Fii}\DF \bK_{\Fii}^{-1}&=\q^{-4}\DF,\\[-1.3\baselineskip]
\end{alignat*}
\begin{gather*}
  \Fi^2 = 0, \quad \Fii^p = 0, \quad \xx^2 \Fi\Fii\Fii -
  \xx(\q+\q^{-1})\Fii\Fi\Fii + \Fii\Fii\Fi = 0,
  \\
  \DB^2 = 0, \quad \DF^p = 0, \quad \xx^2 \DB\DF\DF -
  \xx(\q+\q^{-1})\DF\DB\DF + \DF\DF\DB = 0
\end{gather*}
(and we can also consistently impose the relations $\bK_{\Fi}^p =
\one$ and $\bK_{\Fii}^p = \one$).

To obtain a Hopf algebra, as explained above, we set $\xx^2=1$, as is
required by condition~\eqref{Sommerh}, and take an Abelian group
$\Gamma$ such that the braided vector space $X$ become an object in
$\GGyd$.  We choose $\Gamma$ to be the Abelian group with two
generators $K$ and $k$, with
\begin{equation*}
  k^{2p}=\one,\qquad K^{2p}=\one,
\end{equation*}
acting and coacting as
\begin{alignat*}{2}
  k\leftii\Fi&=-\Fi,\qquad & k\leftii\Fii&=\xx\q \Fii,
  \\
  K\leftii\Fi &=\xx
  \q \Fi,\qquad &K\leftii\Fii&=\q^{-2}\Fii,
  \\
  \Fi&\mapsto k^{-1}\tensor \Fi,\quad &\Fii&\mapsto K^{-1}\tensor \Fii.
\end{alignat*}
Then the ``half-way algebra'' is a subalgebra in the Hopf algebra
$\algUa = \Nich(X^*)\tensor\Nich(X)\tensor k\Gamma$, embedded via
$\bK_{\Fi}=k^{-2}$ and $\bK_{\Fii}=K^{-2}$.  The generators $\phi_i$
introduced in~\bref{double-B} are now $k\DB$ and $K\DF$.  We change
the normalization in order to obtain more conventional commutation
relations in what follows: we set
\begin{equation*}
  \Ei = \ffrac{-1}{\q - \q^{-1}}\,k\DB
  \qquad\text{and}\qquad
  \Eii = \ffrac{1}{\q - \q^{-1}}\,K\DF.
\end{equation*}
The Hopf algebra $\algUa$ with the generators chosen this way is fully
described below for the convenience of further reference.

\subsubsection{The Hopf algebra $\balgUa$}\label{algUa-def} It
follows that the double-bosonization algebra $\algUa$ is the algebra
on generators $\Fi$, $\Fii$, $k$, $K$, $\Ei$, $\Eii$ with the
following relations.  First, $\algUa$ contains a $\Ures$ subalgebra
(which is also a Hopf subalgebra) generated by $\Eii$, $K$, and
$\Fii$, with the relations
\begin{equation}\label{Ures}
  \begin{gathered}
    K\Fii=\q^{-2}\Fii K,\quad \Eii \Fii-\Fii \Eii
    =\mfrac{K-K^{-1}}{\q-\q^{-1}},\quad K\Eii=\q^{2}\Eii K,
    \\
    \Fii^p=0,\qquad \Eii^p=0,\qquad K^{2p}=\one.
  \end{gathered}
\end{equation}
Next, $\Ures$ and $k$ generate a subalgebra, denoted by $\Uresk$ in
what follows, with further relations
\begin{gather}\label{Ustar}
  k\Eii =\xx\q^{-1}\Eii k,\quad k\Fii=\xx\q \Fii k,\quad
  k^{2p}=\one,\quad kK = Kk.
\end{gather}
The other relations in $\algUa$ are
\begin{equation}\label{the-other}
  \begin{gathered}
  k\Fi =-\Fi k,\quad K\Fi =\xx\q \Fi K, \quad k\Ei =-\Ei k,\quad
  K\Ei = \xx\q^{-1}\Ei K,
  \\[-2pt]
  \Fi^2 = 0, \quad \Fi \Ei -\Ei \Fi =\mfrac{k-k^{-1}}{\q-\q^{-1}},
  \quad \Ei^2 = 0,
  \\
  \Fii \Ei -\Ei \Fii=0, \qquad \Fi \Eii -\Eii \Fi =0,
  \\
  \Fii\Fii\Fi - \xx(\q+\q^{-1})\Fii\Fi\Fii + \Fi\Fii\Fii = 0, \quad
  \Eii\Eii\Ei - \xx(\q+\q^{-1}) \Eii\Ei\Eii + \Ei\Eii\Eii = 0.
\end{gathered}
\end{equation}
The coproduct, antipode, and counit are given by
\begin{align}
  \label{Delta-a}
  &\begin{alignedat}{2}
    \Delta(\Fii)&=\Fii\otimes1 + K^{-1}\otimes \Fii,\qquad
    &\Delta(\Eii)&=\Eii\tensor K+\one\tensor\Eii,
    \\
    \Delta(\Fi)&= \Fi\otimes1 + k^{-1}\otimes \Fi,
    &\Delta(\Ei)&=\Ei\tensor k+\one\tensor\Ei,
  \end{alignedat}
  \\
  &\begin{alignedat}{2}
    S(\Fi)&=-k\Fi,\quad S(\Fii)=-K\Fii,    
    \quad&\zS(\Ei)&=-\Ei k^{-1}, \quad \zS(\Eii)=-\Eii K^{-1},
  \end{alignedat}
  \\
  &\begin{alignedat}{2}
    \epsilon(\Fi)&=0,\quad\epsilon(\Fii)=0,
    &\quad
    \epsilon(\Ei)&=0,\quad \epsilon(\Eii)=0,
  \end{alignedat}
\end{align}
with $k$ and $K$ group-like.

\subsection{The $\balgUs$ algebra}\label{sec:algUs}
In the graded dual $\Nich(\Xs)^*$, we define $\Da,\Db\in \Xs^*$ by
requiring that the only nonzero evaluations that they have with the
PBW basis in $\Nich(\Xs)$ be $\eval{\Da}{a}=1$ and $\eval{\Db}{b}=1$.
In accordance with~\eqref{q-hat}, the braiding matrix of $(\Da,\Db)$
coincides with~\eqref{qij-symm}.  It is also easy to see that the
coproduct in~\bref{sec:sDelta} immediately implies the relations
\begin{gather*}
  \Da^2=0,\qquad \Db^2 = 0,\qquad
  (\Da \Db)^p - \xx^{-p} (\Db \Da)^p = 0.
\end{gather*}
The ``half-way algebra,'' in addition to the relations in $\Nich(\Xs)$
and $\Nich(\Xs^*)$, has the relations
\begin{alignat*}{2}
  \Da\, a + a\, \Da &= 1 - \bK_a,
  \qquad &\Da\, b + \xx\q\, b\, \Da &= 0,
  \\
  \Db\, a + \xx^{-1}\q\, a\, \Db &= 0,
  \quad &\Db\, b + b\, \Db &= 1 - \bK_b,
\end{alignat*}
where
\begin{alignat*}{2}
  \bK_a a\bK_a^{-1} &= a,\quad & \bK_a b\bK_a^{-1} &= \q^2 b,
  \\
  \bK_b a\bK_b^{-1} &= \q^2 a,\quad &\bK_b b\bK_b^{-1} &= b
\end{alignat*}
(and we can set $\bK_a^p=1$ and $\bK_b^p=1$).

\subsubsection{The Hopf structure of $\balgUs$}\label{algUs-def}
To obtain a Hopf algebra, as in~\bref{double-B}, we assume that
$\xi^2=1$ in accordance with condition~\eqref{Sommerh} and take
$\Gamma$ to be the Abelian group with generators $K_1$ and $K_2$,
$K_i^{2p}=\one$, acting and coacting on the basis $F_1=a$ and $F_2=b$
in $\Xs$~as
\begin{alignat*}{2}
  K_1\leftii F_1&=-F_1,\quad
  &K_1\leftii F_2&=-\xx\q^{-1} F_2,
  \\
  K_2\leftii F_1&=-\xx
  \q^{-1} F_1,\quad
  &K_2\leftii F_2&=-F_2,
\end{alignat*}
and $F_i\mapsto K_i^{-1}\tensor F_i$.  In accordance
with~\bref{double-B}, the double bosonization
$\algUs=\Nich(\Xs^*)\tensor\Nich(\Xs)\tensor k\Gamma$ is then the algebra
on generators $F_1$, $F_2$, $\phi_1$, $\phi_2$, $K_1$, $K_2$ with the
relations
\begin{alignat*}{2}
  K_1^{2p}&=\one,\quad K_2^{2p}=\one,
  &K_1 K_2 &= K_2 K_1,
  \\
  F_1^2&=0,\quad F_2^2=0,
  &(F_1 F_2)^p &- \xx^{-p} (F_2 F_1)^p = 0,
  \\
  K_1F_1K_1^{-1}&=-F_1,\quad
  &K_1F_2K_1^{-1}&=-\xx\q^{-1} F_2,
  \\
  K_2F_1K_2^{-1}&=-\xx
  \q^{-1} F_1,\quad
  &K_2F_2K_2^{-1}&=- F_2,
  \\
  F_1\phi_1-\phi_1F_1&=K_1-K_1^{-1},
  &\quad F_1\phi_2-\phi_2F_1&=0,
  \\
  F_2\phi_1-\phi_1F_2&=0,
  &\qquad 
  F_2\phi_2-\phi_2F_2&=K_2-K_2^{-1},
  \\
  K_1\phi_1 K_1^{-1}&=-\phi_1,&\quad
  K_1\phi_2 K_1^{-1}&=-\xx\q\phi_2,
  \\
  K_2\phi_1 K_2^{-1}&=-\xx\q\phi_1,&\quad K_2\phi_2 K_2^{-1}&=-\phi_2,
  \\
  \phi_1^2&=0,\quad \phi_2^2= 0,\quad
  &(\phi_1 \phi_2)^p &- \xx^{-p} (\phi_2 \phi_1)^p = 0
\end{alignat*}
and with the Hopf-algebra structure defined by
\begin{alignat}{2}
  \label{Delta-s}
  \Delta(F_i)&=F_i\otimes 1 + K_i^{-1}\otimes F_i,\quad
  &\Delta(\phi_i)&=\phi_i\tensor K_i + \one\tensor\phi_i,
  \\
  S(F_i)&=-K_iF_i,\quad
  &S(\phi_i)&=-\phi_i K_i^{-1},
  \\
  \epsilon(F_i)&=0,\quad
  &\epsilon(\phi_i)&=0,
\end{alignat}
$i=1,2$ (and $K_1$ and $K_2$ group-like).  Unlike in the
''asymmetric'' case, we do not renormalize the generators to produce
$\q - \q^{-1}$ in denominators.

\subsection{Isomorphism}\label{sec:iso}
From~\bref{algUa-def} and \bref{algUs-def}, we have double
bosonizations $\algUa$ and $\algUs$ of the ``asymmetric'' and
``symmetric'' Nichols algebras.  They turn out to be ``essentially the
same''---related somewhat simpler than their CFT{} counterparts
in~\bref{s:Wak}.
\begin{it}%
\begin{enumerate}
\item\label{iso-assoc}
    The algebras $\algUa$ and $\algUs$ are isomorphic as associative
    algebras.  Explicitly, the isomorphism $\ssigma:\algUa\to\algUs$
    is given by
  \begin{alignat*}{2}
    \ssigma(F_1)&= (\q - \q^{-1}) (\xi B F - \q^{-1} F B),
    &\quad
    \ssigma(F_2)&=-(\q - \q^{-1}) C k^{-1},
    \\
    \ssigma(\phi_1)&= E C - \xi \q C E,
    &\quad
    \ssigma(\phi_2) &= B k
    \\
    \ssigma(K_1) &= K k,
    &\quad
    \ssigma(K_2) &= k^{-1}
  \end{alignat*}
  and the inverse map is
  \begin{alignat*}{2}
    \ssigma^{-1}(\Fi)&=-\fffrac{1}{(\q - \q^{-1})^2} (F_1 F_2 + \xi \q
    F_2 F_1), &\quad \ssigma^{-1}(\Fii)&=\phi_2 K_2,
    \\
    \ssigma^{-1}(\Ei)&=-\fffrac{1}{\q - \q^{-1}} F_2 K_2^{-1}, &\quad
    \ssigma^{-1}(\Eii)&= -\q^{-1} (\phi_1 \phi_2 + \xi \q \phi_2
    \phi_1),
    \\
    \ssigma^{-1}(K)&= K_1K_2, &\quad \ssigma^{-1}(k)&= K_2^{-1}.
  \end{alignat*}
  
\item\label{simi}The two coalgebra strcutures, $\Delta_{\asymm}$
  defined in~\eqref{Delta-a} and $\Delta_{\symm}$ defined
  in~\eqref{Delta-s}, are related as stated in~\eqref{PHI} with
  \begin{align*}
    \Phi &= \one\tensor\one + (\q - \q^{-1}) \Fi  k\tensor \Ei k^{-1}.
  \end{align*}

\item The antipodes are related as $US_{\asymm}(\ssigma(x))U^{-1} =
  \ssigma S_{\symm}(x)$ with $U=\one-(\q - \q^{-1}) \Fi\Ei k$.
\end{enumerate}
\end{it}

In proving that $\ssigma$ is an algebra morphism, we must of course
verify that relations are mapped into relations.  It is immediate to
see that $\ssigma(F_1)\ssigma(F_2)$ vanishes due to the
``$\Fii\Fii\Fi$'' relation in~\eqref{the-other} (and its
consequence~\eqref{FBFB}).  To see how $(F_1 F_2)^p - \xi^{-p} (F_2
F_1)^p$ is mapped by $\ssigma$, we first inductively establish the
identities
\begin{align*}
  \ssigma:\fffrac{1}{(\q - \q^{-1})^{2 n}}
  \bigl((F_1 F_2)^n - \xi^n (F_2 F_1)^n\bigr) &\mapsto
  (-1)^n  F^{n}
  + (-1)^n \xi^n \q^{-2 n + 1} (\q^{n} + 1)
   C B F^{n} k^{-1}
  \\
  &\quad{}+ (-1)^{n + 1} \xi^{n + 1} \q^{-2 n} (\q^{n} + 1)
   C F B F^{n - 1} k^{-1},
\end{align*}
whence it indeed follows that the right-hand side vanishes at $n=p$
due to the relations in $\algUa$ and the fact that $\q^p+1=0$. \ For
the $\phi_i$, everything is totally similar.

That the above $\Phi$ does the job in~\eqref{PHI} is verified on the
$\algUs$ generators straightforwardly.  The map relating the antipodes
is then standard, $U=\Phi^{(1)}S(\Phi^{(2)})$~\cite{[MaOe]}.

\subsubsection{}
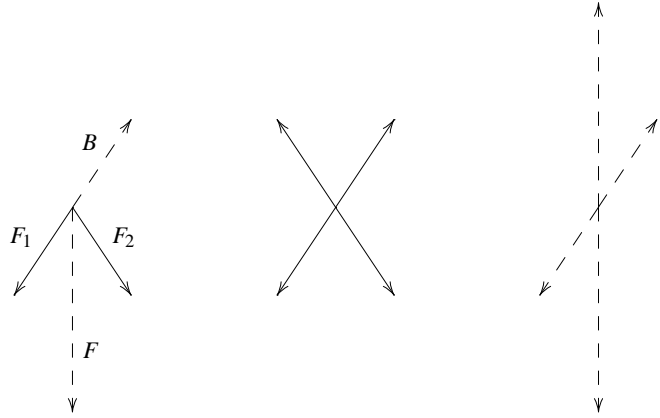
\begin{figure}[tb]
  \centering
  \begin{equation*}
    \xymatrix@R=3pt@C=3pt{
      &&&&&&&&&&&&&&&&&&&&\nothing="phi3"&&&
      \\
      &&&&&&&&&&&&&&&&&&&&&&&
      \\
      &&&&&&&&&&&&&&&&&&&&&&&
      \\
      &&&&&&&&&&&&&&&&&&&&&&&
      \\
      &&&&\nothing="phi1"&&&&&\nothing="phi22"
      &&&&\nothing="phi12"&&&&&&
      &&&\nothing="phi13"&
      \\
      &&&&&&&&&&&&&&&&&&&&&&&
      \\
      &&&&&&&&&&&&&&&&&&&&&&&
      \\
      &&\nothing="top"&&&&&&&&&
      \nothing="top2"&&&&&&&&&\nothing="top3"&&&
      \\
      &&&&&&&&&&&&&&&&&&&&&&&
      \\
      &&&&&&&&&&&&&&&&&&&&&&&
      \\
      \nothing="F1"&&&&\nothing="F2"&&&&&
      \nothing="F12"&&&&\nothing="F22"&&&&&\nothing="F13"&&&&&
      \\
      &&&&&&&&&&&&&&&&&&&&&&&
      \\
      &&&&&&&&&&&&&&&&&&&&&&&
      \\
      &&&&&&&&&&&&&&&&&&&&&&&
      \\
      &&\nothing="F"&&&&&&&&&&&&&&&&&&\nothing="F3"&&&
      \\
      \ar_{F_1}"top"+<0pt,0pt>;"F1"+<0pt,0pt>
      \ar^{F_2}"top"+<0pt,0pt>;"F2"+<0pt,0pt>
      \ar^(.6){\Fi}@{{}{--}{>}}"top"+<0pt,0pt>;"phi1"+<0pt,0pt>
      \ar^(.7){\Fii}@{{}{--}{>}}"top"+<0pt,0pt>;"F"+<0pt,0pt>
      \ar"top2"+<0pt,0pt>;"F12"+<0pt,0pt>
      \ar"top2"+<0pt,0pt>;"F22"+<0pt,0pt>
      \ar"top2"+<0pt,0pt>;"phi22"+<0pt,0pt>
      \ar"top2"+<0pt,0pt>;"phi12"+<0pt,0pt>
      \ar@{{}{--}{>}}"top3"+<0pt,0pt>;"phi13"+<0pt,0pt>
      \ar@{{}{--}{>}}"top3"+<0pt,0pt>;"F13"+<0pt,0pt>
      \ar@{{}{--}{>}}"top3"+<0pt,0pt>;"F3"+<0pt,0pt>
      \ar@{{}{--}{>}}"top3"+<0pt,0pt>;"phi3"+<0pt,0pt>
    }
  \end{equation*}
  \caption{\small\textsc{Left:} generators of the ``symmetric'' (solid
    lines) and ``asymmetric'' (dashed lines) Nichols algebras. \
    \textsc{Middle and Right:} two systems of generators in the common
    double bosonization $\algU$ of the two Nichols algebras.%
  }
  \label{fig:screenings}
\end{figure}

The isomorphism of associative algebras $\algUa\simeq\algUs$ is not
unexpected if we recall that the Nichols algebras $\Nich(\Xa)$ and
$\Nich(\Xs)$ are related by a Weyl pseudoreflection, which (at the
level of bosonizations $\Nich(X)\Smash k\Gamma$, to be precise)
roughly amounts to the following procedure~\cite{[Ag-0804-standard]}:
pick up a $\Nich(X)$ generator $F_{\ell}$; drop
$F_{\ell}$ (the corresponding $1$-dimensional subspace) and ``add''
the dual $\dF_\ell$ instead; and replace each $F_i$, $i\neq\ell$, with
$(\mathrm{ad}_{F_{\ell}})^{\text{max}} F_i$, where ``max'' means
taking the maximum power (of the braided adjoint) that does not yet
lead to identical vanishing.  For $\Nich(\Xs)$, whose fermionic
generators $F_1$ and $F_2$ are shown on the left in
Fig.~\ref{fig:screenings},
this means dropping $F_1$, introducing $B=\dF_1$ instead, and
replacing $F_2$ with $\Fii=[F_1,F_2]_{\q}$ (a braided commutator);
these $B$ and $F$ are generators of~$\Nich(\Xa)$.  In the double
bosonization of each Nichols algebra, the generators are the original
$\Nich(X)$ generators and their opposite ones, and the two Nichols
algebras related by a Weyl reflection yield two systems of generators
in \textit{the same}~$\algU$.

\subsubsection{}\label{rhorho}
We note that Fig.~\ref{fig:screenings} is reproduced in the structure
of Eqs.~\eqref{rho-s} and \eqref{rho-a}.



\subsection{Simple $\balgU$-modules}
In view of~\bref{sec:iso}, the representation theories of $\algUa$ and
$\algUs$ are equivalent.  We choose $\algUa$ with $\xi=1$ and let it
be denoted simply by $\algU$ (the algebra with $\xi=-1$ has an
equivalent representation category~\cite{[BPW]}).

We quote some of the results established in~\cite{[nich-sl2-2]}.  The
algebra $\algU$ has $4p^2$ simple modules, which are labeled as
\begin{gather*}
  \repZ^{\alpha,\beta}_{s,r},\qquad \alpha,\beta=\pm,\quad
  s=1,\dots,p,\quad r=0,\dots,p-1,
\end{gather*}
and have the dimensions
\begin{equation}\label{dim-simp-mod}
  \dim\repZ^{\alpha,\beta}_{s,r}=
  \begin{cases}
    2s-1,&r=0,\quad 1\leq s\leq p,\\
    2s+1,&r=s,\quad 1\leq s\leq p-1,\\
    4s,&r\neq0,s,\quad 1\leq s\leq p-1,\\
    4p,&1\leq r\leq p-1,\quad s=p.
  \end{cases}
\end{equation}
On the highest-weight vector of $\repZ^{\alpha,\beta}_{s,r}$, \ $k$
and $K$ have the respective eigenvalues
\begin{equation}\label{eigens}
  \beta\q^{-r}\qquad\text{and}\qquad \alpha\q^{s-1}.
\end{equation}

\section{\textbf{$\balgU$ and the $\bWWtridwa$ algebra}\label{sec:7}}
We consider the triplet--triplet $W$-algebra $\WWtritri$ (which
corresponds to $j=0$ in the ``asymmetric'' case, for definiteness).
Introducing an Abelian group $\Gamma$ such that $X\in\GGyd$ implies an
effect that has no clear analogue in the known $\Wpone/$Virasoro case:
not all of the $W$-algebra commutes with~$\Gamma$.

\subsection{$\myboldsymbol{\Gamma}$ in terms of free fields}
In the asymmetric free-field realization, the generators of $\Gamma$
can be represented in terms of zero modes of the fields as
\begin{equation}\label{Kk}
  k=e^{-i\pi(\Aone)_0}=e^{\frac{1}{2}i\pi(\Aii_0+\Afour_0)},\qquad
  K=e^{-i\pi(\Atwo)_0} = e^{-i\pi\Aii_0}
\end{equation}
(see the field redefinition in~\bref{ortho}).  It follows that $k$
\textit{anti}commutes with $\Jplus(z)$ and $\Jminus(z)$.

\subsection{$\bWWtridwa\myboldsymbol{{}\subset{}}\bWWtritri$}
To maintain the idea that the algebras ``on the Hopf side'' and ``on
the CFT{} side'' centralize each other, we have to choose a subalgebra
in~$\WWtritri$ that commutes with $\Gamma$.  We recall
from~\bref{a:j=0} that the (not necessarily minimal) set of fields
generating the $\WWtritri$ algebra was
\begin{gather*}
  \Jminus_0\Jminus_0\W^a(z),\quad
  \Jminus_0\W^a(z),\quad \W^a(z),\qquad a={+},0,{-},
\end{gather*}
which gave rise to $\Jminus(z)$, $\Jnaught(z)$, $\Jplus(z)$ in their
OPEs.  The subalgebra that centralizes $\Gamma$ is generated by the
fields
\begin{gather*}
  \Jminus_0\Jminus_0\W^a(z),\quad \W^a(z),\qquad a={+},0,{-},
  \\
  \intertext{and} (\Jminus(z))^2,\quad\Jnaught(z),\quad(\Jplus(z))^2,
  \quad T_{\text{Sug}}(z)
\end{gather*}
($(J^{\pm}(z))^2$ do occur in the OPEs of $\Jminus_0\Jminus_0\W^a(z)$
and $\W^b(z)$, but we do not discuss the minimal set of generators
now, emphasizing the ``$J$-squaredeness'' instead).  We let this
algebra be denoted by $\WWtridwa$, because with the middle terms
dropped, the triplets under the horizontal $s\ell(2)$ become
``doublets'' with respect to zero modes of $(\Jminus(z))^2$ and
$(\Jplus(z))^2$.

The representation theories of $\WWtridwa$ and $\WWtritri$ are not
\textit{very} different.  It remains to be seen whether$/$how
reintroducing $\Jplus(z)$ and $\Jminus(z)$ as such, not their squares,
spoils some presumably nice properties of the $\WWtridwa$ theory.

\subsection{The $\bWWtridwa$ action on vertices}\label{J2}
We recall the vertex operators
$\Vsr_{s,r;\vartheta}^{\nu,\mu}[n,m](z)$ introduced in~\bref{V[n,m]}.
We saw there that acting with $\W^{+}_{\ell}$ and $\W^{-}_{\ell}$ maps
over the values of~$n$.  We next observe that the modes of
$(\Jplus(z))^2$ and $(\Jminus(z))^2$ map over the values of~$m$.

The annihilation conditions with respect to the modes of
$(\Jplus(z))^2$ and $(\Jminus(z))^2$ are
\begin{align*}
  \bigl((\Jplus)^2\bigr)_{2\vartheta-1+\ell}\Vsr_{s,r;\vartheta}^{\nu,\mu}[n,m](z) &=0,\quad \ell\geq 0,
  \\
  \bigl((\Jminus)^2\bigr)_{3-2\vartheta+\ell}\Vsr_{s,r;\vartheta}^{\nu,\mu}[n,m](z) &= 0,\quad \ell\geq 0,
\end{align*}
and the maximum modes that are generically nonvanishing act as
\begin{align*}
  \bigl((\Jplus)^2\bigr)_{2\vartheta-2}\Vsr_{s,r;\vartheta}^{\nu,\mu}[n,m](z) &=
  \Vsr_{s,r;\vartheta}^{\nu,\mu}[n,m+1](z),
  \\
  \bigl((\Jminus)^2\bigr)_{2-2\vartheta}\Vsr_{s,r;\vartheta}^{\nu,\mu}[n,m](z)
  &=
  \bigl(\fffrac{r - s}{p} + \mu + \nu + 2 n + 2 m\bigr)
  \bigl(\fffrac{r - s}{p} + \mu + \nu + 1 + 2 n + 2 m\bigr)  
  \\
  &\quad\times{}
  \bigl(\fffrac{r}{p} + \mu + 2 m\bigr)
  \bigl(\fffrac{r}{p} + \mu - 1 + 2 m\bigr)
  \Vsr_{s,r;\vartheta}^{\nu,\mu}[n,m-1](z).
\end{align*}
The vanishing conditions for the brackets in the last formula, as well
as in~\eqref{Wminus-act}, allow making some simple observations about
the occurrence of $\bWWtridwa$ submodules.

\subsection{$\bWWtridwa$ submodules and the correspondence with
  $\balgU$ representations}
For the algebra $\WWW\equiv\WWtridwa$, we
describe a set of its (conjecturally simple) modules
$\modZ_{s,r;\vartheta}^{\alpha,\beta}$ and set them in correspondence
with the simple $\algU$-modules $\repZ_{s,r}^{\alpha,\beta}$.  We set
$\alpha=(-1)^\nu$ and $\beta=(-1)^\mu$ (with $\nu$ and $\mu$ taking
values $0$ or $1$) hereafter in this subsection.

Let $\modF_{s,r;\vartheta}^{\alpha,\beta}$ be the space spanned by
$P(\partial\Aii,\partial\Afour,\partial\Athree)
\Vsr_{{s},{r};\vartheta}^{\nu,\mu}[n,m](z)$, where $n,m\in\oZ$ and $P$
are differential polynomials.  It bears a natural $\WWW$-action 
and, as a $\WWW$-module, plays the role of a Verma module, with the
$\Vsr_{{s},{r};\vartheta}^{\nu,\mu}[n,m](z)$, $n,m\in\oZ$, being the
extremal vectors.  The $\WWW$-module structure of
$\modF_{s,r;\vartheta}^{\alpha,\beta}$ depends on the parameters
$1\leq s\leq p$ and $0\leq r\leq p-1$, and we now list some
characteristic irreducibility/reducibility cases.

\begin{enumerate}
\item\label{item:1} For $s=p$ and $r\neq0$, we conclude
  from~\bref{V[n,m]} and~\bref{J2} that each extremal vector is
  reachable by the $\WWW$ action from any other extremal vector; we
  conjecture that
  $\modZ_{p,r;\vartheta}^{\alpha,\beta}=\modF_{p,r;\vartheta}^{\alpha,\beta}$
  is irreducible in this case.

\item\label{item:2} For $1\leq s\leq p-1$ and $r\neq0$, $r\neq s$, it
  follows from~\bref{V[n,m]} that extremal vectors with $n<0$ are
  unreachable from extremal vectors with $n\geq0$.  There is therefore
  a proper submodule
  $\modZ_{s,r;\vartheta}^{\alpha,\beta}\subset\modF_{s,r;\vartheta}^{\alpha,\beta}$
  generated from $\Vsr_{{s},{r};\vartheta}^{\nu,\mu}[0,0](z)$.

\item\label{item:3} For $1\leq s\leq p-1$ and $r=s$, it follows
  from~\bref{V[n,m]} that extremal vectors with $n<0$ are unreachable
  from the vectors with $n\geq0$; from~\bref{J2}, moreover, we see
  that in the case $(\nu=0,\,\mu=0)$, extremal vectors with $n<-m$ are
  unreachable from the vectors with $n\geq-m$, and in the three
  remaining cases $(\nu=0,\,\mu=1)$, $(\nu=1,\,\mu=0)$, and
  $(\nu=1,\,\mu=1)$,  extremal vectors with $n<-m-1$ are
  unreachable from those with $n\geq-m-1$.  The smallest submodule
  $\modZ_{s,s;\vartheta}^{\alpha,\beta}\subset\modF_{s,s;\vartheta}^{\alpha,\beta}$
  has the extremal vectors as shown in Fig.~\ref{fig:submodule},
  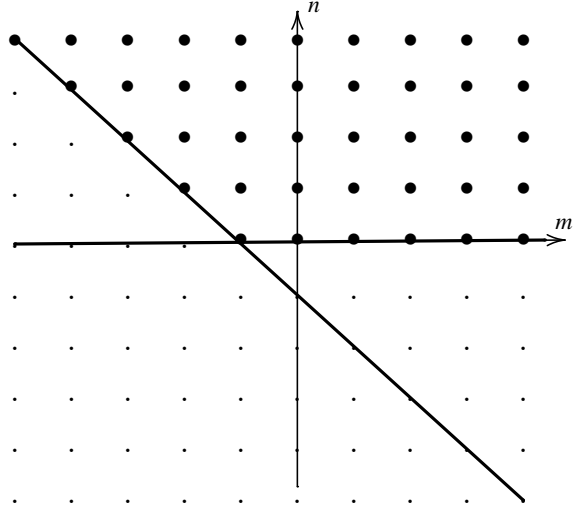
\begin{figure}[tb]
    \centering
    \begin{equation*}
      \xymatrix@C10pt@R10pt{
        &&&&&&&&&
        \\
        *{\bullet}\ar@{-}@*{[|(3)]}[]+<0pt,1pt>;[9,9]+<0pt,-2pt>&*{\bullet}&*{\bullet}&*{\bullet}&*{\bullet}&*{\bullet}&*{\bullet}&*{\bullet}&*{\bullet}&*{\bullet}
        \\
        .&*{\bullet}&*{\bullet}&*{\bullet}&*{\bullet}&*{\bullet}&*{\bullet}&*{\bullet}&*{\bullet}&*{\bullet}
        \\
        .&.&*{\bullet}&*{\bullet}&*{\bullet}&*{\bullet}&*{\bullet}&*{\bullet}&*{\bullet}&*{\bullet}
        \\
        .&.&.&*{\bullet}&*{\bullet}&*{\bullet}&*{\bullet}&*{\bullet}&*{\bullet}&*{\bullet}
        \\
        .\ar@{-}@*{[|(3)]}[]+<0pt,-1.5pt>;[0,9]+<8pt,0pt>\ar@*{[|(1.5)]}^(.98){m}[]+<0pt,-1.5pt>;[0,10]&.&.&.&*{\bullet}&*{\bullet}&*{\bullet}&*{\bullet}&*{\bullet}&*{\bullet}&
        \\
        .&.&.&.&.&.&.&.&.&.
        \\
        .&.&.&.&.&.&.&.&.&.
        \\
        .&.&.&.&.&.&.&.&.&.
        \\
        .&.&.&.&.&.&.&.&.&.
        \\
        .&.&.&.&.&.\ar@*{[|(1.5)]}_(.98){n}[-10,0]&.&.&.&.
      }
    \end{equation*}
    \caption[States in the submodule]{\small The filled dots are the
      extremal vectors
      $\protect\Vsr_{{s},{s};\vartheta}^{\nu,\mu}[n,m](z)$ that are in
      the submodule $\modZ_{s,s;\vartheta}^{\alpha,\beta}$, in any of
      the three cases $(\nu=0,\,\mu=1)$, $(\nu=1,\,\mu=0)$, and
      $(\nu=1,\,\mu=1)$. \ The submodule is the intersection of two
      submodules, respectively represented by dots to the right of
      (and including) the tilted line $n+m=-1$ and by dots above (and
      including) the horizontal line $n=0$. \ (In the case
      $(\nu=0,\,\mu=0)$, the tilted line is $n+m=0$.)}
    \label{fig:submodule}
  \end{figure}
  and can be generated from~$\Vsr_{{s},{s};\vartheta}^{\nu,\mu}[0,0](z)$.
  
\item\label{item:4} For $1\leq s\leq p-1$ and $r=0$, it follows
  from~\bref{V[n,m]} that extremal vectors with $n<0$ are unreachable
  from the vectors with $n\geq0$, and from~\bref{J2}, that extremal
  vectors with $m<0$ are unreachable from the vectors with $m\geq 0$.
  (The picture is similar to that in Fig.~\ref{fig:submodule}, but
  with the tilted line becoming the vertical line $m=0$.) \ In this
  case, a proper submodule
  $\modZ_{s,0;\vartheta}^{\alpha,\beta}\subset\modF_{s,0;\vartheta}^{\alpha,\beta}$
  also can be generated from
  $\Vsr_{{s},{0};\vartheta}^{\nu,\mu}[0,0](z)$.

\end{enumerate}

In cases \ref{item:1}, \ref{item:2}, \ref{item:3}, and \ref{item:4},
the respective submodules $\modZ_{p,r;\vartheta}^{\alpha,\beta}$,
$\modZ_{s,r;\vartheta}^{\alpha,\beta}$,
$\modZ_{s,s;\vartheta}^{\alpha,\beta}$, and
$\modZ_{s,0;\vartheta}^{\alpha,\beta}$ are conjecturally simple
$\WWW$-modules.  We also conjecture that the
$\modZ_{s,r;\vartheta}^{\alpha,\beta}$ with $1\leq s\leq p$, $0\leq
r\leq p-1$, $\vartheta\in\oZ$, and $\alpha,\beta=\pm$ are all simple
modules of~$\WWW$ (we do not define
$\modZ_{p,0;\vartheta}^{\alpha,\beta}$ here; the case $s=p$, $r=0$
will be considered elsewhere).

The correspondence with the simple $\algU$-modules is suggested by the
$\Gamma$-action on the $\WWW$-modules.  The $\Gamma$ generators
represented as in~\eqref{Kk} act on the extremal states as
\begin{align*}
  k\Vsr_{s,r;\vartheta}[n,m](z)
  &=(-1)^\mu\q^{-r}\Vsr_{s,r;\vartheta}^{\nu,\mu}[n,m](z),
  \\
  K\Vsr_{s,r;\vartheta}[n,m](z)
  &=(-1)^\nu\q^{s-1}\Vsr_{s,r;\vartheta}^{\nu,\mu}[n,m](z),
\end{align*}
which is to be compared with~\eqref{eigens}.  Because $\WWW$ and
$\algU$ commute, every vector in
$\modZ_{s,r;\vartheta}^{\alpha,\beta}$ is a highest-weight vector of
the $\algU$-module $\repZ_{s,r}^{\alpha,\beta}$.  This suggests the
correspondence
$\modZ_{s,r;\vartheta}^{\alpha,\beta}\to\repZ_{s,r}^{\alpha,\beta}$
(for $\vartheta\in\oZ$) between $\WWW$- and $\algU$-modules (quite
similar to~\cite{[FGST2],[NT]}), thus paving the way toward the Morita
equivalence of categories.  Moreover, it can be expected that
$\WWW$-fusion of the $\modZ_{s,r;\vartheta}^{\alpha,\beta}$ is closely
related to tensor products of the $\repZ_{s,r}^{\alpha,\beta}$
modules.

\section{\textbf{Conclusions and an outlook}}
The project pursued in this paper is by no means completed at this
point.  We described two Nichols algebras $\Nich(\Xa)$ and
$\Nich(\Xs)$ quite explicitly, and calculated their action on
one-vertex Yetter--Drinfeld modules, but we have not yet considered
their Yetter--Drinfeld categories in greater detail (including the
multivertex modules).  We identified generators of extended chiral
algebras and found candidates for irreducible representations in the
``best'' case ($j=0$ in the asymmetric realization), but we have not
yet computed the characters of these representations.  We outlined the
Hamiltonian reduction of the extended algebras to the previously known
one-boson ``logarithmic'' chiral algebras, but we have not explicitly
described how the fields of these last algebras are combined and
``dressed'' by the additional scalars into the extended algebras
discussed in this paper.  We derived an ordinary Hopf algebra $\algU$
from a $\Nich(X)$ by double bosonization (and noticed that $\algU$ is
``the same'' for $\Nich(\Xa)$ and $\Nich(\Xs)$), but we have not yet
articulated the functorial properties of the correspondence between
$\Nich(X)$ and $\algU$, for $X\in\HHyd$.  We noticed an encouraging
correspondence between simple $\algU$ and $\WWtridwa$ modules, but we
have not yet stated the expected categorial equivalence, with the
spectral flow taken into account.

Further prospects are therefore numerous and can be rather exciting.
\begin{enumerate}\addtolength{\itemsep}{2pt}
\item The characters of representations of the extended algebras
  constructed here can be calculated relatively straightforwardly; in
  addition to their ``nominal'' significance, they may encode some
  nontrivial combinatorial properties, being related to generation
  functions of planar partitions subjected to some set of
  constraints~\cite{[FJMM]}.


\item The associative algebra isomorphism $\ssigma$ and the twist
  $\Phi$ in~\bref{sec:iso} can be regarded as algebraic counterparts
  of the ``nonlocal change of variables''
  $\mapa{\cdot}\ccirc\mapsinv{\cdot}$ in~\bref{s:Wak}.  The similarity
  may not necessarily be superficial.

\item With the $(p,p')$ logarithmic minimal models~\cite{[FGST3]}
  deducible by Hamiltonian reduction from a theory with manifest
  $\hSL2$ symmetry, it is interesting how much this last can help in
  elucidating a number of subtle properties of the $(p,p')$
  models~\cite{[GRW],[RGW],[W-0907]}.  In particular, setting $p=2$
  and $j=1$ in the ``asymmetric'' case gives an $\hSL2_{-\half}$ model
  (with central charge $-1$), with the underlying $(2,3)$ minimal
  model (also see~\cite{[VGJS]}).
  
\item The case of \textit{three} fermionic screenings, for which some
  Nichols-algebra details were already worked out
  in~\cite{[Ag-1008-presentation]} (where the Nichols algebra
  generators were not called screenings, however), is certainly
  interesting from the CFT{} standpoint~\cite{[FS-D]}.  In higher
  rank, moreover, Nichols algebras may depend on several ``$p$-type''
  parameters, which is another interesting possibility of going beyond
  logarithmic CFTs based on rescalings of classic root lattices.

\item Going beyond finite-dimensional Nichols algebras---extending
  them by divided powers of nonfermionic generators---is certainly of
  interest in logarithmic CFTs (cf.~\cite{[BFGT]}), and possibly also
  in the theory of Nichols algebras.  A double bosonization of our
  $\Nicha$ and $\Nichs$ ``with divided powers'' can be regarded as a
  Lusztig-type extension of $\algU$.

\item Nichols algebras (and their Yetter--Drinfeld modules) have
  already appeared in~\cite{[BR]}, in the context of integrable
  deformations of conformal field theories.  Links with the results in
  this paper and possible generalizations are quite intriguing.

\item Lattice models related to logarithmic CFTs
  (see~\cite{[GV],[GRS],[GST]} and the references therein) are another
  direction where the ``braided'' standpoint may be welcome.  A spin
  chain based on our $\algU$ can also be useful in describing the spin
  quantum Hall effect~\cite{[EFSa-SaSch]}.

\item An intriguing development is to investigate modular properties
  of the extended-algebra characters and compare them with the modular
  group representation realized on the center of the Hopf algebra
  $\algU$.  Based on our previous experience~\cite{[FGST],[FGST3]},
  the exact coincidence of modular group representations can be
  expected here (which would nevertheless be a nontrivial result).  An
  even more ambitious program is to compare the modular group
  representation generated from $W$-algebra characters with the one
  realized in braided terms, much in the spirit of the recent
  development in~\cite{[FSS],[FFSS]}.  A higher-genus mapping class
  group is closely related to $n$th cohomology spaces of $\algU$ with
  coefficients in $\algU^{\tensor \ell}\tensor Y_1\tensor\dots\tensor
  Y_M$, where the $Y_i$ (vertex-operator algebra representations and
  simultaneously Nichols algebra representations) label boundary
  conditions imposed at the holes of the Riemann surface.  It would be
  very interesting to reformulate this in the framework of the
  corresponding Nichols algebra (cf.\ the complex constructed
  in~\cite{[STbr]}) by taking its multiple tensor products with
  itself~and~with~the~dual.
\end{enumerate}

\subsubsection*{Acknowledgments}
We are grateful to D.~B\"ucher, B.~Feigin, A.~Gainutdinov, I.~Runkel,
H.~Saleur, C.~Schweigert, J.~Shiraishi, and A.~Tsuchiya for valuable
discussions. \ Part of the motivation to start this project came from
discussions with C.~Schweigert. \ AMS thanks him and the Department of
Mathematics, Hamburg University, for the hospitality.  IYuT is
grateful to H.~Saleur for kind hospitality at IPhT.  This paper was
supported in part by the RFBR grant 10-01-00408.

\appendix
\section{\textbf{Some constructions and conventions for Nichols
    algebras}}
\subsection{$\bNich\myboldsymbol{(X)}$ action and coaction on
  multivertex Yetter--Drinfeld modules}\label{app:Nich}
We use the standard graphic notation for basic Hopf-algebra structures
and braiding.  The diagrams are read from top down; the convention for
braiding is $\begin{tangles}{l} \hx
\end{tangles}$ \ (and $\begin{tangles}{l} \hxx
\end{tangles}$\ is the inverse braiding).  All our conventions are
fully explained in~\cite{[STbr]}.  A warning (also articulated
in~\cite{[STbr]}) is that diagrams of two types are actually in use:
those where a line denotes a Hopf algebra (such as $\Nich(X)$) or its
module, etc., and those where a line is a copy of a braided space
(such as~$X$).

We let $X=(X,\Psi)$ denote a braided vector space and $\Nich(X)$ its
Nichols algebra.  The reader can always regard $X$ as an object in a
braided category such that the braiding induced on $X$ coincides
with~$\Psi$.

For any braided vector spaces $Y_j$ in the same category, there is a
Yetter--Drinfeld $\Nich(X)$-module structure on $\Nich(X)\tensor
Y_1\tensor\Nich(X)\tensor Y_2\tensor\dots\tensor\Nich(X)\tensor Y_N$,
given by the left adjoint action and by the coaction via
deconcatenation up to the first $Y$ space (the ``$N$-vertex''
Yetter--Drinfeld module, considered in more detail in~\cite{[STbr]}).
For one-vertex modules, in particular, the action and coaction are
\begin{equation}\label{1-vertex}
  \begin{tangles}{l}
    \vstr{200}
    \lu\object{\raisebox{18pt}{\kern-4pt\tiny$\blacktriangleright$}}
    \step[1]\id\object{\kern-10pt\rule[18pt]{12pt}{4pt}}\step[.1]\id
  \end{tangles}\ \
  =\ \ 
  \begin{tangles}{l}
    \hcd\step[.9]\id\step[.1]\id\\
    \vstr{50}\id\step[1]\hx\\
    \id\step[.9]\id\step[.1]\id\step[1]\O{S}\\
    \vstr{50}\id\step[1]\hx\\
    \hcu\step[.9]\id\step[.1]\id
  \end{tangles}
  \quad\qquad\text{and}\qquad
  \begin{tangles}{l}
    \vstr{200}
    \ld
    \step[1]\id\object{\kern-10pt\rule[18pt]{12pt}{4pt}}\step[.1]\id
  \end{tangles}\ \
  =\ \ 
  \begin{tangles}{l}
    \hcd\step[.9]\id\step[.1]\id\\
    \id\step[1]\id\step[.9]\id\step[.1]\id
  \end{tangles}
\end{equation}
where the black horizontal strip indicates that the $\Nich(X)$
(co)action applies to the tensor product as a whole---the tensor
product of a copy of $\Nich(X)$ (single vertical line) and
$V^{\{\alpha, \beta\}}$ (double vertical
line).

\subsection{A (very) brief reminder on $\bHHyd$ and Radford's
  biproduct formula}
\subsubsection{}\label{app:HHyd}
For an (ordinary) Hopf algebra $H$ and its module comodule $U$, the
left--left Yetter--Drinfeld axiom is
\begin{equation*}
  (h\leftii u)\mone \tensor(h\leftii u)\zero
  =h' u\mone \hA(h''')\tensor h''\leftii u\zero,
\end{equation*}
where $h\mapsto h'\tensor h''$ is the $H$ coproduct and $\hA$ is the
antipode of~$H$, and $h\leftii u$ and $\delta:u\mapsto u\mone\tensor
u\zero$ define the left $H$-module and left $H$-comodule structures.
The category $\HHyd$ of left--left Yetter--Drinfeld $H$-modules is
pre-braided, and braided if $\hA$ is bijective, with the braiding and
its inverse given by
\begin{align}\label{YD-braiding}
  \Psi: \cU\tensor \cV&\to \cV\tensor \cU:
  u\tensor v\mapsto u\mone\leftii v\tensor u\zero,
  \\
  \Psi^{-1}: \cV\tensor \cU&\to \cU\tensor \cV: v\tensor u\mapsto
  v\zero\tensor \hA^{-1}(v\mone)\leftii u.
\end{align}

\subsubsection{}\label{app:Radford}
For a Hopf-algebra object $\cR$ in $\HHyd$, the smash product
$\cR\Smash H$ is made into an ordinary Hopf algebra by Radford's
formula~\cite{[Radford-bos]}, dubbed bosonization when rediscovered
in~\cite{[Majid-bos]} (and actually placed into the context of braided
categories there).  The multiplication in $\cR\Smash H$ is the
standard
\begin{equation*}
  (r\Smash h)(t\Smash g) = r(h'\leftii t)\Smash h'' g,
\end{equation*}
where $h\leftii t$ is the left $H$-action on its (Yetter--Drinfeld)
modules and $h\mapsto h'\tensor h''$ is the coproduct of~$H$, and
Radford's coproduct is
\begin{equation*}
  r\Smash h
  \mapsto
  (r\bup{1}\Smash r\bup{2}{}\mone h')\tensor
  (r\bup{2}{}\zero\Smash h''),
\end{equation*}
where $r\mapsto r\mone\tensor r\zero$ is the $H$-coaction and
$r\mapsto r\bup{1}\tensor r\bup{2}\in\cR\tensor\cR$ is the coproduct
of~$\cR$.  The bialgebra is furthermore made into a Hopf algebra by
defining the antipode
\begin{equation*}
  \pmb{\A}(r\Smash h) = (1\Smash \hA(r\mone h))
  \bigl(\A(r\zero)\Smash 1\bigr),
\end{equation*}
where $\hA$ is the antipode of $H$ and $\A$ is the antipode of $\cR$.

\subsubsection{}\label{app:commcocomm}
A special case of $\HHyd$ is where $H$ is commutative and
cocommutative, $H=k\Gamma$ with a finite Abelian group $\Gamma$.  Then
Yetter--Drinfeld $H$-modules are just $\Gamma$-graded vector spaces
$\mathscr{X}=\bigoplus_{g\in\Gamma}\mathscr{X}_g$ with the left
comodule structure $\delta:x\mapsto g\tensor x$ for all $x\in
\mathscr{X}_g$, and with $\Gamma$ acting on each~$\mathscr{X}_g$.  The
action is diagonalizable, and hence
$\mathscr{X}=\bigoplus_{\chi\in\widehat\Gamma} \mathscr{X}^{\chi}$,
where $\widehat\Gamma$ is the group of characters of $\Gamma$ and
$\mathscr{X}^{\chi} = \{x\in \mathscr{X}\mid g\leftii x = \chi(g) x$
for all $g\in\Gamma\}$.  Then
\begin{equation*}
  \mathscr{X}=\bigoplus_{g\in\Gamma,\;\chi\in\widehat\Gamma} \mathscr{X}^{\chi}_g,
\end{equation*}
where $\mathscr{X}^{\chi}_g=\mathscr{X}^{\chi}\cap \mathscr{X}_g$. \
Therefore, each Yetter--Drinfeld $H$-module $\mathscr{X}$ has a basis
$(x_\alpha)$ such that, for some $g_\alpha\in\Gamma$ and
$\chi_\alpha\in\widehat\Gamma$, \ $\delta x_\alpha = g_\alpha\tensor
x_\alpha$ and $g\leftii x_\alpha = \chi_\alpha(g) x_\alpha$ for
all~$g$. \ Braiding~\eqref{YD-braiding} then becomes
$\Psi:x_\alpha\tensor x_\beta\mapsto
\chi_\beta(g_\alpha)\,x_\beta\tensor x_\alpha$.  For the Nichols
algebra generators $F_i\in X$, in particular, with
\begin{equation*}
  \delta F_i = g_i\tensor F_i\quad\text{and}\quad
  g\leftii F_i=\chi_i(g)F_i,
\end{equation*}
we recover braiding~\eqref{qij} with $q_{i,k}=\chi_k(g_i)$.

\section{\textbf{Twisted Verma modules and singular vectors of
    $\hSL2$}}\label{app:sl2}
\subsection{The $\hSL2$ algebra}\label{app:sl2-conv}
Our conventions for the $\hSL2_k$ algebra are
\begin{align*}
  [\Jnaught_m, J^{\pm}_n] &= \pm J^{\pm}_{m + n},
  ,\\
  [\Jnaught_m, \Jnaught_n] &= 
  \ffrac{k}{2} m \delta_{m + n, 0},\\
  [\Jplus_m, \Jminus_n] &= 
  k m \delta_{m + n, 0} + 2 \Jnaught_{m + n},
\end{align*}
where $J^{\pm,0}(z)=\sum_{n\in\oZ}J_n z^{-n-1}$.  The zero modes
$J^{\pm,0}_0$ generate an $s\ell(2)$ Lie algebra, which we call the
\textit{zero-mode} $s\ell(2)$ to distinguish it from other classical
and quantum $s\ell(2)$ algebras.

The Sugawara energy--momentum tensor constructed from the $\hSL2_k$
currents is
\begin{equation}\label{eq:Sug}
  T_{\text{Sug}}(z) \eqdef\ffrac{1}{2 (k + 2)}\bigl( \Jplus\Jminus(z)
  + \Jminus\Jplus(z) + \Jnaught\Jnaught(z)\bigr).
\end{equation}
Several different energy--momentum tensors occur in the text, and we
refer to dimensions of (primary) fields determined by the OPE with
$T_{\text{Sug}}(z)$ as the Sugawara dimensions.

\subsection{Twisted verma modules of $\hSL2$}
\label{app:twisted-verma}
We fix our conventions regarding twisted Verma
modules~\cite{[FSST]}. For $\lambda\in\oC$ and
$\vartheta\in\oZ$, the twisted Verma module
$\rep{M}_{\lambda;\vartheta}$ is freely generated by
$\Jplus_{\leq\vartheta-1}$, $\Jminus_{\leq-\vartheta}$, and
$\Jnaught_{\leq-1}$ from a twisted highest-weight vector
$\ket{\lambda;\vartheta}$ defined by the conditions
\begin{equation}\label{twistedhw}
  \begin{split}
    \Jplus_{\geq\vartheta}\,\ket{\lambda;\vartheta}&=
    \Jnaught_{\geq1}\,\ket{\lambda;\vartheta}=
    \Jminus_{\geq-\vartheta+1}\,\ket{\lambda;\vartheta}=0,\\
    \bigl(\Jnaught_{0}+\ffrac{k}{2}\,\vartheta\bigr)\,\ket{\lambda;\vartheta}&=
    \lambda\,\ket{\lambda;\vartheta}.
  \end{split}
\end{equation}
Setting $\vartheta=0$ gives the usual (``untwisted'') Verma modules.  We
write $\ket{\lambda}=\ket{\lambda;0}$ and, similarly,
$\hVerma_{\lambda}=\rep{M}_{\lambda;0}$.  The highest-weight vector
$\ket{\lambda;\vartheta}$ of a twisted Verma module has the Sugawara
dimension
\begin{equation}\label{dim-Sug-tw}
  \Delta_{\lambda;\vartheta}=
  \ffrac{\lambda^2+\lambda}{{k\!+\!2}}-\vartheta \lambda +
  \ffrac{k}{4}\vartheta^2.
\end{equation}

We write $\ket{\alpha}\doteq\ket{\lambda;\vartheta}$ whenever
conditions~\eqref{twistedhw} are satisfied for a state~$\ket{\alpha}$.

We introduce the $(\mbox{charge},\mbox{dimension})$ bigrade for vectors
in a Verma module $\hVerma_{\lambda}$ in an obvious way, by assigning the
grade $(\lambda,\Delta_{\lambda})$ to the highest-weight vector
$\ket{\lambda}$ and setting $\gr\jp_{-n}=(1,n)$,
$\gr\jm_{-n}=(-1,n)$, and $\gr\Jnaught_{-n}=(0,n)$.  Then, e.g.,
$\gr\Jplus_{-1}\ket{\lambda}=(\lambda+1,\Delta_{\lambda}+1)$.  For the
twisted highest-weight state $\ket{\lambda;\vartheta}$, we have
\begin{equation*}
  \gr \ket{\lambda;\vartheta} =
  (\lambda-\ffrac{k}{2}\,\vartheta, \Delta_{\lambda;\vartheta}).
\end{equation*}

Twists, although producing nonequivalent modules, do not alter the
submodule structure, and we can therefore reformulate a classic result
as follows.

\begin{Thm}[\cite{[KK],[MFF]}]\label{MFFthm}
  \addcontentsline{toc}{subsection}{\thesubsection. \ \ Singular
    vectors} A singular vector exists in the twisted Verma module
  $\hVerma_{\lambda;\vartheta}$ of $\hSL2_k$ if and only if $\lambda$ can
  be written as $\lambda=\jp(r,s)$ or $\lambda=\jm(r,s)$ with $r,s
  \in\oN$, where
\begin{equation*}
    \jp(r,s)=\ffrac{r\!-\!1}{2}-(k\!+\!2)\ffrac{s-1}{2},\qquad
    \jm(r,s)=-\ffrac{r\!+\!1}{2}+(k\!+\!2)\ffrac{s}{2}.
  \end{equation*}
  Whenever $\lambda=\jp(r,s)$ or $\lambda=\jm(r,s)$, the corresponding
  singular vector is given by
  \begin{multline}\label{mffplus}
    \MFFplus{r,s;\vartheta}=
    (\Jminus_{-\vartheta})^{r+(s-1)(k+2)}
    (\Jplus_{\vartheta-1})^{r+(s-2)(k+2)}(\Jminus_{-\vartheta})^{r+(s-3)(k+2)}
    \ldots\\
    \cdot(\Jplus_{\vartheta-1})^{r-(s-2)(k+2)}
    (\Jminus_{-\vartheta})^{r-(s-1)(k+2)}\ket{\jp(r,s);\vartheta}
  \end{multline}
  or
  \begin{multline}\label{mffminus}
    \MFFminus{r,s;\vartheta}=
    (\Jplus_{\vartheta-1})^{r+(s-1)(k+2)}(\Jminus_{-\vartheta})^{r+(s-2)(k+2)}
    (\Jplus_{\vartheta-1})^{r+(s-3)(k+2)}\ldots\\
    {}\cdot(\Jminus_{-\vartheta})^{r-(s-2)(k+2)}
    (\Jplus_{\vartheta-1})^{r-(s-1)(k+2)}\ket{\jm(r,s);\vartheta}.
  \end{multline}
\end{Thm}
The dependence of $\lambda^{\pm}(r,s)$ on $k$ is not indicated for the
brevity of notation.

\subsubsection{}
We recall that the above formulas yield polynomial expressions in the
currents via repeated application of the formulas
\begin{align}
  (\Jminus_0)^\alpha\,\Jplus_m &= -\alpha
  (\alpha - 1) \Jminus_m(\Jminus_0)^{\alpha-2} - 2 \alpha
  \Jnaught_m\,(\Jminus_0)^{\alpha-1} + \Jplus_m\,(\Jminus_0)^{\alpha}
  ,\notag
  \\
  (\Jminus_0)^{\alpha}\,\Jnaught_m &= \alpha
  \Jminus_m(\Jminus_0)^{\alpha-1} + \Jnaught_m\,(\Jminus_0)^{\alpha}
  ,\notag\\
  (\Jplus_{-1})^\alpha\,\Jminus_m &= -\alpha (\alpha - 1)
  \Jplus_{m-2}(\Jplus_{-1})^{\alpha-2} - k\,\alpha\,\delta_{m - 1, 0}
  (\Jplus_{-1})^{\alpha-1}
  \label{properties}
  \\*
  &\qquad{}+ 2 \alpha
  \Jnaught_{m-1}\,(\Jplus_{-1})^{\alpha-1} +
  \Jminus_m\, (\Jplus_{-1})^{\alpha},\notag
  \\
  (\Jplus_{-1})^{\alpha}\,\Jnaught_m &= -\alpha
  \Jplus_{m-1}(\Jplus_{-1})^{\alpha-1} +
  \Jnaught_m\,(\Jplus_{-1})^{\alpha},\notag
\end{align}
and, when necessary, of their images under the Lie algebra
homomorphism $\Jplus_n\mapsto\Jplus_{n+\vartheta}$,
$\Jnaught_n\mapsto\Jnaught_{n}+\frac{k}{2}\vartheta\delta_{n,0}$,
$\Jminus_n\mapsto\Jminus_{n-\vartheta}$.  Formulas~\eqref{properties} are
derived for positive integer $\alpha $ and are then continued to
arbitrary complex~$\alpha $.

\subsubsection{}
Singular vectors
$\MFFname^{\pm}(r,s;\vartheta)=\MFFname^{\pm}(r,s;\vartheta|\lambda)$
constructed on a twisted highest-weight state $\ket{\lambda;\vartheta}$
lie in the grades
\begin{align*}
  \gr\MFFplus{r,s;\vartheta|\lambda}&=
  (\lambda-r-\ffrac{k}{2}\,\vartheta,
  \Delta_{\lambda;\vartheta} + r(s-1+\vartheta)),
  \\
  \gr\MFFminus{r,s;\vartheta|\lambda}&=
  (\lambda+r-\ffrac{k}{2}\,\vartheta,
  \Delta_{\lambda;\vartheta} + r(s-\vartheta)).
\end{align*}

\subsubsection{}For $s=1$, the above singular vectors do not require
any algebraic rearrangements and take the simple form
\begin{equation*}
  \MFFplus{r,1;\vartheta|\lambda}=(\Jminus_{-\vartheta})^{r}\ket{\lambda;\vartheta},
  \qquad 
  \MFFminus{r,1;\vartheta|\lambda}=(\Jplus_{\vartheta-1})^{r}\ket{\lambda;\vartheta}.
\end{equation*}

\section{\textbf{Triplet--triplet algebra: examples}}
\subsection{Asymmetric realization}\label{app:example-a}
To give examples of $\WWtritri$ algebras, we first recall that the
$\W^+(z)$ generator is always given by
\begin{equation*}
    \W^+(z) = e^{p\Aii(z) +  \frac{p}{1 - 2 p}\Athree(z) + \frac{p}{1 - 2 p}\Afour(z)}
\end{equation*}
and the entire multiplet in the left part of Fig.~\ref{fig:a-TRIPLET}
is
\begin{align*}
  \Jminus_0\W^+(z)&=(\partial\Aii(z) + \partial\Afour(z)) e^{p\Aii(z)},
  \\
  (\Jminus_0)^2\W^+(z)&=
  \bigl(
  \fffrac{1}{2} \partial\Aii \partial\Aii(z)
  \!-\! \partial^2\Aii(z)
  \!+\! \partial\Aii\partial\Afour(z)
  \!+\! \fffrac{1}{2} \partial\Afour\partial\Afour(z)
  \!-\! \partial^2\Afour(z)    
  \!\bigr)
  e^{p\Aii(z) + \frac{p}{2 p - 1}(\Athree(z) + \Afour(z))}.
\end{align*}

\subsubsection{}For $p=2$, the other two fields in~\eqref{a:the-W-j=0}
are given by (omitting the conventional~$(z)$ argument of gradients of
the scalars)
\begin{align*}
  \W^0(z) &=
  \bigl(-\fffrac{1}{3} \partial^3\Aii  + 
  2 \partial^2\Aii \partial\Aii  - 
  \fffrac{4}{3} (\partial\Aii)^3\bigr)
  e^{-\frac{2}{3}\Athree(z) - \frac{2}{3}\Afour(z)},
  \\
  \W^-(z) &=
  \bigl(8 \partial^2\Aii 
    + 8 (\partial\Aii)^2\bigr)
    e^{-2\Aii(z) -\frac{2}{3}\Athree(z) -\frac{2}{3}\Afour(z)}.
\end{align*}
Also for $p=2$, the middle terms in the two zero-mode triplets are
given by
\begin{align*}  
  \Jminus_0 \W^{0}(z) &=
  (\partial^2\Aii)^2 + 
  \partial^3\Aii \partial\Aii
  - 2 \partial^2\Aii (\partial\Aii)^2
  +\bigl[-\fffrac{1}{3} \partial^3\Aii
  + 2 \partial^2\Aii \partial\Aii - 
  \fffrac{4}{3} (\partial\Aii)^3\bigr] \partial\Afour
  - \fffrac{1}{6} \partial^4\Aii
  \\
  \intertext{and}
  \Jminus_0 \W^{-}(z) &=
  \bigl(4 \partial^3\Aii  + 
  8 \partial^2\Aii \partial\Afour  - 
  8 (\partial\Aii)^3 + 
  8 (\partial\Aii)^2 \partial\Afour\bigr) e^{-2\Aii(z)},
\end{align*}
and the leftmost terms, by somewhat lengthier expressions.

\subsubsection{}
For $p=3$, the $\W^0(z)$ and $\W^-(z)$ fields in~\eqref{a:the-W-j=0}
are
\begin{align*}
  \W^0(z) &=\Bigl(-\fffrac{1}{40} \partial^5\Aii  + 
  \fffrac{3}{4} \partial^3\Aii \partial^2\Aii 
  + 
  \fffrac{3}{8} \partial^4\Aii \partial\Aii 
  - 
  \fffrac{27}{8} (\partial^2\Aii)^2 \partial\Aii 
  - 
  \fffrac{9}{4} \partial^3\Aii (\partial\Aii)^2
  \\ &\qquad{}
  + 
  \fffrac{27}{4} \partial^2\Aii (\partial\Aii)^3 
  - 
  \fffrac{81}{40} (\partial\Aii)^5\Bigr)
  e^{-\frac{3}{5}\Athree(z) -\frac{3}{5}\Afour(z)},
  \\
  \W^-(z) &=\Bigl(\!\fffrac{3}{2} \partial^4\Aii 
  - 
  \fffrac{117}{2} (\partial^2\Aii)^2
  + 
  27 \partial^3\Aii \partial\Aii 
  - 
  54 \partial^2\Aii (\partial\Aii)^2 
  - 
  \fffrac{81}{2} (\partial\Aii)^4\!\Bigr)
  e^{-3\Aii(z) -\frac{3}{5}\Athree(z) -\frac{3}{5}\Afour(z)}.
\end{align*}

\subsection{Symmetric realization}\label{app:example-s}
For the ``symmetric'' realization in Sec.~\ref{CFTsymm}, the expanded
expressions for the algebra generators are much more bulky than for
the ``asymmetric'' realization, and we therefore restrict ourself to
only $p=2$ for illustration.

The field $\W^{-}(z)$ (see \eqref{s:wminus-x3}) is
\begin{align*}
  \W^{-}(z) &= \bigl(-\partial^4\Aone 
  + 
  42 (\partial^2\Aone)^2
  + 
  24 \partial^2\Aone \partial^2\Atwo 
  - 
  24 \partial^3\Aone \partial\Aone
  - 
  12 \partial^3\Aone \partial\Atwo
  \\
  &\quad{}+ 
  24 (\partial\Aone)^2\partial^2\Atwo 
  + 
  24 \partial^2\Aone (\partial\Aone)^2
  + 
  24 \partial^2\Aone \partial\Aone\partial\Atwo
  + 
  24 \partial^2\Aone (\partial\Atwo)^2
  \\
  &\quad{}+ 
  24 (\partial\Aone)^4
  + 
  48 (\partial\Aone)^3\partial\Atwo
  + 
  24 (\partial\Aone)^2(\partial\Atwo)^2
  \bigr)
  e^{-\frac{8}{3}\Aone(z) -\frac{4}{3}\Atwo(z)
    -\frac{2}{3}\Athree(z)}.
  \\
  \intertext{It looks simpler in the Wakimoto realization
    in~\bref{s:Wak}:}
  \W^{-}(z) &=
  \bigl(18 \partial\beta \partial\beta 
  - 12 \partial^2\beta \beta  - 
  24 \partial\Azero \partial\beta \beta  + 
  24 \partial^2\Azero \beta \beta  + 
  24 \partial\Azero \partial\Azero \beta \beta\bigr) e^{-2\Azero(z)}.
\end{align*}

Also, the middle element of the ``middle'' triplet is explicitly given
by
\begin{align*}
  \Jminus_0\W^{0}(z) &=
  8 (\partial\Aone)^4 + 
  16 (\partial\Aone)^3\partial\Atwo - 
  16 \partial\Aone(\partial\Atwo)^3 + 
  24 \partial\Aone\partial^2\Atwo \partial\Atwo
  - 4 \partial\Aone\partial^3\Atwo
  \\ &\quad{}
  - 8 (\partial\Atwo)^4 - 
  24 \partial^2\Aone (\partial\Aone)^2 - 
  24 \partial^2\Aone \partial\Aone\partial\Atwo + 
  6 (\partial^2\Aone)^2
  + 24 \partial^2\Atwo(\partial\Atwo)^2
  \\ &\quad{}
  - 6 (\partial^2\Atwo)^2 + 
  8 \partial^3\Aone \partial\Aone + 4 \partial^3\Aone \partial\Atwo - 
  8 \partial^3\Atwo \partial\Atwo - \partial^4\Aone + 
  \partial^4\Atwo.
\end{align*}

\parindent0pt

\end{document}